\documentclass{amsart}

\usepackage{amscd,amsmath,xypic,amssymb,combelow,tikz-cd,etoolbox,calligra,mathrsfs,enumitem,hyperref,mathtools,epigraph}

\usepackage[russian,english]{babel}

\makeatletter
\patchcmd{\@settitle}{\uppercasenonmath\@title}{}{}{}
\makeatother

\xyoption{all}
\CompileMatrices

\makeatletter
\@addtoreset{equation}{section}
\makeatother

\newtheorem{theorem}[subsection]{Theorem}

\newtheorem{proposition}[subsection]{Proposition}
\newtheorem{lemma}[subsection]{Lemma}
\newtheorem{corollary}[subsection]{Corollary}
\newtheorem{conjecture}[subsection]{Conjecture}

\newtheorem{definition}[subsection]{Definition}

\newtheorem{claim}[subsection]{Claim}

\newtheorem{remark}[subsection]{Remark}

\def\loccitt{\emph{loc. cit.}}
\def\loccit{\emph{loc. cit. }}

\def\fg{{\mathfrak{g}}}

\def\fY{{\mathfrak{Y}}}
\def\fZ{{\mathfrak{Z}}}

\def\BA{{\mathbb{A}}}
\def\BC{{\mathbb{C}}}

\def\BN{{\mathbb{N}}}
\def\BF{{\mathbb{F}}}
\def\BP{{\mathbb{P}}}
\def\BR{{\mathbb{R}}}
\def\BQ{{\mathbb{Q}}}
\def\BZ{{\mathbb{Z}}}

\def\woo{\widehat{\otimes}}
\def\too{\widetilde{\otimes}}

\def\CA{{\mathcal{A}}}
\def\CB{{\mathcal{B}}}

\def\CE{{\mathcal{E}}}

\def\CI{{\mathcal{I}}}

\def\CK{{\mathcal{K}}}
\def\CL{{\mathcal{L}}}
\def\CM{{\mathcal{M}}}
\def\CN{{\mathcal{N}}}
\def\CO{{\mathcal{O}}}

\def\CR{{\mathcal{R}}}
\def\CS{{\mathcal{S}}}

\def\CU{{\mathcal{U}}}
\def\CV{{\mathcal{V}}}
\def\CW{{\mathcal{W}}}

\def\Hom{\textrm{Hom}}

\def\tCV{\widetilde{\CV}}

\def\tv{\widetilde{v}}

\def\tV{\widetilde{V}}
\def\tW{\widetilde{W}}

\def\tX{\widetilde{X}}
\def\tY{\widetilde{Y}}
\def\tA{\widetilde{A}}
\def\tB{\widetilde{B}}

\def\vs{\varsigma}

\def\and{\textrm{ }\&\textrm{ }}

\def\sym{\textrm{sym}}
\def\Sym{\textrm{Sym}}

\def\tB{\widetilde{B}}

\def\nn{{{\BN}}^I}
\def\zz{{{\BZ}}^I}

\def\UU{{\mathscr{U}_q(\widehat{\fg}_Q)}}

\def\bCN{\bar{\CN}}

\def\bk{\boldsymbol{k}}

\def\bm{\boldsymbol{m}}
\def\bn{\boldsymbol{n}}
\def\bp{\boldsymbol{p}}
\def\bv{\boldsymbol{v}}
\def\bw{\boldsymbol{w}}
\def\bx{\boldsymbol{x}}
\def\by{\boldsymbol{y}}
\def\bz{\boldsymbol{z}}

\def\bS{\boldsymbol{S}}
\def\bT{\boldsymbol{T}}

\def\bs{{\boldsymbol{\vs}}}
\def\bth{{\boldsymbol{\theta}}}
\def\b0{{\boldsymbol{0}}}

\def\op{\text{op}}

\def\oii{\overrightarrow{ii}}
\def\oij{\overrightarrow{ij}}
\def\oji{\overrightarrow{ji}}
\def\loc{\text{loc}}

\def\stab{\text{Stab}}

\def\attr{\text{Attr}}
\def\eattr{\emph{Attr}}

\def\Tan{\text{Tan}}
\def\eTan{\emph{Tan}}
\def\Nor{\text{Nor}}
\def\eNor{\emph{Nor}}

\def\mindeg{\text{min deg}}
\def\maxdeg{\text{max deg}}

\def\sdet{\text{sdet}}

\def\Gr{\text{Gr}}

\begin{document}

\title[Quantum loop groups and $K$-theoretic stable envelopes]{\Large{\textbf{Quantum loop groups and $K$-theoretic stable envelopes}}}

\author[Andrei Negu\cb t]{Andrei Negu\cb t}
 
\address{École Polytechnique Fédérale de Lausanne (EPFL), Lausanne, Switzerland \newline \text{ } \ \ Simion Stoilow Institute of Mathematics (IMAR), Bucharest, Romania} 

\email{andrei.negut@gmail.com}

\maketitle

\begin{abstract} We develop the connection between the preprojective $K$-theoretic Hall algebra (\cite{SV Hilb, YZ}) of a quiver $Q$ and the quantum loop group associated to $Q$ via stable envelopes of Nakajima quiver varieties (\cite{OS}). \\


\end{abstract}

$$$$

\section{Introduction}

\medskip

\subsection{} Let us consider a quiver $Q$ with vertex set $I$ and edge set $E$; loops and multiple edges are allowed. In the present paper, we make the convention that the set $\BN$ contains 0. Given $\bv,\bw \in \nn$, one may consider the Nakajima quiver variety
\begin{equation}
\label{eqn:nakajima intro}
\CM_{\bv,\bw}
\end{equation}
parameterizing $\bv$-dimensional representations of the preprojective algebra of $Q$, with framing given by $\bw$ (see Subsection \ref{sub:nakajima} for details). We will consider the (localized, equivariant) algebraic $K$-theory groups of Nakajima quiver varieties, and assemble them into
\begin{equation}
\label{eqn:k-theory intro}
K(\bw) = \bigoplus_{\bv \in \nn} K(\bw)_{\bv}, \qquad \text{where} \qquad K(\bw)_{\bv} = K_{T \times G_{\bw}}(\CM_{\bv,\bw})_{\text{loc}}
\end{equation}
(see Subsection \ref{sub:torus} for details on the $T \times G_{\bw}$ action on $\CM_{\bv,\bw}$). A variant of the $\bw = 0$ case of \eqref{eqn:nakajima intro} yields certain stacks $\fY_{\bv}$; as we will recall in Subsection \ref{sub:k-ha}, the (localized, equivariant) algebraic $K$-theory groups of these stacks can be assembled into the Schiffmann-Vasserot preprojective $K$-theoretic Hall algebra
\begin{equation}
\label{eqn:k-ha intro}
\CA^+ = \bigoplus_{\bv \in \nn} \CA_{\bv} \qquad \text{where} \qquad \CA_{\bv} =  K_T(\fY_{\bv})_{\text{loc}}
\end{equation}
defined over the field 
$$
\BF =  \text{Frac}(\text{Rep}_T) = \BQ(q^{\frac 12}, t_e)_{e \in E}
$$
The double of the $K$-theoretic Hall algebra is defined as the $\BF$-vector space
\begin{equation}
\label{eqn:double intro}
\CA = \CA^+ \otimes \CA^0 \otimes \CA^-
\end{equation}
where $\CA^- = \CA^{+,\op}$ (whose graded direct summands will be denoted by $\CA_{-\bv}$, in order to differentiate them from the direct summands which appear in \eqref{eqn:k-ha intro}), and 
$$
\CA^0 = \BF \left[a_{i,\pm d}, b_{i,\pm d}, q^{\pm \frac {v_i}2}, q^{\pm \frac {w_i}2} \right]_{i \in I, d \geq 1}
$$
The operators $a_{i,\pm d}, b_{i,\pm d}$ (respectively $q^{\pm \frac {v_i}2}, q^{\pm \frac {w_i}2}$) act on the $K$-theory groups of Nakajima quiver varieties as multiplication by Adams operations in the tautological bundles $\CV_i,W_i$ (respectively the $q^{\frac 12}$-th powers of the ranks of these bundles), see Subsection \ref{sub:power sum} for details. The vector space \eqref{eqn:double intro} is made into a $\BF$-algebra by imposing the commutation relations \eqref{eqn:rel double 0}--\eqref{eqn:rel double 5} between the three tensor factors. Inspired by the main construction of \cite{Nak}, there is an action (versions of which were studied e.g. in \cite{SV Hilb, YZ})
\begin{equation}
\label{eqn:action 1 intro}
\CA \curvearrowright K(\bw)
\end{equation}
for all $\bw \in \nn$, which we will review in Subsection \ref{sub:action} and prove in Theorem \ref{thm:action}.

\medskip

\begin{remark}
\label{rem:intro}

The only essential feature of localized $K$-theoretic Hall algebras that we will use is that they are generated by the direct summands $\bv = (\dots,0,1,0,\dots)$ of \eqref{eqn:k-ha intro}. This was proved in \cite{Wheel}, and it also holds when the torus $T$ is replaced by any subtorus $H \subset T$ that satisfies the assumption of Remark \ref{rem:smaller torus}, namely the fact that $H$ contains a rank one subtorus which scales all the edge maps of Nakajima quiver varieties with strictly positive weight. Under this assumption, all the results of the present paper hold with $T$-equivariance replaced by $H$-equivariance.

\end{remark}

\medskip

\subsection{}

The $K$-theoretic stable envelope construction (which we will review in Section \ref{sec:stable}) yields homomorphisms
\begin{equation}
\label{eqn:stable basis intro}
K(\bw) \otimes K(\bw') \xrightarrow{\stab, \stab'} K(\bw+\bw')
\end{equation}
for any $\bw,\bw' \in \nn$. These maps were first constructed in cohomology in \cite{MO}, then they were axiomatically defined and studied in $K$-theory in \cite{OS}, then constructed in elliptic cohomology in \cite{AO}, and finally constructed in a broader generality than Nakajima quiver varieties in \cite{O1,O2}. The numerous properties enjoyed by the stable envelopes \eqref{eqn:stable basis intro} imply that the compositions
\begin{equation}
\label{eqn:r matrix intro}
R_{\bw,\bw'} : K(\bw) \otimes K(\bw') \xrightarrow{\stab} K(\bw+\bw') \xrightarrow{{\stab'}^T} K(\bw) \otimes K(\bw')
\end{equation}
(as $\bw,\bw'$ vary over $\nn$) satisfy the quantum Yang-Baxter equation. Following the Faddeev-Reshetikhin-Takhtajan (FRT) principle laid out in \cite{MO}, one defines the quantum group
\begin{equation}
\label{eqn:def u intro}
\UU \subset \prod_{\bw \in \nn} \text{End}(K(\bw))
\end{equation}
as the algebra generated by all possible matrix coefficients of \eqref{eqn:r matrix intro} with respect to the tensor factor $K(\bw')$, as $\bw'$ ranges over $\nn$, together with central elements that correspond to multiplication by the Schur functors of the framing vector spaces. Tautologically, the quantum group comes endowed with an action for all $\bw \in \nn$
\begin{equation}
\label{eqn:action 2 intro}
\UU \curvearrowright K(\bw)
\end{equation}
The main purpose of the present paper is to compare the actions \eqref{eqn:action 1 intro} and \eqref{eqn:action 2 intro}. This will be accomplished in Theorem \ref{thm:main}, which we will prove conditionally on a certain coproduct computation, that we will now explain. 

\medskip

\subsection{} As we will review in Subsection \ref{sub:frt principle}, the FRT formalism implies that $\UU$ is a Hopf algebra. In particular, the coproduct \footnote{See Definition \ref{def:completion} for the meaning of the completion $\too$.}
\begin{equation}
\label{eqn:cop intro}
\Delta : \UU \rightarrow \UU \too \UU
\end{equation}
is defined for all $x \in \UU$ by the commutativity of the following diagram
$$
\xymatrixcolsep{5pc}\xymatrix{
K(\bw) \otimes K(\bw') \ar[d]_{\stab} \ar[r]^{\Delta(x)} & K(\bw) \otimes K(\bw')  \ar[d]^{\stab} \\
K(\bw+\bw') \ar[r]^{x} & K(\bw+\bw')}
$$
for all $\bw,\bw' \in \nn$. As we will recall in Theorem \ref{thm:generate} (see also Remark \ref{rem:intro}), the algebra $\CA^+$ is generated by symbols $\{e_{i,d}\}_{i \in I, d\in \BZ}$ corresponding to elements in the $\bv = (\dots,0,1,0,\dots)$ direct summands of \eqref{eqn:k-ha intro}. The analogous generators of $\CA^-$ will be denoted by $\{f_{i,d}\}_{i \in I, d \in \BZ}$. As a consequence of this fact, we will prove in Proposition \ref{prop:generate} that the algebra $\CA$ is generated by the symbols
\begin{equation}
\label{eqn:generators intro}
\Big\{ e_{i,0}, f_{i,0}, a_{i,\pm 1}, b_{i,\pm d}, q^{\pm \frac {v_i}2}, q^{\pm \frac {w_i}2} \Big\}_{i \in I, d \geq 1}
\end{equation}
(see \eqref{eqn:action e}--\eqref{eqn:action v and w} for formulas for the action of the aforementioned generators on $K(\bw)$). Thus, one may ask for explicit formulas for the coproduct \eqref{eqn:cop intro} on the generators of the algebra $\CA$. To this end, we will prove in Proposition \ref{prop:coproduct e and f} that
\begin{equation}
\label{eqn:coproduct e intro} 
\Delta(e_{i,0}) = h_{i,0}^{-1} \otimes e_{i,0} + e_{i,0} \otimes 1 
\end{equation}
\begin{equation}
\label{eqn:coproduct f intro} 
\Delta(f_{i,0}) = f_{i,0} \otimes h_{i,0} + 1 \otimes f_{i,0} 
\end{equation}
for all $i \in I$, where we let
\begin{equation}
\label{eqn:g to h}
h_{i,0} = q^{\frac 12 (w_i - 2v_i + \sum_{e= \oij} v_j + \sum_{e= \oji} v_j)}
\end{equation}
Moreover, $\Delta(b_{i,\pm d}) = b_{i,\pm d} \otimes 1 + 1 \otimes b_{i,\pm d}$, $\Delta(q^{\frac {v_i}2}) = q^{\frac {v_i}2} \otimes q^{\frac {v_i}2}$, $\Delta (q^{\frac {w_i}2} ) = q^{\frac {w_i}2} \otimes q^{\frac {w_i}2}$.

\medskip

\begin{conjecture}
\label{conj:coproduct h intro}

For any sign $\pm \in \{+,-\}$ and any $i \in I$, we have
\begin{equation}
\label{eqn:coproduct h pm intro}
\Delta(a_{i,\pm 1}) \in a_{i, \pm 1} \otimes 1 + 1 \otimes a_{i ,\pm 1} + \CA^{<} \too \CA^{>}
\end{equation}
where 
\begin{equation}
\label{eqn:a big and small}
\CA^{>} = \bigoplus_{\bn \in \nn \backslash \b0} \CA_{\bn} \otimes \CA^0 \quad \text{and} \quad \CA^{<} = \bigoplus_{\bn \in \nn \backslash \b0} \CA^0 \otimes \CA_{-\bn}
\end{equation}
(above, $\b0 = (0,\dots,0) \in \nn$). 

\end{conjecture}

\medskip

\noindent The coproduct formulas \eqref{eqn:coproduct h pm intro} are explicitly known only in the case when $Q$ is a Dynkin diagram of ADE type (in the cohomological case, they were connected to stable envelopes in \cite{Mc}) and conjecturally when $Q$ is the Jordan quiver.

\medskip

\begin{remark}
	
When $Q$ is of ADE type, $\CA$ should be isomorphic to the quantum affine group $U_q(\widehat{\fg}_Q)|_{c=1}$ as a bialgebra, where the latter is endowed with the Drinfeld-Jimbo coproduct. Under this isomorphism, the elements $a_{i,\pm d}$ should correspond to the root vectors for $\pm d$ times the imaginary root. As the coproduct of the latter is very complicated (see \cite[Section 7]{Da}) we expect the same to be true of $\Delta(a_{i,\pm d})$. 
	
\end{remark}

\medskip

\subsection{} The following is our main result.

\medskip

\begin{theorem}
	\label{thm:main}
	
	If Conjecture \ref{conj:coproduct h intro} holds, there exists an isomorphism 
	\begin{equation}
		\label{eqn:iso main}
		\CA \cong \UU
	\end{equation}
	which intertwines the actions \eqref{eqn:action 1 intro} and \eqref{eqn:action 2 intro} for all $\bw \in \nn$.
	
\end{theorem}

\medskip

\noindent The proof of the Theorem above will be given in the following five steps. 

\medskip

\noindent \textbf{Step 1:} The map
\begin{equation}
\label{eqn:iota intro}
\CA \rightarrow \prod_{\bw \in \nn} \text{End}(K(\bw))
\end{equation}
induced by \eqref{eqn:action 1 intro} is injective (see Proposition \ref{prop:inj}). 

\medskip

\noindent \textbf{Step 2:} Any vector and covector of $K(\bw)$ can be written as a linear combination (with coefficients in $\BF_{\bw} = \BF[u_{i1}^{\pm 1},\dots, u_{iw_i}^{\pm 1}]^{\sym}_{i \in I}$) of vectors and covectors of the form
$$
K(\bw)_{\b0} \xrightarrow{e} K(\bw)_{\bm} \qquad \text{and} \qquad K(\bw)_{\bn} \xrightarrow{f} K(\bw)_{\b0}
$$
respectively, for various $e \in \CA_{\bm} \otimes \CA_0$ and $f \in\CA_0 \otimes \CA_{-\bn}$ (see Proposition \ref{prop:any coefficient}).

\medskip

\noindent \textbf{Step 3:} The image of the map \eqref{eqn:iota intro} is contained in the subalgebra $\UU$ of \eqref{eqn:def u intro}, thus yielding an injective algebra homomorphism
\begin{equation}
\label{eqn:map intro}
\CA \hookrightarrow \UU
\end{equation}
We will prove the aforementioned claim by exhibiting the generators \eqref{eqn:generators intro} of $\CA$ as certain explicit matrix coefficients of the $R$-matrix \eqref{eqn:r matrix intro} (see Proposition \ref{prop:maps to subalgebra}). 

\medskip

\noindent \textbf{Step 4:} The coproduct \eqref{eqn:cop intro} on $\UU$ preserves the subalgebra $\CA$ of \eqref{eqn:map intro}
\begin{equation}
\label{eqn:delta preserves a intro}
\Delta(\CA) \subset \CA \too \CA
\end{equation}
This statement uses Proposition \ref{prop:generate}, formulas \eqref{eqn:coproduct e intro}--\eqref{eqn:coproduct f intro} and Conjecture \ref{conj:coproduct h intro}.

\medskip

\noindent \textbf{Step 5:} Using the general FRT formalism, we will deduce from \eqref{eqn:delta preserves a intro} and Step 2 the fact that arbitrary matrix coefficients of \eqref{eqn:r matrix intro} lie in the algebra $\CA$, thus establishing the surjectivity of the map \eqref{eqn:map intro}. 

\medskip

\subsection{} The isomorphism $\CA \cong \UU$ is one of those mathematical statements that most specialists in the field believe, but no proof has yet been written down. For example, conjectural versions of this statement have appeared in \cite[Conjecture 1.2]{Pad 2} and in \cite{VV}. There has been more work on the subject in cohomology, starting with \cite{Dav} and \cite{SV Gen, SV Coh}, and culminating with proofs of the non-localized cohomological version of Theorem \ref{thm:main} in \cite{bd, SV new}. However, none of the above methods involve the program outlined in Steps 1-5 above, and we believe that proving the cohomological analogue of Conjecture \ref{conj:coproduct h intro} would be of independent interest.

\medskip

\noindent Most of our results require the algebras $\CA$ and $\UU$ to be localized, i.e. defined over a field (which can be the field of rational functions in generic equivariant parameters, or more generally any specialization of said equivariant parameters that satisfies Remark \ref{rem:smaller torus}). A more interesting and (in the author's opinion) much more difficult question would be to establish the isomorphism \eqref{eqn:iso main} for appropriately defined integral forms of the algebras in question. 

\medskip

\subsection{}

I would like to thank Michael McBreen and Andrei Okounkov for teaching me almost everything I know about stable envelopes, and for their help and advice on this project. I also thank the anonymous referee for numerous helpful comments. I gratefully acknowledge NSF grant DMS-$1845034$, as well as support from the MIT Research Support Committee.

\medskip

\subsection{}
\label{sub:notation 1}

Let us briefly review the main notation used in the present paper, particularly that which pertains to $K$-theory. Given any vector spaces (or more generally locally free sheaves, equivariant or not) $V$ and $W$, we will write
\begin{equation}
	\label{eqn:k-th}
V - W
\end{equation}
for their formal difference. As vector spaces (or more generally locally free sheaves, equivariant or not) have a notion of tensor product, it will be convenient to write
\begin{equation}
\label{eqn:conv 1}
\frac {V'}V \quad \text{instead of} \quad V' \otimes V^\vee
\end{equation}
As a $K$-theory class, the formula above is additive in $V$ and $V'$
$$
\frac {V' \pm' W'}{V \pm W} = \frac {V'}V \pm \frac {V'}W \pm' \frac {W'}V \pm \pm' \frac {W'}W
$$
However, we may not simplify terms from \eqref{eqn:conv 1} multiplicatively, i.e. $\frac {V' \otimes W}{V \otimes W} \neq \frac {V'}V$. The quantity 
$$
\wedge^\bullet \left( \frac Vz \right) = \sum_{i = 0}^{\text{dim }V } \frac {\wedge^i(V)}{(-z)^i}
$$
is additive in $V$, so we may extend it to formal differences \eqref{eqn:k-th}
\begin{equation}
\label{eqn:notation exterior}
\wedge^\bullet\left(\frac {V-W}z\right) =  \frac {\displaystyle\wedge^\bullet\left(\frac Vz\right)}{\displaystyle\wedge^\bullet\left(\frac Wz\right)}
\end{equation}
We interpret \eqref{eqn:notation exterior} as a rational function (or a power series, depending on the context) in $z$. Finally, recalling that $\det V = \wedge^r (V)$ if $r = \dim V$, let us consider
\begin{equation}
\label{eqn:sdet}
\sdet \ V = (-1)^r \det V
\end{equation}
for any vector space (or locally free sheaf, equivariant or not) $V$, and extend this notation additively to formal differences \eqref{eqn:k-th}
$$
\sdet(V-W) = \frac {\sdet \ V}{\sdet \ W}
$$ 

\medskip

\subsection{}
\label{sub:notation 2}

In the present paper, all our schemes will be quasiprojective over $\BC$. Given such a scheme $X$ and a two-step complex of locally free sheaves on $X$
$$
\CU = \Big[ \CW \longrightarrow \CV\Big]
$$
(with $\CV$ in homological degree 0) we define the projectivization
\begin{equation}
\label{eqn:diagram projectivization}
\xymatrix{
\BP_X(\CU) \ar[rd]_{\pi} \ar@{^{(}->}[r]
&\BP_X(\CV) \ar[d]^{\rho} \\
&X}
\end{equation}
as the subscheme of $\BP_X(\CV) = \text{Proj}_X(\text{Sym}_X(\CV))$ cut out by the composition
$$
\rho^*(\CW) \longrightarrow \rho^*(\CV) \stackrel{\text{taut}}\twoheadrightarrow \CO(1)
$$
where $\CO(1)$ denotes the tautological line bundle on $\BP_X(\CV)$. If we let $\pi$ denote the composed map in \eqref{eqn:diagram projectivization}, then we have the following formula for the push-forwards of the (algebraic $K$-theory classes of the) powers of the tautological line bundle
\begin{equation}
\label{eqn:push power series}
\pi_* \left[ \delta \left( \frac {\CO(1)}z \right) \right] = \wedge^\bullet\left(-\frac {\CU}z \right) \Big|^z_{\infty - 0}
\end{equation}
where $\delta(z) = \sum_{k \in \BZ} z^k$ is the formal delta series, and the right-hand side of \eqref{eqn:push power series} denotes the formal series obtained as the difference of the power series expansions (firstly as $z \rightarrow \infty$ and secondly as $z \rightarrow 0$) of the rational function 
$$
\wedge^\bullet\left( - \frac {\CU}z \right) = \frac {\displaystyle\wedge^\bullet\left(\frac {\CW}z\right)}{\displaystyle\wedge^\bullet\left(\frac {\CV}z\right)}
$$
Formula \eqref{eqn:push power series} is well-known, but we refer to \cite[Proposition 5.19]{Surf} for an overview.

\medskip

\subsection{}
\label{sub:notation 3}

In the present paper, we let $K(X)$ denote the 0-th algebraic $K$-theory group of a smooth algebraic variety $X$. Given algebraic varieties $X$ and $Y$, a class
$$
\Gamma \in K(X \times Y)
$$
with proper support over $X$ is called a correspondence. To such a correspondence, one associates the operator
$$
\Phi_\Gamma : K(Y) \rightarrow K(X), \qquad \alpha \mapsto \text{proj}_{X*} \Big(\Gamma \cdot \text{proj}_Y^* (\alpha) \Big)
$$
where $\text{proj}_X$, $\text{proj}_Y$ denote the standard projections from $X \times Y$ onto the factors. Given correspondences $\Gamma \in K(X \times Y)$ and $\Gamma' \in K(Y \times Z)$, we have
$$
\Phi_{\Gamma \circ \Gamma'} = \Phi_{\Gamma} \circ \Phi_{\Gamma'}
$$
where the composition of correspondences is the operation defined by
\begin{equation}
\label{eqn:composition of correspondences}
\Gamma \circ \Gamma' = \text{proj}_{X \times Z *} \Big( \text{proj}_{X \times Y}^*(\Gamma) \cdot \text{proj}_{Y \times Z}^*(\Gamma') \Big)
\end{equation}
The assumption that the support of $\Gamma$ is proper over $X$ ensures that the push-forward map above is well-defined. 

\medskip

\noindent Let us now assume that $X$ is a smooth algebraic variety acted on by a torus $T$. We may consider the equivariant algebraic $K$-theory group $K_T(X)$, and all the constructions above apply. In particular, the equivariant localization theorem states that elements of $K_T(X)$ are completely determined (up to torsion) by their restriction to the connected components $F \subset X^T$. In particular, the composition rule \eqref{eqn:composition of correspondences} states that
\begin{equation}
\label{eqn:composition of equivariant correspondences}
\Gamma \circ \Gamma' \Big|_{F \times F'} = \sum_{G} \frac {\Gamma|_{F \times G} \cdot \Gamma' |_{G \times F'}}{\wedge^\bullet(\Tan \ Y|_{G})^\vee}
\end{equation}
for all connected components $F$, $F'$ of the fixed point sets of $X^T$, $Z^T$, where the sum in the right-hand side runs over all the connected components $G$ of $Y^T$.

\medskip

\noindent We may drop the properness assumption on the support of correspondences if we only work with $T$-varieties $X$ with proper fixed point set $X^T$ (this will be the case throughout the present paper). In this case, one may regard $K_T(X)$ as a module over $K_T(\text{point}) = \text{Rep}_T$, and consider the localized equivariant $K$-theory group
$$
K_T(X)_{\text{loc}} = K_T(X) \bigotimes_{\text{Rep}_T} \text{Frac} \left(\text{Rep}_T\right)
$$
The right-hand side of \eqref{eqn:composition of equivariant correspondences} is well-defined in $K_T(F \times F')_{\text{loc}}$ for any $\Gamma$ and $\Gamma'$, even without the assumption that the correspondences have proper support. 

\medskip

\subsection{}
\label{sub:notation 4}

For a proper smooth algebraic variety $X$, its Euler characteristic pairing
$$
K(X) \otimes K(X) \xrightarrow{(\cdot,\cdot)} \BZ
$$
is defined by
$$
( \alpha ,  \beta ) = \int_X \alpha \cdot \beta
$$
(above and throughout the present paper, $\int_X$ denotes push-forward from $X$ to a point, or equivalently, the operation of taking Euler characteristic). Any correspondence $\Gamma \in K(X \times Y)$ for proper $X$ and $Y$ yields adjoint operators with respect to the aforementioned pairing. In more details, consider the operators
$$
\Phi_{\Gamma} : K(Y) \leftrightharpoons K(X) : \Phi_{\Gamma^T}
$$
where $\Gamma^T$, called the transpose correspondence, is the element of $K(Y \times X)$ obtained from $\Gamma$ by permuting the factors $X$ and $Y$. Then we have
\begin{equation}
	\label{eqn:adjunction}
\Big(\alpha, \Phi_{\Gamma}(\beta) \Big) = \Big( \Phi_{\Gamma^T}(\alpha), \beta \Big)
\end{equation}
for all $\alpha \in K(X)$ and $\beta \in K(Y)$. 

\medskip

\noindent The constructions and statements above pertain equally well to the situation of a smooth algebraic $T$-variety $X$ with proper fixed point set $X^T$, as long as we consider the localized version of the Euler characteristic pairing
$$
K_T(X)_{\loc} \otimes K_T(X)_{\loc} \xrightarrow{(\cdot, \cdot)} \text{Frac}\left(\text{Rep}_T\right)
$$
This is because expressions such as $([\CO_X], [\CO_X])$ are not well-defined elements of $\text{Rep}_T$ if $X$ is not proper, but they are well-defined elements of $\text{Frac}(\text{Rep}_T)$ due to the equivariant localization formula. Similarly, the adjunction \eqref{eqn:adjunction} only requires $X^T$ and $Y^T$ to be proper, but the pairing $(\cdot,\cdot)$ must be valued in $\text{Frac}(\text{Rep}_T)$.

\bigskip

\section{The $K$-theoretic Hall algebra and Nakajima quiver varieties} 
\label{sec:k-ha}

\medskip

\subsection{} 

Let us fix a quiver $Q$, i.e. an oriented graph with vertex set $I$ and edge set $E$; loops and multiple edges are allowed. We will consider the monoid $\nn$, where throughout the present paper, the set $\BN$ will be thought to contain 0. Particularly important elements of $\nn$ are the zero vector $\b0 = (0,\dots,0)$, as well as 
$$
\bs^i = (0,\dots,0,1,0,\dots,0) \in \nn
$$
where the unique 1 lies on the $i$-th spot. We will define the dot product
\begin{equation}
\label{eqn:dot product}
\zz \otimes \zz \xrightarrow{\cdot} \BZ, \qquad \bw \cdot \bv = \sum_{i\in I} w_iv_i
\end{equation}
the bilinear form
\begin{equation}
\label{eqn:bilinear form}
\zz \otimes \zz \xrightarrow{\langle \cdot,\cdot\rangle } \BZ, \qquad \langle \bv,\bw \rangle = \bv \cdot \bw - \sum_{\oij \in E} v_iw_j
\end{equation}
and the symmetric form
\begin{equation}
\label{eqn:symmetric form}
\zz \otimes \zz \xrightarrow{(\cdot,\cdot)} \BZ, \qquad (\bv,\bw) = \langle \bv,\bw\rangle+\langle \bw,\bv\rangle
\end{equation}
for any $\bv = (v_i)_{i \in I}$ and $\bw = (w_i)_{i \in I} \in \zz$. Given $\bv = (v_i)_{i \in I} \in \nn$, a \textbf{quiver representation} of dimension vector $\bv$ is a collection of vector spaces
\begin{equation}
\label{eqn:collection}
V_\bullet = (V_i)_{i \in I}
\end{equation}
with $\dim V_i = v_i$ for all $i \in I$, endowed with a choice of linear maps 
$$
\left( X_e : V_i \rightarrow V_j \right)_{e = \oij \in E} 
$$
Multiple edges from $i$ to $j$ correspond to multiple linear maps from $V_i$ to $V_j$ as part of the datum. A \textbf{double quiver representation} is a collection \eqref{eqn:collection} together with a choice of linear maps going in both directions, namely
\begin{equation}
\label{eqn:doubled}
\left(  V_i \overset{X_e}{\underset{Y_e}{\rightleftharpoons}} V_j  \right)_{e = \oij \in E} 
\end{equation}
A \textbf{framed double quiver representation} is a collection of \eqref{eqn:collection} and \eqref{eqn:doubled} as above, together with a choice of linear maps
\begin{equation}
\label{eqn:framed}
\left(W_i \overset{A_i}{\underset{B_i}{\rightleftharpoons}} V_i\right)_{i \in I}
\end{equation}
for certain fixed vector spaces $W_i$ of dimension $w_i$, $\forall i \in I$. The tuple $\bw = (w_i)_{i \in I}$ is called the framing vector of the framed double quiver representation in question.

\medskip

\subsection{}
\label{sub:nakajima}

The main geometric objects studied in the present paper are the \textbf{Nakajima quiver varieties} (\cite{Nak 0}) associated to any $\bv = (v_i)_{i \in I}, \bw = (w_i)_{i \in I} \in \nn$. To define these, consider the affine space of framed double quiver representations
$$
M_{\bv,\bw} = \bigoplus_{e = \oij \in E} \Big[ \Hom(V_i, V_j) \oplus  \Hom(V_j, V_i) \Big] \bigoplus_{i \in I} \Big[ \Hom(W_i, V_i) \oplus \Hom(V_i, W_i) \Big]
$$
where $V_i$ (respectively $W_i$) are vector spaces of dimension $v_i$ (respectively $w_i$) for all $i \in I$. Points of the affine space above will be denoted by quadruples
\begin{equation}
\label{eqn:quadruples}
(X_e,Y_e,A_i,B_i)_{e \in E, i\in I}
\end{equation}
where $X_e, Y_e, A_i, B_i$ denote homomorphisms in the four types of Hom spaces that enter the definition of $M_{\bv,\bw}$. Consider the action of
\begin{equation}
\label{eqn:gauge group}
G_{\bv} = \prod_{i \in I} GL(V_i)
\end{equation}
on $M_{\bv,\bw}$ given by conjugation of $X_e, Y_e$, left multiplication of $A_i$ and right multiplication of $B_i$. It is easy to see that $G_{\bv}$ acts freely on the open locus of \textbf{stable} \footnote{In general, stability conditions for Nakajima quiver varieties are indexed by $\bth \in \mathbb{R}^I$; the one studied herein corresponds to $\bth \in \BR_+^I$.} representations
\begin{equation}
\label{eqn:stable points}
M_{\bv,\bw}^s \subset M_{\bv,\bw}
\end{equation}
i.e. quadruples \eqref{eqn:quadruples} such that there does not exist a collection of subspaces $\{V_i' \subseteq V_i\}_{i \in I}$ (other than $V_i' = V_i$ for all $i \in I$) which is preserved by the maps $X_e$ and $Y_e$, and contains $\text{Im }A_i$ for all $i \in I$. Let us consider the quadratic moment map
\begin{equation}
\label{eqn:moment}
M_{\bv,\bw} \xrightarrow{\mu_{\bv,\bw}} \text{Lie } G_{\bv} = \bigoplus_{i \in I} \text{Hom}(V_i,V_i)
\end{equation}
given by
$$
\mu_{\bv,\bw} ( (X_e,Y_e,A_i,B_i)_{e\in E, i\in I}) = \sum_{e \in E} \Big(X_e Y_e - Y_e X_e \Big) + \sum_{i \in I} A_iB_i
$$
If we write $\mu^{-1}_{\bv,\bw}(0)^s = \mu^{-1}_{\bv,\bw}(0) \cap M_{\bv,\bw}^s$, then there is a geometric quotient
\begin{equation}
\label{eqn:quiver var}
\CM_{\bv,\bw} = \mu^{-1}_{\bv,\bw}(0)^s \Big / G_{\bv}
\end{equation}
which is called the \textbf{Nakajima quiver variety} for the quiver $Q$, associated to the dimension vector $\bv$ and the framing vector $\bw$. \\

\subsection{}
\label{sub:torus}

There are two commuting algebraic group actions on Nakajima quiver varieties. The first of these is the torus action
\begin{equation}
\label{eqn:small torus}
T = \BC^* \times (\BC^*)^{E} \curvearrowright \CM_{\bv,\bw}
\end{equation}
given by
\begin{equation}
\label{eqn:small torus acts}
\left( \bar{q}, \bar{t}_e \right)_{e\in E} \cdot (X_e,Y_e,A_i,B_i)_{e \in E, i \in I} = \left(\frac {X_e}{\bar{t}_e}, \frac {\bar{t}_e Y_e}{\bar{q}}, A_i , \frac {B_i}{\bar{q}} \right)_{e\in E, i \in I}
\end{equation}
and the second of these is the group action
\begin{equation}
\label{eqn:big torus}
G_{\bw} = \prod_{i \in I} GL(W_i) \curvearrowright \CM_{\bv,\bw}
\end{equation}
given by
$$
\left(\bar{U}_{i} \right)_{i\in I} \cdot (X_e,Y_e,A_i,B_i)_{e \in E, i \in I} = \left(X_e, Y_e, A_i \bar{U}_{i}^{-1}, \bar{U}_{i} B_i \right)_{e\in E, i \in I}
$$
We will denote the elementary characters of the torus \eqref{eqn:small torus} by $\{q,t_e\}_{e \in E}$, and the elementary characters of (a maximal torus of) \eqref{eqn:big torus} by $\{u_{i1},\dots,u_{iw_i}\}_{i \in I}$. Thus, the equivariant algebraic $K$-theory group $K_{T \times G_{\bw}}(\CM_{\bv,\bw})$ is a module over the ring
\begin{equation}
\label{eqn:base ring}
K_{T \times G_{\bw}}(\text{point})  = \BZ[q^{\pm \frac 12}, t_e^{\pm 1}, u_{ia}^{\pm 1}]^{\text{sym}}_{e\in E,i\in I, 1 \leq a \leq w_i}
\end{equation}
(where ``sym" means symmetric in the parameters $u_{i1},\dots,u_{iw_i}$ for each $i \in I$ separately). We fix a square root $q^{\frac 12}$ in \eqref{eqn:base ring} to ensure that subsequent formulas will be well-defined; rigorously speaking, this can be achieved by replacing the torus \eqref{eqn:small torus} by a two-fold cover. We will localize our $K$-theory groups, i.e. consider
\begin{equation}
\label{eqn:k v w}
K(\bw)_{\bv} = K_{T \times G_{\bw}}(\CM_{\bv,\bw}) \bigotimes_{K_{T \times G_{\bw}}(\text{point})}  \BF_{\bw}
\end{equation}
where we set
\begin{equation}
\label{eqn:ground field}
\BF = \BQ(q^{\frac 12},t_e)_{e \in E}
\end{equation}
and $\BF_{\bw} = \BF[u_{i1}^{\pm 1},\dots, u_{iw_i}^{\pm 1}]_{i\in I}^{\sym}$. In other words, we only localize our $K$-theory groups with respect to the $T$ action, and not with respect to the $G_{\bw}$ action.

\medskip

\begin{remark}
\label{rem:smaller torus}

The contents of the present paper still hold if we replace $T$ of \eqref{eqn:small torus} by any subtorus $H \subset T$ which contains a one-parameter subgroup $\alpha$ on which the characters 
\begin{align*}
&\BC^* \stackrel{\alpha}{\hookrightarrow} H \subset T \xrightarrow{t_e} \BC^* \\
&\BC^* \stackrel{\alpha}{\hookrightarrow} H \subset T \xrightarrow{\frac q{t_e}} \BC^*
\end{align*}
have strictly positive weight for all $e \in E$ (in the present context, this condition is equivalent to Assumption \begin{otherlanguage*}{russian}Ъ\end{otherlanguage*} of \cite{Wheel}).
 
\end{remark}

\medskip

\subsection{}

We will use the notation $\CV_i$ for the rank $v_i$ \textbf{tautological vector bundle} on $\CM_{\bv,\bw}$, whose fiber over a framed double quiver representation is the vector space $V_i$. The defining maps of a quiver representation induce maps of vector bundles
$$
X_e : t_e \CV_i \rightarrow \CV_j, \qquad Y_e : \frac q{t_e} \CV_j \rightarrow  \CV_i $$
$$
A_i : W_i \rightarrow \CV_i, \qquad \ \ B_i : q \CV_i \rightarrow W_i
$$
for any $i,j \in I$ and any edge $e = \oij$ (above and henceforth, we abuse the notation $W_i$ to denote the equivariant vector bundle $\CO_{\CM_{\bv,\bw}} \otimes W_i$). The reason we twist the vector bundles by various characters of $T$ in the formulas above is to ensure that the resulting maps are $T$-equivariant. It is well-known (\cite{Nak 0}) that the tangent bundle of Nakajima quiver varieties has the following formula in equivariant algebraic $K$-theory
\begin{equation}
\label{eqn:tangent bundle}
[\Tan \ \CM_{\bv,\bw} ] = \sum_{e = \oij} \left(\frac { \CV_j}{t_e\CV_i} + \frac {t_e\CV_i}{q\CV_j}\right) - \sum_{i \in I} \left(1+\frac 1q\right)\frac {\CV_i}{\CV_i} + \sum_{i \in I} \left( \frac {\CV_i}{W_i} + \frac {W_i}{q\CV_i} \right)
\end{equation}
(see Subsection \ref{sub:notation 1} for the notation of ``ratios" of locally free sheaves). Coupled with the well-known smoothness of $\CM_{\bv,\bw}$ (\cite{Nak 0}), formula \eqref{eqn:tangent bundle} implies that
\begin{equation}
\label{eqn:dim nakajima}
\dim \CM_{\bv,\bw} = 2\bw \cdot \bv - (\bv, \bv)
\end{equation} 
where the dot product and the symmetric form were defined in \eqref{eqn:dot product} and \eqref{eqn:symmetric form}. 

\medskip

\subsection{}
\label{sub:k-ha}

When $\bw = \b0$, the construction of Subsection \ref{sub:nakajima} yields the empty set $\forall \bv \neq \b0$, because of the stability  \eqref{eqn:stable points}. However, we may define instead the stack quotient
$$
\fY_{\bv} = \mu_{\bv,\b0}^{-1}(0) \Big/ GL_{\bv}
$$
whose equivariant $K$-theory groups are $K_T(\fY_{\bv}) = K_{T \times GL_{\bv}}(\mu_{\bv,\b0}^{-1}(0))$. Closed points of $\fY_{\bv}$ are double quiver representations \eqref{eqn:doubled} satisfying the equation 
$$
\sum_{e \in E} \Big(X_e Y_e - Y_e X_e \Big) = 0
$$
We will be interested in the localized equivariant $K$-theory groups
\begin{equation}
\label{eqn:k groups}
\CA_{\bv} = K_T(\fY_{\bv}) \bigotimes_{K_T(\text{point})} \BF
\end{equation}
where we recall that $\BF = \text{Frac}(K_T(\text{point}))$. As shown in \cite{SV Hilb, YZ}, the direct sum
\begin{equation}
\label{eqn:k-ha}
\CA^+ = \bigoplus_{\bv \in \nn} \CA_{\bv}
\end{equation}
can be endowed with the so-called \textbf{preprojective $K$-theoretic Hall algebra} (K-HA) structure. To define the product, one first considers the stack of extensions
\begin{equation}
\label{eqn:ext diagram 1}
\xymatrix{& \fZ_{\bv',\bv''} \ar[ld]_{p'} \ar[d]_{p} \ar[rd]^{p''} & \\
\fY_{\bv'} & \fY_{\bv'+\bv''} & \fY_{\bv''}}
\end{equation}
where if the three stacks on the bottom row parameterize double quiver representations $V'_\bullet$, $V_\bullet$, $V''_\bullet$ of dimensions $\bv'$, $\bv'+\bv''$, $\bv''$ respectively, then the stack on the top row parameterizes short exact sequences of double quiver representations
\begin{equation}
\label{eqn:ses}
0 \longrightarrow V'_\bullet \longrightarrow V_\bullet \longrightarrow V''_\bullet \longrightarrow 0
\end{equation}
(the maps $p$, $p'$, $p''$ in \eqref{eqn:ext diagram 1} record $V_\bullet$, $V'_\bullet$, $V''_\bullet$, respectively). Secondly, one considers the tautological vector bundles $\CV_i'$ on $\fY_{\bv'}$ and $\CV''_i$ on $\fY_{\bv''}$ which parameterize the vector spaces which appear in \eqref{eqn:ses}. Then following \loccitt, one defines the multiplication in \eqref{eqn:k-ha} as the operation
\begin{equation}
\label{eqn:multiplication k-ha}
\CA_{\bv'} \otimes \CA_{\bv''} \xrightarrow{*} \CA_{\bv'+\bv''}
\end{equation}
$$
\alpha \otimes \beta \mapsto p_{*} \left(\sdet \left[ \sum_{i \in I} \frac {\CV_i'}{q^{\frac 12}\CV_i''} - \sum_{e = \oij \in E} \frac {t_e \CV_i'}{q^{\frac 12}\CV_j''} \right] \cdot (p' \times p'')^!(\alpha \boxtimes \beta) \right)
$$
(the precise refined pull-back $(p' \times p'')^!$ required above is explained in detail in \cite{YZ}, and it is the most technically involved part of the construction). The expression $\sdet[\dots]$ in \eqref{eqn:multiplication k-ha} is defined according to \eqref{eqn:sdet}, and it is chosen to ensure compatibility with formulas \eqref{eqn:action plus} and \eqref{eqn:action minus} below, as well as with the contents of Section \ref{sec:stable}. We remark that none of the aforementioned results require localization (i.e. tensoring with $\BF$ in \eqref{eqn:k groups}). However, the following crucially requires localization. 

\medskip

\begin{theorem}
\label{thm:generate}

(\cite[Theorem 1.2 and Corollary 2.16]{Wheel}) As a $\BF$-algebra, $\CA^+$ is generated by its direct summands $\CA_{\bs^i}$, as $i$ runs over $I$. The same holds if the torus $T$ in the definition of $\CA^+$ is replaced by any subtorus $H$ as in Remark \ref{rem:smaller torus}.

\end{theorem}

\medskip

\subsection{} 
\label{sub:algebra}

As $\fY_{\b0} = \text{point}$, we have $\CA_{\b0} = \BF$ sitting inside $\CA^+$ as its ground field. Similarly, it is easy to see that for all $i \in I$, we have
\begin{equation}
\label{eqn:simple stack}
\fY_{\bs^i} = \BA^{2g_i} \Big / \BC^*
\end{equation}
where $g_i$ denotes the number of loops at the vertex $i$, and the $\BC^*$ action is trivial. If we identify the tautological line bundle on $\text{point}/\BC^*$ with $z_i$, we obtain an identification
$$
\CA_{\bs^i} = \BF[z_i^{\pm 1}]
$$
In this language, Theorem \ref{thm:generate} can be restated as the fact that $\CA^+$ is generated by
\begin{equation}
\label{eqn:generators e}
e_{i,d} = z_i^d \cdot [\CO_{\text{origin in }\BA^{2g_i}}] \in \CA_{\bs^i} \subset \CA^+
\end{equation}
for any $i \in I$ and $d \in \BZ$. We will also encounter the opposite algebra $\CA^- = \CA^{+,\op}$, and consider its direct sum decomposition analogous to \eqref{eqn:k-ha}
\begin{equation}
\label{eqn:k-ha opposite}
\CA^- = \bigoplus_{\bv \in \nn} \CA_{-\bv}
\end{equation}
Although $\CA_{\bv} = \CA_{-\bv}$ as $\BF$-vector spaces, we will denote them differently to emphasize the fact that they lie in different algebras $\CA^\pm$, and denote the generator $z_i^d\cdot [\CO_{\text{origin}}]$ of $\CA^-$ by $f_{i,d}$. It will sometimes be convenient to place the aforementioned generators into formal series
\begin{equation}
\label{eqn:formal series}
e_i(z) = \sum_{d \in \BZ} \frac {e_{i,d}}{z^d} \qquad \text{and} \qquad f_i(z) = \sum_{d \in \BZ} \frac {f_{i,d}}{z^d}
\end{equation}
Finally, we consider the polynomial ring
\begin{equation}
\label{eqn:polynomial ring}
\CA^0 = \BF \left[a_{i,\pm d}, b_{i,\pm d}, q^{\pm \frac {v_i}2}, q^{\pm \frac {w_i}2} \right]_{i \in I, d \geq 1}
\end{equation}
from which we construct the formal power series
\begin{equation}
\label{eqn:formal power series}
h_i^\pm(z) = h_{i,0}^{\pm 1} + \sum_{d = 1}^{\infty} \frac {h_{i,\pm d}}{z^{\pm d}} := (q^{\pm \frac 12})^{w_i - 2v_i + \sum_{e= \oij} v_j + \sum_{e= \oji} v_j}
\end{equation}
$$
\exp \left(\sum_{d=1}^{\infty} \frac {b_{i,\pm d} - a_{i,\pm d}(1 + q^{\pm d}) + \sum_{e = \oij} a_{j,\pm d} q^{\pm d} t_e^{\mp d} + \sum_{e = \oji} a_{j,\pm d}t_e^{\pm d}}{d z^{\pm d}} \right)
$$
for every $i \in I$.

\medskip

\begin{definition}
\label{def:kha}

The \textbf{double $K$-theoretic Hall algebra} is the $\BF$-algebra
\begin{equation}
\label{eqn:double shuffle}
\CA = \CA^+ \otimes \CA^0 \otimes \CA^{-}
\end{equation}
with the multiplication between the factors governed by the relations
\begin{equation}
\label{eqn:rel double 0}
q^{\pm \frac {w_i}2} \text{ and } b_{i,\pm d} \text{ are central}
\end{equation}
\begin{equation}
\label{eqn:rel double 1}
e_i(z)q^{\pm \frac {v_j}2} = q^{\pm \frac {v_j}2} e_i(z) \cdot q^{\mp \frac {\delta_{ij}}2} 
\end{equation}
\begin{equation}
\label{eqn:rel double 2}
f_i(z)q^{\pm \frac {v_j}2} = q^{\pm \frac {v_j}2} f_i(z) \cdot q^{\pm \frac {\delta_{ij}}2}
\end{equation}
\begin{equation}
\label{eqn:rel double 3}
[e_i(z), a_{j,d}] = e_i(z) \cdot \delta_{ij} z^{d}(q^{- d}-1) 
\end{equation}
\begin{equation}
\label{eqn:rel double 4}
[f_i(z), a_{j,d}] = f_i(z) \cdot \delta_{ij} z^{d}(1 - q^{- d}) 
\end{equation}
and
\begin{equation}
\label{eqn:rel double 5}
[e_{i,d}, f_{j,k}] =  \delta_{ij} \cdot \gamma_i \begin{cases} - h_{i,d+k} &\text{if } d+k > 0 \\ h_{i,0}^{-1} - h_{i,0} &\text{if } d+k = 0 \\ h_{i,d+k} &\text{if } d+k<0 \end{cases}
\end{equation}
for all $i,j \in I$ and $d,k \in \BZ$, where
\begin{equation}
\label{eqn:constant}
\gamma_i = \frac {\prod_{e = \oii} \left[ \left(q^{\frac 12}-t_eq^{-\frac 12}\right)\Big(1-t_e\Big) \right]}{q^{\frac 12} - q^{-\frac 12}}
\end{equation}
and $h_{i,\pm d}$ denote the coefficients of the power series \eqref{eqn:formal power series}.

\end{definition}

\medskip

\begin{proposition}
\label{prop:generate}

The algebra $\CA$ is generated by 
\begin{equation}
\label{eqn:generators}
\Big\{ e_{i,0}, f_{i,0}, q^{\pm \frac {v_i}2}, a_{i,\pm 1} \Big\}_{i \in I}
\end{equation}
together with the central elements $\left\{q^{\pm \frac {w_i}2}, b_{i,\pm d}\right\}_{i\in I, d\geq 1}$.

\end{proposition}

\medskip

\begin{proof} Let $\CB \subseteq \CA$ denote the subalgebra generated by \eqref{eqn:generators} and the central elements. Using relations \eqref{eqn:rel double 3}--\eqref{eqn:rel double 4}, one can successively show that $e_{i,d}, f_{i,d} \in \CB$ for all $i \in I$, $d \in \BZ$. Theorem \ref{thm:generate} therefore implies that $\CA^\pm \subset \CB$. As for the tensor factor $\CA^0$, it follows from \eqref{eqn:rel double 5} that $h_{i,\pm d} \in \CB$ for all $i\in I$ and $d \geq 1$. Since $q^{\pm \frac {v_i}2}$ and the central elements are already in $\CB$, then it follows from \eqref{eqn:formal power series} that 
$$
a_{i,\pm d}(1 + q^{\pm d}) - \sum_{e = \oij} a_{j,\pm d}q^{\pm d}  t_e^{\mp d} - \sum_{e = \oji} a_{j,\pm d} t_e^{\pm d} \in \CB
$$
for all $i \in I$ and $d \geq 1$. We claim that the $|I| \times |I|$ square matrix $C$ defined by
\begin{equation}
\label{eqn:square matrix}
C_{ij} = (1 + q^{\pm d}) \delta_{ij} - \sum_{e = \oij} q^{\pm d} t_e^{\mp d} - \sum_{e = \oji}t_e^{\pm d}
\end{equation}
is invertible, which would prove that $a_{i,\pm d} \in \CB$ for all $i \in I$, $d \geq 1$, and thus $\CB = \CA$. The aforementioned invertibility holds because there is a single summand in the Leibniz formula for $\det C$ which produces the monomial $q^{\pm d|I|}$, namely the one obtained by multiplying out the diagonal elements. All the other summands in the Leibniz formula produce powers of $q$ which are smaller in absolute value (this argument also applies in the setting of Remark \ref{rem:smaller torus}), hence $\det C$ cannot be 0. \end{proof}

\medskip

\subsection{} 
\label{sub:action}

Fix a framing vector $\bw \in \nn$. By analogy with \eqref{eqn:ext diagram 1}, one defines the stack
\begin{equation}
\label{eqn:ext diagram 2}
\xymatrix{& \CN_{\bv,\bv+\bn,\bw} \ar[ld]_{p} \ar[d]_{\pi_+} \ar[rd]^{\pi_-} & \\
\fY_{\bn} & \CM_{\bv+\bn,\bw} & \CM_{\bv,\bw}}
\end{equation}
parameterizing short exact sequences
\begin{equation}
\label{eqn:ses 2}
0 \longrightarrow K_\bullet \longrightarrow V^+_\bullet \longrightarrow V^-_\bullet \longrightarrow 0
\end{equation}
where $V^+_\bullet, V^-_\bullet$ are stable framed double quiver representations with dimension vectors $\bv+\bn, \bv$ (respectively) and framing vector $\bw$, while $K_\bullet$ is a unframed double quiver representation of dimension vector $\bn$; the maps $p$, $\pi_+$, $\pi_-$ in \eqref{eqn:ext diagram 2} remember $K_\bullet$, $V^+_\bullet$, $V^-_\bullet$, respectively. The maps
\begin{equation}
\label{eqn:action plus}
\CA_{\bn} \otimes K(\bw)_{\bv} \rightarrow K(\bw)_{\bv+\bn}
\end{equation}
$$
\alpha \otimes \beta \mapsto \pi_{+*} \left(\sdet \left[ \sum_{e = \oji \in E} \frac {t_e\CV_j^+}{q^{\frac 12}\CK_i} - \sum_{i\in I} \frac {\CV_i^+}{q^{\frac 12} \CK_i} + \sum_{i\in I}\frac {W_i}{q^{\frac 12}\CK_i} \right] \cdot (p \times \pi_-)^!(\alpha \boxtimes \beta) \right) 
$$
\begin{equation}
\label{eqn:action minus}
\CA_{-\bn} \otimes K(\bw)_{\bv+\bn} \rightarrow K(\bw)_{\bv} 
\end{equation}
$$
\alpha \otimes \beta \mapsto \pi_{-*} \left(\sdet \left[ \sum_{e = \oij \in E} \frac {q^{\frac 12}\CV^-_j}{t_e \CK_i} - \sum_{i\in I} \frac {q^{\frac 12}\CV_i^-}{\CK_i} \right] \cdot (p \times \pi_+)^!(\alpha \boxtimes \beta) \right)
$$
(above, $\CK_{i},\CV_i^+, \CV_i^-$ denote the tautological vector bundles on $\fY_{\bn},\CM_{\bv+\bn,\bw},\CM_{\bv,\bw}$, respectively, compatible with the notation in \eqref{eqn:ses 2}) yield actions
\begin{equation}
\label{eqn:two actions}
\CA^\pm \curvearrowright K(\bw) := \bigoplus_{\bv \in \nn} K(\bw)_{\bv}
\end{equation}
The interested reader may find details (particularly the subtle aspect of properly defining the refined pull-back maps $(p \times \pi_\pm)^!$) in \cite[Section 5]{YZ}. We note that our particular choices of line bundles $\sdet[\dots]$ in \eqref{eqn:action plus}--\eqref{eqn:action minus} are non-standard, but they do not affect the construction in any significant way, and are chosen so as to match the stable envelope formulas that will appear in Section \ref{sec:stable}. 

\medskip

\begin{remark}

The map $\pi_-$ in \eqref{eqn:ext diagram 2} is not proper in general. As a consequence, the map \eqref{eqn:action minus} is only defined in localized equivariant $K$-theory. For certain choices of $K$-theory class $\alpha$, the map $\pi_-$ will be proper over the support of $p^!(\alpha)$, which implies that \eqref{eqn:action minus} is actually defined in integral equivariant $K$-theory; this will be the case for the operators $f_{i,d}$ in the next Subsection. \end{remark}

\medskip

\subsection{}
\label{sub:action simple}

In the particular case $\bn = \bs^i$ for $i \in I$, the stack $\CN_{\bv,\bv+\bs^i,\bw}$ is actually a quasiprojective scheme. Moreover, it is smooth, as shown in \cite{Nak 0} (the proof therein applies to the case when the number $g_i$ of loops at the vertex $i \in I$ is 0, but the general case is analogous). The natural map 
\begin{equation}
\label{eqn:fiber}
\CN_{\bv,\bv+\bs^i,\bw} \xrightarrow{p} \fY_{\bs^i} = \BA^{2g_i} \Big / \BC^*
\end{equation}
gives rise to the tautological line bundle $\CL_i$ on $\CN_{\bv,\bv+\bs^i,\bw}$, by pulling back the tautological line bundle on $\text{point}/\BC^*$. In terms of closed points \eqref{eqn:ses 2}, the line bundle $\CL_i$ simply parameterizes the one-dimensional vector space $K_i$. We will actually consider the fiber of \eqref{eqn:fiber} over the origin of $\BA^{2g_i}$, i.e. work with
\begin{equation}
\label{eqn:punctual}
\bar{\CN}_{\bv,\bv+\bs^i,\bw} = \left\{ \left(0 \rightarrow \BC^{\delta_{\bullet i}} \rightarrow V_\bullet^+ \rightarrow V_\bullet^- \rightarrow 0, \text{ such that } X,Y\Big|_{\BC^{\delta_{\bullet i}}} = 0 \right) \right\}
\end{equation}
and consider the projections $\bar{\pi}_-, \bar{\pi}^+ : \bar{\CN}_{\bv,\bv+\bs^i,\bw} \rightarrow \CM_{\bv,\bw}, \CM_{\bv+\bs^i, \bw}$. Under the action \eqref{eqn:two actions}, the generators of $\CA^\pm$ act on $K(\bw)$ via Nakajima's correspondences
\begin{align}
&e_{i,d} = \bar{\pi}_{+*} \left( \CL_i^d \cdot \sdet \left[ \sum_{e = \oji \in E} \frac {t_e\CV_j^+}{q^{\frac 12}\CL_i} - \frac {\CV_i^+}{q^{\frac 12}\CL_i} + \frac {W_i}{q^{\frac 12}\CL_i} \right] \cdot \bar{\pi}_-^* \right) \label{eqn:action e} \\
&f_{i,d} = \bar{\pi}_{-*} \left( \CL_i^d \cdot \sdet \left[ \sum_{e = \oij \in E} \frac {q^{\frac 12}\CV^-_j}{t_e \CL_i} - \frac {q^{\frac 12} \CV_i^-}{\CL_i} \right] \cdot \bar{\pi}_+^* \right)  \label{eqn:action f}
\end{align}
for any $i \in I$ and $d \in \BZ$. 

\medskip

\subsection{}
\label{sub:power sum}

For any $d \in \BZ \backslash 0$, the $d$-th Adams operation of a $K$-theory class $X = \sum_i \pm x_i$ (where we formally assume the $x_i$'s are line bundles) is the $K$-theory class
$$
p_d(X) = \sum_i \pm x_i^d
$$
With this in mind, let us define the operators
\begin{equation}
\label{eqn:action a}
a_{i,d} : K(\bw)_{\bv} \rightarrow K(\bw)_{\bv}, \qquad a_{i,d} = \text{ multiplication by } p_d(\CV_i(1-q^{-1}))
\end{equation}
\begin{equation}
\label{eqn:action b}
\text{ } b_{i,d} : K(\bw)_{\bv} \rightarrow K(\bw)_{\bv}, \qquad \ b_{i,d} = \text{ multiplication by } p_d(W_i(1-q^{-1}))
\end{equation}
as well as
\begin{equation}
\label{eqn:action v and w}
q^{\frac {v_i}2}, q^{\frac {w_i}2} \text{ act on }K(\bw)_{\bv} \text{ by multiplication with the obvious power of }q^{\frac 12}
\end{equation}
As a consequence of the formulas above, the power series \eqref{eqn:formal power series} act on $K(\bw)_{\bv}$ as 
\begin{multline}
\label{eqn:multiplication}
h^\pm_i(z) (\alpha) = \\ = q^{\frac {\bw \cdot \bs^i - ( \bv, \bs^i )}2} \wedge^\bullet \left( \frac {q^{-1}-1}z \left[ W_i - ( 1+q) \CV_i  + \sum_{e = \oij}     \frac {q\CV_j}{t_e} + \sum_{e = \oji} t_e \CV_j \right] \right) \cdot \alpha
\end{multline}
expanded as a power series in $z^{\pm 1}$. Taking the leading order term of the expression above in $z^{- 1}$, we infer that
\begin{equation}
\label{eqn:multiplication leading}
h_{i, 0}(\alpha) = q^{\frac {\bw \cdot \bs^i - ( \bv, \bs^i )}2} \cdot \alpha
\end{equation}
The following theorem is a well-known generalization of \cite[Theorem 9.4.1]{Nak}.

\medskip

\begin{theorem}
\label{thm:action}

For any $\bw \in \nn$, there is an action 
$$
\CA \curvearrowright K(\bw)
$$
for which $e_{i,d}, f_{i,d}, a_{i, d}, b_{i,d}$ act via \eqref{eqn:action e}, \eqref{eqn:action f}, \eqref{eqn:action a}, \eqref{eqn:action b}, respectively.

\end{theorem}
 
\medskip 

\begin{proof} The fact that the operators $e_{i,d}$ (respectively the operators $f_{i,d}$) respect the relations in the algebra $\CA^+$ (respectively $\CA^-$) is built into the very definition of $K$-theoretic Hall algebras, see \cite{SV Hilb, YZ}. In notation closer to ours, the interested reader may find a proof of this fact in \cite[Theorem 4.16]{R-matrix}. Meanwhile, relations \eqref{eqn:rel double 0}--\eqref{eqn:rel double 4} are all straightforward exercises that we leave to the reader. For instance, \eqref{eqn:rel double 3} is an immediate consequence of the following relation
$$
p_d(\CV^-_j(1-q^{-1})) - p_d(\CV^+_j(1-q^{-1})) = \delta_{ij} \CL_i^d (q^{-d}-1) \in K_{T \times G_{\bw}}(\bar{\CN}_{\bv,\bv+\bs^i,\bw})
$$
where $\CL_i$ denotes the tautological line bundle on $\bar{\CN}_{\bv,\bv+\bs^i,\bw}$, while $\CV_j^-$ and $\CV_j^+$ denote the pull-backs of the $j$-th tautological vector bundle from $\CM_{\bv,\bw}$ and $\CM_{\bv+\bs^i,\bw}$, respectively. Thus, it remains to prove \eqref{eqn:rel double 5}, to which we now turn.

\medskip

\noindent To this end, let us consider the \textbf{universal} complex (\cite[Subsection 4.2]{Nak 0}) of vector bundles on $\CM_{\bv,\bw}$, with the middle term sitting in homological degree 0
\begin{equation}
\label{eqn:universal complex}
\CU_i = \left[q \CV_i \xrightarrow{(B_i,- X_e,Y_{e'})} W_i \oplus \bigoplus_{e = \oij} \frac q{t_e} \CV_j   \oplus \bigoplus_{e' = \oji}  t_{e'} \CV_j \xrightarrow{(A_i,Y_e ,X_{e'})} \CV_i  \right]
\end{equation}
The stability condition forces the right-most map in the complex above to be pointwise surjective, and thus $\CU_i$ is quasi-isomorphic to the two-step complex of locally free sheaves 
$$
\left[ q \CV_i \xrightarrow{(B_i,-X_e,Y_{e'})} \text{Ker }(A_i,Y_e,X_{e'})_{e = \oij ,e' = \oji} \right]
$$
The relevance of the universal complex to us is the fact (due to \cite{Nak} when $g_i = 0$, but the general case is analogous) that the maps
\begin{align*}
&\bar{\pi}_- : \bar{\CN}_{\bv,\bv+\bs^i,\bw} \longrightarrow \CM_{\bv,\bw}  \\
&\bar{\pi}_+ : \bar{\CN}_{\bv-\bs^i,\bv,\bw} \longrightarrow \CM_{\bv,\bw} 
\end{align*}
are the projectivizations of the complex \eqref{eqn:universal complex} and its shifted dual, respectively
\begin{align}
&\bar{\CN}_{\bv,\bv+\bs^i,\bw} = \BP_{\CM_{\bv,\bw}}(\CU_i)  \label{eqn:projectivization 1} \\
&\bar{\CN}_{\bv-\bs^i,\bv,\bw} = \BP_{\CM_{\bv,\bw}}(q\CU_i^\vee[1]) \label{eqn:projectivization 2}
\end{align}
(see \eqref{eqn:diagram projectivization} for the definition of the projectivizations above).  Moreover, according to the notation \eqref{eqn:notation exterior}, the exterior algebra of $\CU_i$ is given by
$$
\wedge^\bullet \left(\frac {\CU_i}z \right) = \frac {\wedge^\bullet\left(\frac {W_i}z\right) \prod_{e = \oij} \wedge^\bullet\left(\frac {q\CV_j}{t_e z}\right) \prod_{e' = \oji} \wedge^\bullet\left(\frac {t_{e'}\CV_j }{z}\right)}{\wedge^\bullet\left(\frac {\CV_i}z\right)\wedge^\bullet\left(\frac {q\CV_i}z\right)}
$$
(as a rational function in $z$ whose coefficients lie in $K(\bw)_{\bv}$).

\medskip

\noindent Recall that the correspondences $e_{i,d}$ and $f_{i,d}$ are supported on the locus $\bar{\CN}_{\bv,\bv+\bs^i,\bw}$ parameterizing pairs of framed double quiver representations 
$$
(V^+_\bullet, V^-_\bullet) \quad \text{such that} \quad V^+_\bullet \stackrel{i}{\twoheadrightarrow} V^-_\bullet
$$
where the notation in the right-hand side means that there exists a surjective map of framed double quiver representations $V^+_\bullet \twoheadrightarrow V^-_\bullet$ which is an isomorphism if $\bullet \neq i$ and has a one-dimensional kernel if $\bullet = i$; moreover, the $X$ and $Y$ maps are required to annihilate the one-dimensional kernel. As explained in Subsection \ref{sub:notation 3}, the composition
$$
e_{i,d} \circ f_{j,k} \qquad \Big( \text{respectively } f_{j,k} \circ e_{i,d} \Big)
$$
is given by a $K$-theory class supported on the locus $(V_\bullet^1, V_\bullet^2) \in  \CM_{\bv+\vs^i-\bs^j,\bw} \times \CM_{\bv,\bw}$ (for various $\bv \in \nn$) with the property that there exists $V'_{\bullet} \in \CM_{\bv-\bs^j, \bw}$ (respectively $V''_{\bullet} \in \CM_{\bv+\bs^i, \bw}$) such that
$$
V^1_\bullet \stackrel{i}{\twoheadrightarrow} V'_\bullet \stackrel{j}{\twoheadleftarrow} V^2_\bullet \qquad \Big( \text{respectively } V^1_\bullet \stackrel{j}{\twoheadleftarrow} V''_{\bullet} \stackrel{i}{\twoheadrightarrow} V^2_\bullet \Big)
$$
If $i \neq j$ or if $i = j$ and $V^1_\bullet \not \cong V^2_\bullet$, then the existence of $V'_\bullet$ as above implies the existence of a unique $V''_\bullet$ as above, and vice-versa (by taking the pull-back and push-out of a diagram of vector spaces, respectively). Therefore, the difference
$$
e_{i,d} \circ f_{j,k} - f_{j,k} \circ e_{i,d} \in K_{T \times G_{\bw}}(\CM_{\bv+\bs^i-\bs^j,\bw} \times \CM_{\bv, \bw})
$$
is 0 if $i \neq j$, and is supported on the diagonal $\Delta \hookrightarrow \CM_{\bv,\bw} \times \CM_{\bv,\bw}$ if $i = j$. Thus, we may focus only on the case $i = j$ for the remainder of the proof. The aforementioned diagonal support condition implies that
\begin{equation}
\label{eqn:commutator support}
[e_{i,d}, f_{i,k}] = \Delta_* (\beta^i_{d,k})
\end{equation}
for some $\beta^i_{d,k} \in K(\bw)_{\bv}$. In light of formula \eqref{eqn:multiplication} prescribing the action of $h_{i,d+k}$ on $K(\bw)_{\bv}$, the required formula \eqref{eqn:rel double 5} boils down to proving the identity
$$
\beta^i_{d,k} = \gamma_i  \left\{ q^{ \frac {\bw \cdot \bs^i - ( \bv, \bs^i )}2}  \wedge^\bullet \left( \frac {q^{-1}-1}z \left[ W_i - (1+q)\CV_i  + \sum_{e = \oij} \frac {q\CV_j}{t_e} + \sum_{e = \oji} t_e \CV_j \right] \right) \right\}_{z^{d+k}}
$$
where $\gamma_i$ is the scalar in \eqref{eqn:constant}, and $\{\dots\}_{z^{d+k}}$ refers to the coefficient of $z^{d+k}$ when the rational function in question is expanded near $z = \infty$ minus the coefficient of $z^{d+k}$ when the rational function in question is expanded near $z = 0$. 

\medskip 

\noindent One can compute $\beta^i_{d,k}$ by applying formula \eqref{eqn:commutator support} to the unit $K$-theory class, so the proof of Theorem \ref{thm:action} boils down to establishing the following identity
\begin{equation}
\label{eqn:remains to show}
e_{i,d} \left( f_{i,k} ( 1 ) \right) - f_{i,k} \left( e_{i,d} (1) \right) = \frac {\prod_{e = \oii} \left[ \left(q^{\frac 12} - t_e q^{-\frac 12} \right)\Big(1-t_e\Big) \right]}{q^{\frac 12} -q^{-\frac 12}}
\end{equation}
$$
\left\{ q^{ \frac {\bw \cdot \bs^i - ( \bv, \bs^i )}2}  \wedge^\bullet \left( \frac {q^{-1}-1}z \left[ W_i - (1+q)\CV_i + \sum_{e = \oij} \frac {q\CV_j}{t_e} + \sum_{e = \oji} t_e\CV_j \right] \right) \right\}_{z^{d+k}} 
$$
To this end, let us apply \eqref{eqn:push power series} to formulas \eqref{eqn:projectivization 1} and \eqref{eqn:projectivization 2}, and we obtain
\begin{align*}
&e_{i,d}(1) = \left\{ q^{\frac {\bw \cdot \bs^i - \langle \bv,\bs^i \rangle}2} \cdot \frac {\wedge^\bullet\left(\frac {W_i}{qz}\right) \prod_{e = \oij} \wedge^\bullet\left(\frac {t_ez}{\CV_j}\right) \prod_{e' = \oji} \wedge^\bullet\left(\frac {t_{e'}\CV_j }{qz}\right)}{\wedge^\bullet\left(\frac {\CV_i}{qz}\right)\wedge^\bullet\left(\frac {z}{\CV_i}\right)} \right\}_{z^d} \\
&f_{i,k}(1) = \left\{ q^{\frac {\langle \bs^i, \bv \rangle}2} \cdot \frac {\wedge^\bullet\left(\frac {\CV_i}y\right)\wedge^\bullet\left(\frac y{q\CV_i}\right)}{\wedge^\bullet\left(\frac {W_i}{y}\right) \prod_{e = \oij} \wedge^\bullet\left(\frac {t_e y}{q \CV_j}\right) \prod_{e' = \oji} \wedge^\bullet\left(\frac { t_{e'}\CV_j}{y}\right)} \right\}_{y^k}
\end{align*}
Applying $f_{i,k}$ and $e_{i,d}$ to the formulas above, respectively (see \cite[Theorem 4.16]{R-matrix} for a more general version of this computation) gives us
$$
e_{i,d} (f_{i,k}(1)) = \left\{ \left\{ q^{\frac {\bw\cdot \bs^i + \langle \bs^i, \bv - \bs^i \rangle - \langle \bv, \bs^i \rangle}2} \cdot  \frac {\prod_{e = \oii} \left[ \left(1 - \frac {t_e y}{q z} \right) \left(1- \frac {t_{e'}z}y \right) \right]}{\left(1- \frac zy \right) \left(1 - \frac {y}{qz} \right)} \right. \right.
$$
\begin{equation}
	\label{eqn:for 1}
\left. \left. \frac {\wedge^\bullet\left(\frac {W_i}{qz}\right)}{\wedge^\bullet\left(\frac {W_i}{y}\right)} \cdot \frac {\wedge^\bullet\left(\frac {\CV_i}y\right)}{\wedge^\bullet\left(\frac {\CV_i}{qz}\right)} \cdot \frac {\wedge^\bullet\left(\frac y{q\CV_i}\right)}{\wedge^\bullet\left(\frac {z}{\CV_i}\right)}  \prod_{e = \oij} \frac {\wedge^\bullet\left(\frac {t_ez}{\CV_j}\right)}{\wedge^\bullet\left(\frac {t_ey}{q\CV_j}\right)} \prod_{e' = \oji} \frac {\wedge^\bullet\left(\frac {t_{e'}\CV_j}{qz}\right)}{\wedge^\bullet\left(\frac {t_{e'}\CV_j}{y}\right)} \right\}_{y^k} \right\}_{z^d}
\end{equation}
and 
$$
f_{i,k}(e_{i,d}(1)) =  \left\{ \left\{ q^{\frac {\bw\cdot \bs^i + \langle \bs^i, \bv \rangle - \langle \bv+\bs^i, \bs^i \rangle}2} \cdot \frac {\prod_{e = \oii} \left[ \left(1 - \frac {t_ey}{qz} \right) \left(1- \frac {t_{e'}z}y \right) \right]}{\left(1- \frac zy \right) \left(1 - \frac {y}{qz} \right)}  \right. \right.
$$
\begin{equation}
	\label{eqn:for 2}
\left. \left. \frac {\wedge^\bullet\left(\frac {W_i}{qz}\right)}{\wedge^\bullet\left(\frac {W_i}{y}\right)} \cdot \frac {\wedge^\bullet\left(\frac {\CV_i}y\right)}{\wedge^\bullet\left(\frac {\CV_i}{qz}\right)} \cdot \frac {\wedge^\bullet\left(\frac y{q\CV_i}\right)}{\wedge^\bullet\left(\frac {z}{\CV_i}\right)}  \prod_{e = \oij} \frac {\wedge^\bullet\left(\frac {t_ez}{\CV_j}\right)}{\wedge^\bullet\left(\frac {t_ey}{q\CV_j}\right)} \prod_{e' = \oji} \frac {\wedge^\bullet\left(\frac {t_{e'}\CV_j }{qz}\right)}{\wedge^\bullet\left(\frac {t_{e'}\CV_j }{y}\right)}  \right\}_{z^d} \right\}_{y^k}
\end{equation}
We note that the right-hand sides of the expressions above are identical, up to changing the order of the expansions (first in $y$ and second in $z$ versus first in $z$ and second in $y$). In general, for any rational function $F(z,y)$, the residue theorem implies that
$$
\Big\{ \Big\{ F(z,w) \Big \}_{y^k} \Big \}_{z^d} - \Big\{ \Big\{ F(z,w) \Big \}_{z^d} \Big \}_{y^k} = \sum_{c \neq 0} \left\{ c^{-k} \underset{y=zc}{\text{Res}} \frac {F(z,y)}y \right \}_{z^{d+k}}
$$
The integrand in \eqref{eqn:for 1} and \eqref{eqn:for 2} has two poles involving $y$ and $z$, namely $y = z$ and $y = zq$. The latter pole does not produce any residue as all rational functions cancel out when we set $y = zq$ (and $\{\text{constant}\}_{z^{d+k}} = 0$), while the former pole has residue given by the right-hand side of \eqref{eqn:remains to show}. \end{proof}

\medskip

\subsection{} As a consequence of Theorem \ref{thm:action}, we have a $\BF$-algebra homomorphism
\begin{equation}
\label{eqn:action hom}
\CA \xrightarrow{\iota} \prod_{\bw \in \nn} \text{End}(K(\bw))
\end{equation}

\medskip

\begin{proposition}
\label{prop:inj}

The homomorphism \eqref{eqn:action hom} is injective. 

\end{proposition}

\medskip

\noindent We will prove Proposition \ref{prop:inj} in Section \ref{sec:appendix}, as it requires us to develop certain equivariant localization tools. For the time being, we will need the following result (which was already proved in cohomology by \cite[Propositions 5.2 and 5.19]{SV Gen}).

\medskip

\begin{proposition}
\label{prop:any coefficient}

Any vector and covector of $K(\bw)$ can be written as
$$
\BF_{\bw} = K(\bw)_{\b0} \xrightarrow{e} K(\bw)_{\bm} \qquad \text{and} \qquad K(\bw)_{\bn} \xrightarrow{f} K(\bw)_{\b0} = \BF_{\bw}
$$
respectively, for various $e \in \CA_{\bm} \otimes \CA^0$ and $f \in \CA^0 \otimes \CA_{-\bn}$.

\end{proposition}

\medskip

\begin{proof} Let us consider the $\bv = \b0$ case of the action map \eqref{eqn:action plus}
$$
\text{act} : \CA_{\bm} \otimes \underbrace{K(\bw)_{\b0}}_{\BF_{\bw}} \rightarrow K(\bw)_{\bm}
$$
which is given by the formula
\begin{equation}
\label{eqn:action map}
\text{act} = \pi_{*}\Big(\text{line bundle} \cdot p^* \Big)
\end{equation}
with $\pi$ and $p$ as in the following diagram
\begin{equation}
\label{eqn:pi and gamma}
\xymatrix{
\CN_{\b0,\bm,\bw} \ar[d]_p \ar[r]^{\pi} & \CM_{\bm,\bw} \\
\fY_{\bm} &}
\end{equation}
The map $\pi$ is simply the closed embedding of the locus $\CN_{\b0,\bm,\bw} = \{B_i = 0, \forall i \in I\}$, while the map $p$ forgets $(A_i)_{i \in I}$. More explicitly,
$$
\xymatrix{
\CN_{\b0,\bm,\bw} \ar[rd]_-p \ar@{^{(}->}[r]^-j & \text{Tot}_{\fY_{\bm}}\left( \bigoplus_{i \in I} \text{Hom}(W_i,\CV_i) \right) \ar[d]^{\tau} \\
& \fY_{\bm}}
$$
where $j$ is the open embedding of the locus of stable triples $(X_e,Y_e,A_i,B_i = 0)$. Since $\tau$ is an affine bundle and $j$ is an open embedding, we have that $\tau^*$ and $j^*$ are surjective maps, and therefore so is $p^*$ in \eqref{eqn:action map}. Thus, to prove the surjectivity of $\text{act}$ (which would imply the required statement that every vector in $K(\bw)_{\bm}$ can be obtained by the action of elements $e \in \CA_{\bm} \otimes \CA^0$ on $K(\bw)_{\b0} = \BF_{\bw}$), it is enough to prove that the direct image map $\pi_*$ is surjective. Since the equivariant localization theorem implies that the direct image map
$$
K_{T \times G_{\bw}}(\CM_{\bm,\bw}^{T})_{\text{loc}} \rightarrow K_{T \times G_{\bw}}(\CM_{\bm,\bw})_{\text{loc}}
$$
is an isomorphism, it suffices to show that the fixed point locus $\CM_{\bm,\bw}^{T}$ is contained in $\CN_{\b0,\bm,\bw}$, i.e. that any quadruple $(X_e, Y_e, A_i, B_i)$ which is fixed by the $T$ action has the property that $B_i = 0$, $\forall i \in I$. However, a quadruple is fixed precisely when there exists a decomposition of its underlying vector spaces
$$
V_i = \mathop{\bigoplus_{\text{characters}}}_{\chi : T \rightarrow \BC^*} V_{i,\chi}
$$
such that the maps $(X_e,Y_e,A_i,B_i)$ are only non-zero between the direct summands
$$
X_e : V_{i,\chi} \rightarrow V_{j,\frac {\chi}{t_e}}, \quad Y_e : V_{j,\chi} \rightarrow V_{i,\frac {\chi t_e}q}
$$
$$
A_i : W_i \rightarrow V_{i,1}, \qquad \quad \ B_i : V_{i,q} \rightarrow W_i
$$
However, the stability condition (i.e. the fact that $\{V_i\}_{i \in I}$ are generated by the maps $X_e, Y_e$ acting on $\{\text{Im }A_i\}_{i \in I}$) implies that $V_{i,\chi} \neq 0$ only if $\chi$ is a product of $t_e$ and $\frac q{t_e}$ raised to non-positive powers. Since the character $q$ cannot be written as such a product due to the genericity of $t_e$ (in fact, this holds even if the $t_e$ are not generic, as long as the assumption in Remark \ref{rem:smaller torus} holds), we have $V_{i,q} = 0$ and thus $B_i = 0$ for all $i\in I$.

\medskip

\noindent Having proved the required statement of Proposition \ref{prop:any coefficient} for vectors, let us turn to the analogous statement for covectors. Recall from \cite[Proposition 7.3.4]{Nak} that
\begin{equation}
\label{eqn:euler pairing nakajima}
K(\bw)_{\bv} \otimes K(\bw)_{\bv} \xrightarrow{(\cdot, \cdot)} \BF_{\bw}, \qquad (\alpha,\beta) = \int_{\CM_{\bv,\bw}} \alpha \cdot \beta
\end{equation}
is a perfect pairing, where $\int_X$ denotes push-forward of equivariant $K$-theory classes from $X$ to a point (since our $K$-theory groups are localized with respect to the torus $T$, the argument of \loccit carries through because $\CM_{\bv,\bw}^T$ is projective, see Proposition \ref{prop:projective}). In other words, for any covector
$$
\lambda : K(\bw)_{\bn} \rightarrow \BF_{\bw}
$$
there exists  $\beta \in K(\bw)_{\bn}$ such that
\begin{equation}
\label{eqn:integrals 1}
\lambda(\alpha) = \int_{\CM_{\bn,\bw}} \alpha \cdot \beta
\end{equation}
From the already proved statement for vectors, there exists $e \in K_{T \times G_{\bw}}(\fY_{\bn})_{\loc} = K_{T}(\fY_{\bn})_{\loc} \otimes \BF_{\bw}$ such that, in the notation of \eqref{eqn:pi and gamma}, we have
$$
\beta = \pi_*(p^*(e))
$$
The usual adjunction of push-forward and pull-back implies
\begin{equation}
\label{eqn:integrals 2}
\int_{\CM_{\bn,\bw}} \alpha \cdot \beta = \int_{\CN_{0,\bn,\bw}} \pi^*(\alpha) \cdot p^*(e) 
\end{equation}
However, if we let $f = e$ times an appropriately chosen line bundle, the right-hand side of \eqref{eqn:integrals 2} is none other than $f(\alpha)$, i.e. the operator $K(\bw)_{\bn} \rightarrow K(\bw)_{\b0} = \BF_{\bw}$ that is induced by the $K$-theory class $f \in K_{T \times G_{\bw}}(\fY_{\bn})_{\loc}$ using the action map \eqref{eqn:action minus}. Combining this with \eqref{eqn:integrals 1} yields the required statement for covectors. \end{proof}

\medskip

\subsection{} 
\label{sub:complex}

We conclude the present Section with an alternative point of view  on the diagram \eqref{eqn:ext diagram 2}, with the goal of answering the following question. While \eqref{eqn:projectivization 1} and \eqref{eqn:projectivization 2} describe the projection maps from $\bCN_{\bv,\bv+\bs^i,\bw}$ to $\CM_{\bv,\bw}$ and $\CM_{\bv+\bs^i,\bw}$ individually, Proposition \ref{prop:lci} below will describe the closed embedding
$$
\bCN_{\bv,\bv+\bs^i,\bw} \hookrightarrow \CM_{\bv,\bw} \times \CM_{\bv+\bs^i, \bw}
$$
Although we will deal with the slightly more general setup of two distinct framing vectors, the construction and proof below are an almost word-for-word adaptation of \cite[Section 5.1]{Nak}.

\medskip

\begin{proposition}
\label{prop:lci}

Consider dimension vectors $\bv,\widetilde{\bv},\bw,\widetilde{\bw} \in \nn$, together with a choice of linear maps $\pi_i : \tW_i \rightarrow W_i$ for all $i \in I$. Consider the following complex
\begin{equation}
\label{eqn:complex}
C_{-1} \xrightarrow{f} C_0 \xrightarrow{g} C_1
\end{equation}
on $\CM_{\bv,\bw} \times \CM_{\widetilde{\bv},\widetilde{\bw}}$, where
\begin{align*}
&C_{-1} = \bigoplus_{i \in I} \emph{Hom}(\tCV_i, \CV_i), \qquad C_1 = \bigoplus_{i \in I} \emph{Hom}(q \tCV_i, \CV_i) \\
&C_0 = \bigoplus_{\oij \in E} \left( \emph{Hom}(t_e \tCV_i, \CV_j) \oplus \emph{Hom}\left( \frac q{t_e} \tCV_j, \CV_i \right) \right) \oplus \bigoplus_{i \in I} \left( \emph{Hom}(\tW_i,\CV_i) \oplus \emph{Hom}(q\tCV_i,W_i) \right)
\end{align*}
and the maps $f,g$ are given by
\begin{align*}
&f \left( \{\psi_i\}_{i \in I} \right) = (X_{\oij}\psi_i - \psi_j \tX_{\oij}, Y_{\oij} \psi_j - \psi_i  \tY_{\oij},  \psi_i \tA_i, B_i \psi_i) \\
&g\left( (x_{\oij}, y_{\oij}, a_i, b_i )\right) = \left(\sum_{\oij \in E} ( X_{\oji} y_{\oji} + x_{\oji} \tY_{\oji} - Y_{\oij} x_{\oij} - y_{\oij} \tX_{\oij} ) + \sum_{i \in I} ( A_i  b_i - a_i \tB_i )\right)
\end{align*}
The map $f$ is fiberwise injective and $g$ is fiberwise surjective, so the middle cohomology of the complex \eqref{eqn:complex} is a vector bundle $\CE$ of rank
\begin{equation}
\label{eqn:rank of complex}
\emph{rank } \CE = \bw \cdot \widetilde{\bv} + \widetilde{\bw} \cdot \bv - (\bv, \widetilde{\bv})
\end{equation}
Moreover, in the notation of the middle vector space in \eqref{eqn:complex}, the assignment 
\begin{equation}
\label{eqn:section}
s \left( (X_{\oij}, Y_{\oij}, A_i, B_i), (\tX_{\oij}, \tY_{\oij}, \tA_i, \tB_i) \right) = (0,0,A_i \pi_i, \pi_i \tB_i)
\end{equation}
yields a section
$$
s \in \Gamma \left(\CM_{\bv,\bw} \times \CM_{\widetilde{\bv},\widetilde{\bw}}, \CE \right)
$$
whose zero locus is the set of pairs of framed double quiver representations which fit into a diagram
\begin{equation}
\label{eqn:locus}
\xymatrix{{\tV}_\bullet \ar@/^/[d] \ar[r]^{\psi_\bullet} & V_\bullet \ar@/^/[d] \\ {\tW}_\bullet \ar@/^/[u] \ar[r]^{\pi_\bullet} & W_\bullet \ar@/^/[u]}
\end{equation}
The stability condition implies that if horizontal arrows $\psi_\bullet$ as in \eqref{eqn:locus} exist, then they are unique. If the maps $\pi_\bullet$ are all surjective, then so are the maps $\psi_\bullet$ in \eqref{eqn:locus}. 

\end{proposition}

\medskip

\begin{proof} It is easy to check that $g \circ f = 0$, and thus \eqref{eqn:complex} is a complex, due to the equations
$$
\sum_{\oij \in E} (X_{\oji}Y_{\oji} - Y_{\oij} X_{\oij}) + \sum_{i\in I} A_iB_i = \sum_{\oij \in E} (\tX_{\oji}\tY_{\oji} - \tY_{\oij} \tX_{\oij}) + \sum_{i\in I} \tA_i \tB_i = 0
$$
Let us show that the map $f$ is fiberwise injective. If $f(\{\psi_i\}_{i\in I}) = 0$ for a certain collection of linear maps 
\begin{equation}
\label{eqn:the psis}
\{\psi_i : \tV_i \rightarrow V_i\}_{i \in I}
\end{equation}
then the collection of subspaces $\{\text{Ker }\psi_i \subseteq \tV_i\}_{i\in I}$ is preserved by the $X,Y$ maps and contains the image of the $A$ maps. By stability, this requires 
$$
\text{Ker }\psi_i = \tV_i \quad \Rightarrow \quad \psi_i = 0, \ \forall i \in I
$$
Similarly, the dual of the map $g$ (using the trace pairing $\text{Hom}(\tV,V)^\vee \cong \text{Hom}(V,\tV)$) is shown to be fiberwise injective by the analogous argument to the one previously presented for $f$, and thus $g$ is surjective. We conclude that the middle cohomology $\CE$ of the complex \eqref{eqn:complex} is a vector bundle, and it is clear that its rank is \eqref{eqn:rank of complex}. Finally, it is easy to see that $s$ takes values in the kernel of $g$ (and thus yields a section of $\CE$) and it vanishes precisely when there exist linear maps as in \eqref{eqn:the psis} such that
$$
X_{\oij}\psi_i = \psi_j \tX_{\oij}, \quad Y_{\oij} \psi_j = \psi_i  \tY_{\oij} 
$$
$$
A_i\pi_i = \psi_i \tA_i, \qquad \quad \ \pi_i \tB_i = B_i \psi_i
$$
for all $\oij \in E$ and all $i \in I$. This precisely yields the top arrows in the commutative diagram \eqref{eqn:locus}. Finally, the surjectivity of the maps $\psi_i$ follows from the stability of the representation $\{V_i\}_{i \in I}$, together with the surjectivity of the maps $\pi_i$.  \end{proof}

\medskip

\subsection{} 
\label{sub:enhanced complex}

When $\bv = \widetilde{\bv}$ and $\bw = \widetilde{\bw}$, Proposition \ref{prop:lci} realizes the diagonal of $\CM_{\bv,\bw}$ as the zero locus of a regular section $s$. However, for general $\bv,\widetilde{\bv},\bw,\widetilde{\bw}$, the section $s$ will not necessarily be regular, as the rank of $\CE$ will be larger than the codimension of the locus \eqref{eqn:locus}. The mildest version of this problem is when $\widetilde{\bv} = \bv+\bs^i$ for some $i \in I$ and $\widetilde{\bw} = \bw$, in which case the locus \eqref{eqn:locus} is precisely $\CN_{\bv, \bv+\bs^i, \bw}$. It was shown in \cite{Nak} that
$$
\dim \CN_{\bv, \bv+\bs^i, \bw} = \bw \cdot (2\bv + \bs^i) - (\bv+\bs^i, \bv) + 2g_i - 1
$$
where $g_i$ denotes the number of loops at the vertex $i$ (strictly speaking, the result of \loccit applies to the case $g_i=0$, but the general case is analogous) hence the codimension of $\CN_{\bv, \bv+\bs^i, \bw}$ in $\CM_{\bv,\bw} \times \CM_{\bv+\bs^i,\bw}$ is equal to
$$
\text{codim } \CN_{\bv, \bv+\bs^i, \bw} = \bw \cdot (2\bv + \bs^i) - (\bv+\bs^i, \bv) - 1
$$
because of the identity $(\bs^i, \bs^i) = 2 - 2g_i$. However, the vector bundle $\CE$ of Proposition \ref{prop:lci} has rank one more than the codimension above, and so we will need to ``amend" it. The solution found in \cite{Nak}, which applies whenever we have $W_\bullet = \tW_\bullet$ in the setting of Proposition \ref{prop:lci}, is to consider the linear map
$$
C_0 \xrightarrow{h} q^{-1} \BC, \qquad (x_{e}, y_{e}, a_j, b_j )_{e \in E, j \in I} \mapsto \text{Tr }\sum_{j \in I} \left(a_j B_j - \tA_j b_j \right)
$$
Because $h \circ f = 0$, the map $h$ induces a map of vector bundles $h':  \CE \rightarrow q^{-1}\CO$. The image of the section $s$ lands in $\CE' = \text{Ker }h'$, because
\begin{equation}
\label{eqn:trace is zero}
\text{Tr}\sum_{j \in I} A_jB_j  = \text{Tr} \sum_{e \in E} \left(Y_{e}X_{e} - X_{e}Y_{e} \right) = 0
\end{equation}
and analogously for $\tA_j\tB_j$ instead of $A_jB_j$. Thus, to conclude the fact that $\CN_{\bv, \bv+\bs^i, \bw}$ is cut out by the induced section $s'$ of $\CE'$, we need to establish the fact that $\CE'$ is a vector bundle, or in other words that the map
$$
C_0 \xrightarrow{(g,h)} \bigoplus_{j \in I} \text{Hom}(q\tCV_j,\CV_j) \oplus q^{-1}\CO
$$
is surjective. This is equivalent to the point-wise injectivity of the dual map
\begin{multline*}
\bigoplus_{j \in I} \text{Hom}(V_j,\tV_j) \oplus \BC \longrightarrow \\ \bigoplus_{e = \overrightarrow{jj'} \in E} \left( \Hom(V_{j'},\tV_j) \oplus \Hom(V_j,\tV_{j'}) \right)  \bigoplus_{j \in I} \left( \Hom(V_j,W_j) \oplus \Hom(W_j,\tV_j) \right)
\end{multline*}
$$
\left( \{\tau_j\}_{j \in I}, c \right) \mapsto \left(\tY_{\overrightarrow{jj'}} \tau_{j'} - \tau_j Y_{\overrightarrow{jj'}}, \tau_{j'} X_{\overrightarrow{jj'}} - \tX_{\overrightarrow{jj'}} \tau_j, c B_j - \tB_j \tau_j, \tau_jA_j - c\tA_j \right)
$$
If a certain collection $(\{\tau_j\}_{j \in I}, c)$ lied in the kernel of the map above, then 

\begin{itemize}[leftmargin=*]

\item if $c \neq 0$, then the collection $\{\text{Im } \tau_j\}_{j \in I}$ is invariant under the $X,Y$ maps and contains $\text{Im }\tA_j$, which implies that all the maps $\tau_j$ are surjective; however, this is impossible due to the fact that $\dim V_i = \dim \tV_i - 1$

\medskip

\item if $c = 0$, then $\{\text{Ker }\tau_j\}_{j \in I}$ is invariant under the $X,Y$ maps and contains $\text{Im }A_j$, which implies that all the maps $\tau_j$ are 0

\end{itemize}

\medskip

\noindent The upshot of the present Subsection is that the structure sheaf of 
$$
\CN_{\bv,\bv+\bs^i,\bw} \hookrightarrow \CM_{\bv,\bw} \times \CM_{\bv+\bs^i,\bw}
$$
is equal (in $K$-theory) to the Koszul complex of the vector bundle $\CE'$
\begin{multline}
\label{eqn:formula simple correspondence}
\left[ \CO_{\CN_{\bv,\bv+\bs^i,\bw}} \right] = \\ \wedge^\bullet \left( \sum_{j \in I} \left( \frac {W_j}{\CV_j} + \frac {q\tCV_j}{W_j} \right)  + \sum_{e = \overrightarrow{jj'} \in E} \left( \frac {t_e \tCV_j}{\CV_{j'}}  + \frac {q \tCV_{j'}}{t_e \CV_j} \right) - \left(1+q\right) \sum_{j \in I} \frac {\tCV_j}{\CV_j} - q\right)
\end{multline}

\bigskip

\section{$K$-theoretic stable envelopes}
\label{sec:stable}

\medskip

\subsection{} We now recall the construction of \cite{AO, MO, O1, O2, OS}, which produces geometrically defined quantum groups. Explicitly, \loccit constructs \textbf{stable envelopes}
\begin{equation}
\label{eqn:stable basis}
K(\bw') \otimes K(\bw'') \xrightarrow{\stab, \stab'} K(\bw)
\end{equation}
for any $\bw = \bw'+\bw''$, which are uniquely determined by the 

\medskip

\begin{itemize}

\item support \eqref{eqn:support} and \eqref{eqn:support prime}

\medskip

\item normalization \eqref{eqn:normalization} and \eqref{eqn:normalization prime}

\medskip

\item degree \eqref{eqn:degree} and \eqref{eqn:degree prime}

\end{itemize}

\medskip

\noindent conditions (\cite{OS}). In the language of \loccitt, the map $\stab$ (respectively $\stab'$) corresponds to slope infinitesimally smaller (respectively larger) than 0 and the natural (respectively opposite) polarization of Nakajima quiver varieties. In the following Subsections, we will show how to define the map \eqref{eqn:stable basis} more generally, in the context of a symplectic variety $X$ endowed with a rank 1 torus $A \cong \BC^*$ action
\begin{equation}
\label{eqn:symplectic action}
A  \curvearrowright X
\end{equation}
which preserves the symplectic form. In this context, \loccit construct maps
\begin{equation}
\label{eqn:stable basis general}
K_A(X^A) \xrightarrow{\stab, \stab'} K_A(X)
\end{equation}
When $X$ is $\CM_{\bv,\bw}$ and one chooses a very specific $A \cong \BC^*$ action (see Subsection \ref{sub:apply}), the maps \eqref{eqn:stable basis general} will give rise to the sought-for maps \eqref{eqn:stable basis}.

\medskip

\subsection{}

The construction of stable envelopes, on which we will now elaborate, also applies to $K$-theory which is equivariant with respect to a bigger torus $T \supset A$. We assume that this bigger torus scales the symplectic form on $X$ by a certain character $q:T \rightarrow \BC^*$. Then the symplectic form gives rise to an isomorphism 
\begin{equation}
\label{eqn:pairing 1}
\Tan \ X  \cong q^{-1} (\Tan \ X)^\vee
\end{equation}
of $T$-equivariant vector bundles. Let us consider a rank 1 torus $A$ which preserves the symplectic form, and consider any connected component of the fixed point locus
$$
F \subset X^A
$$
We may consider the decomposition of the normal bundle
$$
\Nor_{F}X = \Nor^+_{F}X \oplus \Nor^-_{F}X
$$
into sub-bundles on which $A$ acts with positive/negative weight, respectively. The isomorphism \eqref{eqn:pairing 1} yields a perfect pairing between the positive and negative normal sub-bundles, i.e.
\begin{equation}
\label{eqn:pairing 2}
\Nor^+_{F}X \cong q^{-1} (\Nor^{-}_{F}X)^\vee
\end{equation}

\medskip

\begin{definition}
\label{def:polarization}

Assume we are given a decomposition of the ($K$-theory class of the) tangent bundle of $X$ as
\begin{equation}
\label{eqn:halves}
[\eTan \ X] = [\eTan^{\frac 12}X] + q^{-1}[\eTan^{\frac 12}X]^\vee
\end{equation}
Then for every fixed component $F \subset X^A$, we may consider the decomposition
$$
[\eTan^{\frac 12} X]\Big|_F = [\eTan^{\frac 12,+} X] \Big|_F + [\eTan^{\frac 12, 0} X] \Big|_F + [\eTan^{\frac 12, -} X] \Big|_F
$$
into positive, zero and negative weights for the $A \cong \BC^*$ action. Let 
$$
[\eTan^{\frac 12, \neq 0} X] \Big|_F = [\eTan^{\frac 12, +} X] \Big|_F + [\eTan^{\frac 12, -} X] \Big|_F
$$
There exist well-defined line bundle classes on $F$ (times $\pm 1$) given by the formulas
\begin{align}
&\varepsilon_F = (-1)^{\emph{rank } \eTan^{\frac 12,+} X|_F} \det \left( \eNor^-_{F}X - \eTan^{\frac 12,\neq 0} X \Big|_F  \right)^{\frac 12} \label{eqn:polarization} \\
&\varepsilon'_F = (-1)^{\emph{rank } \eTan^{\frac 12,-} X|_F} \det \left( \eNor^-_{F}X - q^{-1} (\eTan^{\frac 12,\neq 0} X)^\vee \Big|_F  \right)^{\frac 12} \label{eqn:polarization prime}
\end{align}
We will call $\varepsilon_F$ and $\varepsilon_{F}'$ \textbf{polarizations}, which are opposite to each other.

\end{definition}

\medskip

\noindent Let us work out polarizations explicitly, with respect to a decomposition of the ($K$-theory class of the) normal bundle into line bundles
\begin{equation}
\label{eqn:normal explicit}
[\Nor_F^- X] = \sum_{i=1}^k u_i + \sum_{j=1}^l v_j \quad \stackrel{\eqref{eqn:pairing 2}}{\Rightarrow} \quad [\Nor_F^+ X] = \sum_{i=1}^k \frac 1{qu_i} + \sum_{j=1}^l \frac 1{qv_j}
\end{equation}
Let us assume that the decomposition is chosen such that
\begin{equation}
\label{eqn:halves explicit}
[\Tan^{\frac 12,-}X] = \sum_{i=1}^k u_i \quad \Rightarrow \quad [\Tan^{\frac 12,+}X] = \sum_{j=1}^l \frac 1{qv_j}
\end{equation}
Then we have
\begin{align}
&\varepsilon_F = (-1)^l \det \left(\sum_{j=1}^l v_j - \sum_{j=1}^l  \frac 1{qv_j}\right)^{\frac 12} = \prod_{j=1}^l (-q^{\frac 12} v_j) \label{eqn:polarization explicit} \\
&\varepsilon_F' = (-1)^k \det \left(\sum_{i=1}^k u_i - \sum_{i=1}^k \frac 1{qu_i}\right)^{\frac 12} = \prod_{i=1}^k (-q^{\frac 12} u_i)  \label{eqn:polarization prime explicit}
\end{align}
and so $\varepsilon_F$ and $\varepsilon'_F$ are simply the determinants of $\Tan^{\frac 12,\pm}X|_F$ (up to a power of $q^{\frac 12}$ and a sign); thus, they are well-defined elements of the equivariant $K$-theory of $F$.

\medskip

\subsection{}

The following definition pertains to the setting of \eqref{eqn:symplectic action}. For $A \cong \BC^* \curvearrowright X$, we will write $t \cdot x \in X$ for the action of any $t \in \BC^*$ on any $x \in X$.

\medskip

\begin{definition}
\label{def:attracting}

The \textbf{attracting subvariety} is 
\begin{equation}
\label{eqn:attracting}
\eattr = \Big\{ (p,x) \text{ s.t. } \lim_{t \rightarrow 0} t \cdot x = p \Big\} \subset X^A \times X 
\end{equation}
For any connected component $F \subset X^A$, the attracting subvariety of $F$ is
\begin{equation}
\label{eqn:attracting to f}
\eattr_F = \eattr \cap (F \times X)
\end{equation}

\end{definition}

\medskip

\noindent In general, $\attr$ is not closed, which causes problems in $K$-theoretic computations. To remedy this problem, \cite{MO} enhances it by defining the \textbf{full attracting subvariety}
\begin{equation}
\label{eqn:full attracting}
\attr^f \subset X^A \times X 
\end{equation}
as the smallest closed $A$-invariant subset of $X^A \times X$ which contains $\text{Attr}$ and satisfies
$$
(p,p') \in \attr^f \text{ and } \lim_{t \rightarrow 0} t \cdot x = p' \quad \text{imply} \quad (p,x) \in \attr^f
$$
To convey the previous equation more intuitively, consider the following notion.

\medskip

\begin{definition}
\label{def:flow}

A closed (respectively half-open) \textbf{flow line} is a copy of $\BP^1 \subset X$ (respectively $\BA^1 \subset X$) which is the closure of the $A \cong \BC^*$ orbit of a single $x \in X$. Flow lines will be considered to be oriented toward $\lim_{t\rightarrow 0} t \cdot x$.

\end{definition}

\medskip

\noindent The full attracting subvariety consists of those pairs $(p,x) \in X^A \times X$ such that one can travel from $x$ to $p$ along a half-open flow line followed by a chain of closed flow lines. As shown in \cite{MO}, there is a partial order on the fixed components of $X^A$ with
\begin{equation}
\label{eqn:partial order}
F' > F
\end{equation}
if there is a chain of closed flow lines from $F$ into $F'$. By analogy with \eqref{eqn:attracting to f}, let
\begin{equation}
\label{eqn:full attracting to f}
\attr^f_F = \attr^f \cap (F \times X)
\end{equation}
for any connected component $F \subset X^A$. 

\medskip

\subsection{}

With the constructions of the preceding subsection in mind, the \textbf{support condition} defining the stable envelopes \eqref{eqn:stable basis} is that they may be represented as
\begin{align}
&\stab = \text{proj}_{2*} \left(\Gamma \cdot \text{proj}_1^*\right) \label{eqn:support} \\
&\stab' = \text{proj}_{2*} \left(\Gamma' \cdot \text{proj}_1^*\right)\label{eqn:support prime}
\end{align}
for some classes $\Gamma, \Gamma' \in K(X^A \times X)$ supported on the full attracting subvariety, where $\text{proj}_{1,2} : X^A \times X \rightarrow X^A,X$ denote the standard projections. 

\medskip

\noindent Intuitively speaking, in a neighborhood of $X^A \times X^A$, the attracting subvariety ``looks like" the total space of the attracting part of the normal bundle of $X$
$$
\text{Tot} (\Nor^+_{X^A} X )
$$
The \textbf{normalization condition} says that the restriction of stable envelopes to $X^A \times X^A$ is the one prescribed by the aforementioned intuition, times the polarization line bundles. Specifically, for any connected component $F \subset X^A$, we require that
\begin{align}
&\Gamma \Big|_{F \times F} = \Delta_{F*} \left( \varepsilon_F \cdot \wedge^\bullet  (\Nor^{-}_F X)^\vee \right) \label{eqn:normalization} \\ 
&\Gamma' \Big|_{F \times F} = \Delta_{F*} \left( \varepsilon'_F \cdot \wedge^\bullet  (\Nor^{-}_F X)^\vee \right) \label{eqn:normalization prime} 
\end{align}
With respect to the notation in \eqref{eqn:normal explicit} and \eqref{eqn:halves explicit}, the fact that the polarization line bundles are given by \eqref{eqn:polarization explicit} and \eqref{eqn:polarization prime explicit} implies that
\begin{align}
&\Gamma \Big|_{F \times F} = \Delta_{F*} \left(q^{\frac l2} \prod_{i=1}^k \left(1-\frac 1{u_i}\right) \prod_{j=1}^l (1-v_j) \right) \label{eqn:normalization explicit} \\ 
&\Gamma' \Big|_{F \times F} = \Delta_{F*} \left(q^{\frac k2} \prod_{i=1}^k \left(1-u_i\right) \prod_{j=1}^l \left(1-\frac 1{v_j} \right) \right) \label{eqn:normalization prime explicit} 
\end{align}

\medskip

\subsection{}

Finally, the degree condition has to do with constraints on the restriction of $\Gamma, \Gamma'$ to $F' \times F$, where $F \neq F'$ are different connected components of $X^A$. Recall the partial order on the connected components of $X^A$ in \eqref{eqn:partial order}. Because $\Gamma$ and $\Gamma'$ are supported on the full attracting subvariety, then
$$
\Gamma\Big|_{F' \times F} = \Gamma'\Big|_{F' \times F} = 0
$$
if $F' \not \geq F$. If $F = F'$, then the restrictions above are governed by \eqref{eqn:normalization}--\eqref{eqn:normalization prime}. If however $F' > F$, then the restriction $\Gamma |_{F' \times F}$ is an $A$-equivariant $K$-theory class on the variety $F' \times F$. However, as $F' \times F$ is point-wise fixed by $A$, we may decompose $\Gamma |_{F' \times F}$ into its $A$-isotypic components 
\begin{equation}
\label{eqn:restriction degree}
\Gamma \Big|_{F' \times F} = \sum_{k \in \BZ} a^k \cdot \gamma_k
\end{equation}
where $a$ denotes the standard character of $A \cong \BC^*$, and each $\gamma_k$ is an $A$-invariant class on $F' \times F$ (if a torus acts trivially on a variety, an invariant $K$-theory class on said variety refers to one which is endowed with the trivial torus action). Define the minimum (respectively maximum) degree as
\begin{align}
&\mindeg_A \left(\Gamma \Big|_{F' \times F} \right) = \text{smallest }k \text{ such that }\gamma_k \neq 0 \label{eqn:mindeg} \\
&\maxdeg_A \left(\Gamma \Big|_{F' \times F} \right) = \text{largest }k \text{ such that }\gamma_k \neq 0 \label{eqn:maxdeg} 
\end{align}
The \textbf{degree condition} states that the minimum/maximum degrees of the restriction of $\Gamma$ to $F' \times F$ are bounded by those of $\Gamma$ to $F \times F$, in the following sense
\begin{multline}
\label{eqn:degree}
\mindeg_A \left(\Gamma \Big|_{F \times F} \right) \leq \mindeg_A \left(\Gamma \Big|_{F' \times F} \right) \\  \maxdeg_A \left(\Gamma \Big|_{F' \times F} \right) < \maxdeg_A \left(\Gamma \Big|_{F \times F} \right)
\end{multline}
For $\Gamma'$, we define the degree condition analogously, but switching the symbols $\leq$ and $<$ in the formula above
\begin{multline}
\label{eqn:degree prime}
\mindeg_A \left(\Gamma' \Big|_{F \times F} \right) < \mindeg_A \left(\Gamma' \Big|_{F' \times F} \right) \\  \maxdeg_A \left(\Gamma' \Big|_{F' \times F} \right) \leq \maxdeg_A \left(\Gamma' \Big|_{F \times F} \right)
\end{multline}
The left-most and right-most integer in \eqref{eqn:degree} and \eqref{eqn:degree prime} can be calculated from \eqref{eqn:normalization} and \eqref{eqn:normalization prime}, respectively. Explicitly, let us assume that we have the decomposition \eqref{eqn:normal explicit} and \eqref{eqn:halves explicit}, and that moreover the line bundles $u_i$ and $v_j$ have $A$-degree equal to $-1$ (tori $A \cong \BC^*$ with this property are called ``minuscule", and this will be the case throughout the present paper). Using \eqref{eqn:normalization explicit} and \eqref{eqn:normalization prime explicit}, we obtain
\begin{align}
&-l \leq \mindeg_A \left(\Gamma \Big|_{F' \times F} \right) \leq  \maxdeg_A \left(\Gamma \Big|_{F' \times F} \right) < k \label{eqn:degree explicit} \\
&-k < \mindeg_A \left(\Gamma' \Big|_{F' \times F} \right) \leq \maxdeg_A \left(\Gamma' \Big|_{F' \times F} \right) \leq l \label{eqn:degree prime explicit}
\end{align}

\medskip

\subsection{}
\label{sub:apply}

Let us discuss the particular case of the setup above when
$$
X = \CM_{\bv,\bw}
$$
and
\begin{equation}
\label{eqn:a}
A : \BC^* \hookrightarrow G_{\bw} \curvearrowright \CM_{\bv,\bw}, \qquad A(t) = \prod_{i \in I} \begin{bmatrix}
t I_{w_i'} & 0 \\ 0 & I_{w_i''} \end{bmatrix}
\end{equation}
corresponds to a fixed decomposition $W_\bullet = W'_\bullet \oplus W''_\bullet$ of the framing vector spaces into subspaces of dimensions $\bw'$ and $\bw''$, respectively, where $\bw'+\bw'' = \bw$. The fixed point set $X^A$ is given by those $\bv$-dimensional framed double quiver representations $V_{\bullet}$ which split up as a direct sum
\begin{equation}
\label{eqn:direct sum}
V_{\bullet} \cong V'_\bullet \oplus V''_\bullet
\end{equation}
with the $A_i,B_i$ maps only being non-zero between $W'_\bullet \leftrightharpoons V'_\bullet$ and $W''_\bullet \leftrightharpoons V''_\bullet$. Thus, we obtain a natural isomorphism
\begin{equation}
\label{eqn:cb}
X^A \cong \bigsqcup_{\bv = \bv'+\bv''} \CM_{\bv',\bw'} \times \CM_{\bv'',\bw''}
\end{equation}
and we will schematically indicate points of $X^A$ as pairs $(V_\bullet', V_\bullet'')$ (although the framing maps $A_i,B_i$ are also part of the data). Note that \eqref{eqn:cb} is a decomposition of $X^A$ into connected components, since Nakajima quiver varieties are themselves connected (see \cite{cb}). The following result is well-known, and we leave it to the reader to derive it from the explicit description of $\lim_{t \rightarrow 0} t \cdot V_{\bullet}$ for a point $V_{\bullet} \in \CM_{\bv,\bw}$. 

\medskip

\begin{lemma}
\label{lem:attr}

The attracting subvariety $\eattr$ parameterizes triples of framed double quiver representations $(V'_\bullet, V''_\bullet, V_\bullet) \in \CM_{\bv',\bw'} \times \CM_{\bv'',\bw''} \times \CM_{\bv,\bw}$ which fit into a short exact sequence
\begin{equation}
\label{eqn:attr}
V'_\bullet \hookrightarrow V_\bullet  \twoheadrightarrow V''_\bullet 
\end{equation}
The arrows in \eqref{eqn:attr} are required to commute with the $X,Y$ maps, and also to commute with the $A,B$ maps via the split short exact sequence
$$
W'_\bullet \hookrightarrow W_\bullet \twoheadrightarrow W''_\bullet
$$
induced by $W_\bullet = W'_\bullet \oplus W''_\bullet$. 

\end{lemma}

\medskip

\noindent Note that if it exists, the short exact sequence in \eqref{eqn:attr} is uniquely determined (up to automorphism), since $V_\bullet'$ is forced to match the subrepresentation of $V_\bullet$ generated by the $X,Y$ maps acting on $\{\text{Im }A_i|_{W'_i}\}_{i \in I}$. Once the subrepresentation $V_{\bullet}'$ is determined, then $V''_{\bullet}$ is also determined as the corresponding quotient. 

\medskip

\begin{lemma}
\label{lem:flow line}
There exists a closed flow line (as in Definition \ref{def:flow})
$$
\text{from } (\tV'_\bullet, \tV''_\bullet) \text{ to }(V'_\bullet, V''_\bullet)
$$ 
if and only if there exist surjective maps of framed double quiver representations
$$
\xymatrix{V'_{\bullet} \ar@{->>}^{\pi'_\bullet}[rr] \ar@/^/[dr]  & & {\tV'_{\bullet}} \ar@/^/[dl] \\ & W'_{\bullet} \ar@/^/[ul] \ar@/^/[ur] &} \qquad \text{and} \qquad \xymatrix{{\tV_{\bullet}''} \ar@{->>}^{\pi''_\bullet}[rr] \ar@/^/[dr]  & & {V''_{\bullet}} \ar@/^/[dl] \\ & W''_{\bullet} \ar@/^/[ul] \ar@/^/[ur] &}
$$
such that the kernels of $\pi'_\bullet$ and $\pi''_\bullet$ are isomorphic double quiver representations. In other words, the latter condition requires the points
$$
(\tV'_\bullet, V'_\bullet) \in \CN_{\bv'-\bn, \bv', \bw'} \qquad \text{and} \qquad (V''_\bullet, \tV''_\bullet) \in \CN_{\bv'', \bv''+\bn, \bw''}
$$
to have the same projection to $\fY_{\bn}$ via the maps in diagram \eqref{eqn:ext diagram 2}.

\end{lemma}

\medskip

\begin{proof} Consider any $V_\bullet \in \CM_{\bv,\bw}$ such that $\lim_{t \rightarrow 0} t \cdot V_\bullet$ and $\lim_{t \rightarrow \infty} t \cdot V_\bullet$ exist in $\CM_{\bv,\bw}$, and let 
$$
V'_\bullet \subset V_\bullet \quad \text{and} \quad \tV''_\bullet \subset V_\bullet
$$
denote the subspaces generated by the $X,Y$ maps acting on $\{\text{Im }A_i|_{W_i'}\}_{i\in I}$ and $\{\text{Im }A_i|_{W_i''}\}_{i\in I}$, respectively. Let us write
$$
K_\bullet = V'_\bullet \cap \tV''_\bullet
$$
and consider decompositions as vector spaces
$$
V'_\bullet = K_\bullet \oplus \tV'_\bullet \quad \text{and} \quad \tV''_\bullet = K_\bullet \oplus V''_\bullet
$$
With respect to these decompositions, the maps $(X_e,Y_e,A_i,B_i)$ take the form

\begin{equation}
\label{eqn:maps and ts 1}
\xymatrix{{\tV'_{\bullet}} \ar@/^/[ddr] \ar[rr] \ar@(ul,ur) \ar@{.>}[ddrrr] & & K_{\bullet} \ar@(ul,ur) \ar@{.>}[ddr] \ar@{.>}[ddl] & & V''_{\bullet} \ar@(ul,ur) \ar@/^/[ddl] \ar[ll] \ar@{.>}[ddlll] \\ \\
& W'_{\bullet} \ar@/^/[uul] \ar@/_1pc/[uur] & & W''_{\bullet} \ar@/^/[uur] \ar@/^1pc/[uul] &}
\end{equation}
Because the torus $A \cong \BC^*$ scales the compositions
\begin{align*}
&\text{proj}_{W'_\bullet} B \circ \Big(\text{any product of }X \text{ and }Y\Big) \circ A\Big|_{W''_\bullet} \\\\
&\text{proj}_{W''_\bullet} B \circ \Big(\text{any product of }X \text{ and }Y\Big) \circ A\Big|_{W'_\bullet}
\end{align*}
by $t$ and $t^{-1}$, respectively, these compositions need to vanish if $\lim_{t \rightarrow 0} t \cdot V_\bullet$ and $\lim_{t \rightarrow \infty} t \cdot V_\bullet$ are to exist. In this case, the dotted arrows in the diagram \eqref{eqn:maps and ts 1} are all 0, so the representation $V_\bullet$ takes the form
\begin{equation}
\label{eqn:maps and ts 2}
\xymatrix{{\tV'_{\bullet}} \ar@/^/[ddr]^t \ar[rr]^1 \ar@(ul,ur)^1  & & K_{\bullet} \ar@(ul,ur)^1 & & V''_{\bullet} \ar@(ul,ur)^1 \ar@/^/[ddl]^1 \ar[ll]_1  \\ \\
& W'_{\bullet} \ar@/^/[uul]^{t^{-1}} \ar[uur]_{t^{-1}} & & W''_{\bullet} \ar@/^/[uur]^1 \ar[uul]^1 &}
\end{equation}
In diagram \eqref{eqn:maps and ts 2}, we displayed above every arrow its character with respect to the $A \cong \BC^*$ action of \eqref{eqn:a}. Of course, we cannot hope to obtain a stable quiver representation by directly sending $t \rightarrow 0$ or $t \rightarrow \infty$ in the diagram above, but we can first modify the quadruple $(X_e,Y_e,A_i,B_i)$ by the action of \eqref{eqn:gauge group}. Specifically, if we rescale $\tV'_{\bullet}$ and $K_\bullet$ by $t$, then the weights of the various arrows take the form
\begin{equation}
\label{eqn:maps and ts 3}
\xymatrix{{\tV'_{\bullet}} \ar@/^/[ddr]^1 \ar[rr]^1 \ar@(ul,ur)^1  & & K_{\bullet} \ar@(ul,ur)^1 & & V''_{\bullet} \ar@(ul,ur)^1 \ar@/^/[ddl]^1 \ar[ll]_t  \\ \\
& W'_{\bullet} \ar@/^/[uul]^1 \ar[uur]_1 & & W''_{\bullet} \ar@/^/[uur]^1 \ar[uul]^t &}
\end{equation}
and so $\lim_{t\rightarrow 0} t \cdot V_\bullet$ equals $V'_\bullet \oplus V''_\bullet$ as framed double quiver representations. 

\medskip

\noindent On the other hand, if we rescale $\tV'_{\bullet}$ in \eqref{eqn:maps and ts 2} by $t$, then the weights of the various arrows take the form
\begin{equation}
\label{eqn:maps and ts 4}
\xymatrix{{\tV'_{\bullet}} \ar@/^/[ddr]^1 \ar[rr]^{t^{-1}} \ar@(ul,ur)^1  & & K_{\bullet} \ar@(ul,ur)^1 & & V''_{\bullet} \ar@(ul,ur)^1 \ar@/^/[ddl]^1 \ar[ll]_1  \\ \\
& W'_{\bullet} \ar@/^/[uul]^1 \ar[uur]_{t^{-1}} & & W''_{\bullet} \ar@/^/[uur]^1 \ar[uul]^1 &}
\end{equation}
and so $\lim_{t\rightarrow \infty} t \cdot V_\bullet$ equals $\tV'_\bullet \oplus \tV''_\bullet$ as framed double quiver representations. The construction above naturally produces short exact sequences of framed double quiver representations
$$
0 \rightarrow K_\bullet \rightarrow V'_\bullet \xrightarrow{\pi'_\bullet} \tV'_\bullet \rightarrow 0 \quad \text{and} \quad 0 \rightarrow K_\bullet \rightarrow \tV''_\bullet \xrightarrow{\pi''_\bullet} V''_\bullet \rightarrow 0
$$
which concludes the proof of Lemma \ref{lem:flow line}. \end{proof}

\medskip

\noindent Lemma \ref{lem:flow line} implies that the partial order \eqref{eqn:partial order} in the situation at hand is
\begin{equation}
\label{eqn:nakajima partial order}
\CM_{\bv',\bw'} \times \CM_{\bv'',\bw''} \geq  \CM_{\bv'-\bn,\bw'} \times \CM_{\bv''+\bn,\bw''}
\end{equation}
for any $\bn \in \nn$, because flowing from a point $\tV_\bullet' \oplus \tV_\bullet''$ to a point $V_\bullet' \oplus V_\bullet''$ is possible only if $V_\bullet'$ is obtained from $\tV_\bullet'$ by adding (respectively $V_\bullet''$ is obtained from $\tV_\bullet''$ by removing) a quiver representation of arbitrary dimension $\bn \in \nn$. As a consequence of Lemmas \ref{lem:attr} and \ref{lem:flow line}, we obtain the following.

\medskip

\begin{corollary}
\label{cor:full attr}

The full attracting subvariety $\eattr^f$ parameterizes triples of framed double quiver representations $(V'_\bullet, V''_\bullet, V_\bullet) \in \CM_{\bv',\bw'} \times \CM_{\bv'',\bw''} \times \CM_{\bv,\bw}$ such that there exist linear maps
\begin{equation}
\label{eqn:full attr}
V'_\bullet \xrightarrow{f} V_\bullet  \xrightarrow{g} V''_\bullet 
\end{equation}
such that the following conditions hold

\medskip

\begin{itemize}[leftmargin=*]

\item The composition $g \circ f$ is $0$ \\

\item The maps $f$ and $g$ commute with the $X,Y$ maps, and also commute with the $A,B$ maps via the split short exact sequence
$$
W'_\bullet \hookrightarrow W_\bullet \twoheadrightarrow W''_\bullet
$$

\item Letting $\tV'_\bullet = \emph{Im }f$ and $\tV''_\bullet = V_{\bullet}/\emph{Im } f$, we require the existence of filtrations
\begin{align*}
&V'_\bullet = E_\bullet^0 \twoheadrightarrow E_\bullet^1 \twoheadrightarrow \dots \twoheadrightarrow E_\bullet^{k-1} \twoheadrightarrow E_\bullet^k = \tV'_\bullet \\
&\tV''_\bullet = F_\bullet^k \twoheadrightarrow F_\bullet^{k-1} \twoheadrightarrow \dots \twoheadrightarrow F_\bullet^1 \twoheadrightarrow F_\bullet^0 = V''_\bullet
\end{align*}
such that the kernels of the maps $E_\bullet^l \twoheadrightarrow E_\bullet^{l+1}$ and $F_\bullet^{l+1} \twoheadrightarrow F_\bullet^l$ are isomorphic, $\forall l$.

\end{itemize}

\end{corollary}

\medskip

\noindent Note that the stability condition forces the map $g$ in \eqref{eqn:full attr} to be surjective.

\medskip

\subsection{} Corollary \ref{cor:full attr} describes the full attracting subvariety which supports the stable envelopes \eqref{eqn:stable basis}. Let us now spell out the normalization and degree conditions. Our starting point is formula \eqref{eqn:tangent bundle} for the tangent bundle to $\CM_{\bv,\bw}$, which implies the following formula for the normal bundle
\begin{multline}
\label{eqn:normal bundle}
\Nor_{\CM_{\bv',\bw'} \times \CM_{\bv'',\bw''}}(\CM_{\bv,\bw}) = \sum_{e = \oij \in E} \left(\underbrace{\frac {\CV'_j}{t_e\CV''_i} + \frac {t_e\CV'_i}{q\CV''_j}}_{+1} + \underbrace{\frac {\CV''_j}{t_e\CV'_i} + \frac {t_e\CV''_i}{q\CV'_j}}_{-1}\right) - \\ - \sum_{i \in I} \left(1+\frac 1q\right)\left(\underbrace{\frac {\CV'_i}{\CV''_i}}_{+1} + \underbrace{\frac {\CV''_i}{\CV'_i}}_{-1} \right) + \sum_{i \in I} \left( \underbrace{\frac {\CV'_i}{W''_i} + \frac {W'_i}{q\CV''_i}}_{+1} + \underbrace{\frac {\CV''_i}{W'_i} + \frac {W''_i}{q\CV'_i}}_{-1} \right)
\end{multline}
(the fact that points of the fixed locus decompose as in \eqref{eqn:direct sum} means that we have natural isomorphisms $\CV_\bullet = \CV_{\bullet}' \oplus \CV_{\bullet}''$ on the fixed locus $\CM_{\bv',\bw'} \times \CM_{\bv'',\bw''}$). The underbraces in \eqref{eqn:normal bundle} indicate the weights of the $A \cong \BC^*$ action in the various summands, hence the sub-bundles $\Nor^+$ and $\Nor^-$ of $\Nor$ are obtained by only keeping the terms with $+1$ and $-1$ in the underbraces, respectively. 

\medskip

\noindent To compute the polarization line bundles, we must fix a choice of decomposition \eqref{eqn:halves}. The one we choose is the ``canonical" one for cotangent bundles, i.e. we only retain the terms that do not involve $q$ in \eqref{eqn:tangent bundle}
$$
[\Tan^{\frac 12} \CM_{\bv,\bw} ] = \sum_{e = \oij \in E} \frac { \CV_j}{t_e\CV_i}  - \sum_{i \in I} \frac {\CV_i}{\CV_i} + \sum_{i \in I} \frac {\CV_i}{W_i} 
$$
With respect to this choice, the terms in formula \eqref{eqn:normal explicit} are
\begin{align}
&\sum_{i=1}^k u_i = \sum_{e = \oij \in E} \frac {\CV''_j}{t_e\CV'_i} -  \sum_{i \in I} \frac {\CV''_i}{\CV'_i} + \sum_{i \in I}  \frac {\CV''_i}{W'_i} \label{eqn:u nakajima} \\
&\sum_{j=1}^l v_j = \sum_{e = \oij \in E} \frac {t_e\CV''_i}{q\CV'_j} -  \sum_{i \in I} \frac {\CV''_i}{q\CV'_i} + \sum_{i \in I} \frac {W''_i}{q\CV'_i} \label{eqn:v nakajima}
\end{align}
In particular, we have
\begin{align}
&k = \bw' \cdot \bv'' - \langle \bv', \bv'' \rangle \label{eqn:k nakajima} \\
&l = \bw'' \cdot \bv' - \langle \bv'',\bv' \rangle \label{eqn:l nakajima}
\end{align}
With this in mind, conditions \eqref{eqn:normalization explicit} and \eqref{eqn:normalization prime explicit} for $F = \CM_{\bv',\bw'} \times \CM_{\bv'', \bw''}$ read
\begin{equation}
\label{eqn:alpha alpha'}
\Gamma \Big|_{F \times F} = \Delta_{F*}(\alpha) \qquad \text{and} \qquad \Gamma' \Big|_{F \times F} = \Delta_{F*}(\alpha')
\end{equation}
where
\begin{align*}
&\alpha = q^{\frac {\bw'' \cdot \bv' - \langle \bv'',\bv' \rangle}2} \wedge^\bullet \left(\sum_{e = \oij \in E} \left( \frac {t_e\CV'_i}{\CV''_j} +  \frac {t_e\CV''_i}{q\CV'_j} \right) - \sum_{i \in I} \left( \frac {\CV'_i}{\CV''_i} +  \frac {\CV''_i}{q\CV'_i} \right) + \sum_{i \in I} \left( \frac {W'_i}{\CV''_i} + \frac {W''_i}{q\CV'_i}\right) \right) \\ 
&\alpha' = q^{\frac {\bw' \cdot \bv'' - \langle \bv', \bv'' \rangle}2} \wedge^\bullet \left(  \sum_{e = \oij \in E} \left( \frac {\CV''_j}{t_e\CV'_i} + \frac {q\CV'_j}{t_e\CV''_i} \right) - \sum_{i \in I} \left(\frac {\CV''_i}{\CV'_i} + \frac {q\CV'_i}{\CV''_i} \right) + \sum_{i \in I} \left( \frac {\CV''_i}{W'_i} + \frac {q\CV'_i}{W''_i} \right) \right)
\end{align*}

\medskip

\subsection{}

Throughout the remainder of the present Section, the symbol $\otimes$ with no subscript will refer to tensor product of $\BF$-vector spaces. Using the stable envelopes \eqref{eqn:stable basis}, the groundbreaking works \cite{AO,MO,O1,O2,OS} consider the composition
\begin{equation}
\label{eqn:r-matrix}
R_{\bw^1,\bw^2} : K(\bw^1) \otimes K(\bw^2) \xrightarrow{\stab} K(\bw^1+\bw^2)  \xrightarrow{{\stab'}^T} K(\bw^1) \otimes K(\bw^2)
\end{equation}
for any $\bw^1,\bw^2 \in \nn$, which satisfies the equality
\begin{equation}
\label{eqn:qybe}
R_{\bw^1,\bw^2} R_{\bw^1,\bw^3} R_{\bw^2,\bw^3} = R_{\bw^2,\bw^3} R_{\bw^1,\bw^3} R_{\bw^1,\bw^2}
\end{equation}
for any $\bw^1,\bw^2,\bw^3 \in \nn$ (one interprets \eqref{eqn:qybe} as an equality of endomorphisms of $K(\bw^1) \otimes K(\bw^2) \otimes K(\bw^3)$, with each $R$-matrix acting in the tensor product of exactly two of the factors, as indicated by the indices in question). Formula \eqref{eqn:qybe} is the famous \textbf{quantum Yang-Baxter equation}.

\medskip

\begin{definition}
\label{def:quantum group}

(\cite{OS}) Consider any $\bw,\bw^{\emph{aux}} \in \nn$ and any matrix coefficient $c$ of 
\begin{equation}
\label{eqn:r-matrix coefficient}
K(\bw^{\emph{aux}}) \otimes K(\bw) \xrightarrow{R_{\bw^{\emph{aux}},\bw}} K(\bw^{\emph{aux}}) \otimes K(\bw)
\end{equation}
in the first tensor factor \footnote{Up to linear combinations, such a matrix coefficient is determined by $e \in K(\bw^{\text{aux}})$ and $f \in K(\bw^{\text{aux}})^\vee$, and is given by composing \eqref{eqn:r-matrix coefficient} with $e \otimes \text{Id}$ on the left and $f \otimes \text{Id}$ on the right.}. As an endomorphism of $K(\bw)$, the matrix coefficient $c$ is a rational function in the equivariant parameters $\{u_{i1},\dots,u_{iw^{\emph{aux}}_i}\}_{i\in I}$ corresponding to the framing $\bw^{\emph{aux}}$. Then we consider the subalgebra
\begin{equation}
\label{eqn:quantum group}
\UU \subset \prod_{\bw \in \nn} \emph{End}(K(\bw))
\end{equation}
generated by all power series coefficients of $c$ as above in the variables $\{u_{i1},\dots,u_{iw^{\emph{aux}}_i}\}$ as well as the operators of multiplication by
\begin{equation}
\label{eqn:central operators of multiplication}
q^{\pm \frac {w_i}2} \qquad \text{and}\qquad p_d \left( W_i(1-q^{-1}) \right)
\end{equation}
for all $i \in I$ and $d \in \BZ \backslash 0$. The latter operators are clearly central elements.

\end{definition}

\medskip

\begin{proposition}
\label{prop:maps to subalgebra}

The map $\iota$ of \eqref{eqn:action hom} lands in the subalgebra $\UU$, i.e. factors through an $\BF$-algebra homomorphism
\begin{equation}
\label{eqn:iota prime}
\CA \rightarrow \UU
\end{equation}
Moreover, by Proposition \ref{prop:inj}, the map \eqref{eqn:iota prime} is injective.

\end{proposition}

\medskip

\begin{proof} We must show that the generators $e_{i,d}, f_{i,d}, a_{i,\pm d}, q^{\pm \frac {v_i}2}, b_{i,\pm d}, q^{\pm \frac {w_i}2}$ (see Theorem \ref{thm:generate} or Proposition \ref{prop:generate}) of $\CA$ lie in $\UU$. As the latter two generators are elements of $\UU$ due to Definition \ref{def:quantum group}, it remains to show that the former four generators can be realized as matrix coefficients of the $R$-matrix \eqref{eqn:r-matrix coefficient} in the first tensor factor. Let $\bw^{\text{aux}} = \bs^i$ for some fixed $i \in I$, and we will denote $R_{\bs^i, \bw}$ by 
\begin{equation}
\label{eqn:r(u)}
R(u) : K(\bs^i) \otimes K(\bw) \xrightarrow{R_{\bs^i, \bw}} K(\bs^i) \otimes K(\bw)
\end{equation}
in order to emphasize the dependence of $R(u)$ on the unique equivariant parameter $u = u_{i1}$ associated to the framing vector $\bs^i$. For any $\bx, \by \in \nn$, consider the various graded components of \eqref{eqn:r(u)} in the first tensor factor, namely
\begin{equation}
\label{eqn:r(u) coefficient}
_{\langle \bx |}R(u)_{|\by \rangle} : K(\bs^i)_{\by} \otimes K(\bw) \longrightarrow K(\bs^i)_{\bx} \otimes K(\bw)
\end{equation}
We will show that $e_{i,d}, f_{i,d}$ and $a_{i,\pm d}, q^{\pm \frac {v_i}2}$ can be obtained as various coefficients of the powers of $u$ in the following expressions, respectively
\begin{equation}
\label{eqn:r(u) coefficient small}
_{\langle \b0 |}R(u)_{|\bs^i \rangle}, \quad _{\langle \bs^i |}R(u)_{|\b0 \rangle} \quad \text{and} \quad _{\langle \b0 |}R(u)_{|\b0 \rangle}
\end{equation}
Let us first deal with $a_{i,\pm d}$ and $q^{\pm \frac {v_i}2}$. Fix any $\bv \in \nn$, and consider the following component of the fixed locus of the torus $A \cong \BC^*$ of \eqref{eqn:a}
$$
F = \CM_{\b0, \bs^i}  \times \CM_{\bv,\bw} \hookrightarrow \CM_{\bv,\bw+\bs^i}
$$
which is minimal with respect to the order \eqref{eqn:nakajima partial order}. Therefore, we have
\begin{equation}
\label{eqn:r(a) comp 1}
_{\langle \b0 |}R(u)_{|\b0 \rangle} = \frac {\Gamma|_{F \times F} \cdot \Gamma'|_{F \times F}}{\wedge^\bullet (\Tan \ \CM_{\bv,\bw+\bs^i}|_F)^\vee} =\Delta_{F*} \left(\frac {\alpha \cdot \alpha'}{\wedge^\bullet ( \Nor_F \CM_{\bv,\bw+\bs^i})^\vee} \right) 
\end{equation}
as a correspondence on $F \times F$, where $\alpha,\alpha'$ are the classes in \eqref{eqn:alpha alpha'}. In the particular case of the fixed component $F$, these classes are given by
$$
\alpha = \wedge^\bullet \left(  \frac u{\CV_i} \right) \qquad \text{and} \qquad \alpha' = q^{\frac {v_i}2} \wedge^\bullet \left( \frac {\CV_i}u \right)
$$
where $\CV_i$ is the $i$-th tautological vector bundle on $\CM_{\bv,\bw}$. Meanwhile, the exterior class of the normal bundle \eqref{eqn:normal bundle} is given by 
$$
\wedge^\bullet ( \Nor_F \CM_{\bv,\bw+\bs^i})^\vee = \wedge^\bullet \left(\frac u{\CV_i} + \frac {q\CV_i}u \right)
$$
Therefore, formula \eqref{eqn:r(a) comp 1} reads
\begin{equation}
\label{eqn:r(a) comp 2}
_{\langle \b0 |}R(u)_{|\b0 \rangle} = \text{multiplication by } q^{\frac {v_i}2} \wedge^\bullet\left(\frac {\CV_i(1-q)}u \right)
\end{equation}
By extracting the coefficient of $u^0$ in the expression above, expanded as $u \rightarrow \infty$ and $u \rightarrow 0$, we obtain the operators $q^{\frac {v_i}2}$ and $q^{-\frac {v_i}2}$, respectively. Meanwhile, we have
\begin{align*}
&\log \frac {\wedge^\bullet\left(\frac {\CV_i}u \right)}{\wedge^\bullet\left(\frac {\CV_i q}u \right)} = \sum_{d = 1}^{\infty} \frac {p_d(\CV_i(q-1))}{d u^d} \\
&\log \frac {\wedge^\bullet\left(\frac u{\CV_i} \right)}{\wedge^\bullet\left(\frac u{\CV_i q} \right)} = \sum_{d = 1}^{\infty} \frac {p_{-d}(\CV_i(q-1))}{d u^{-d}} 
\end{align*}
and so we may obtain $a_{i,\pm d} = p_{\pm d}(\CV_i(1-q^{-1}))$ for any $d > 0$ as a polynomial in the various power series coefficients of \eqref{eqn:r(a) comp 2}, as we needed to show. 

\medskip

\noindent We will now show that $e_{i,d}$ and $f_{i,d}$ can be obtained from \eqref{eqn:r(u) coefficient small}. To this end, fix any $\bv \in \nn$ and let us consider the following components of the fixed locus 
\begin{equation}
\label{eqn:f and f prime}
F = \CM_{\b0,\bs^i} \times \CM_{\bv+\bs^i,\bw} \hookrightarrow \CM_{\bv+\bs^i,\bw+\bs^i} \hookleftarrow \CM_{\bs^i,\bs^i} \times \CM_{\bv,\bw} = F'
\end{equation}
With respect to the partial order \eqref{eqn:nakajima partial order}, the component $F$ is minimal, and it is the only component strictly smaller than $F'$. We will represent points of $F$ and $F'$ as
$$
(0, \tV_\bullet) \qquad \text{and} \qquad (\BC^{\delta_{i\bullet}}, V_\bullet)
$$
respectively, where $\{V_\bullet\}_{\bullet \in I}$ and $\{\tV_\bullet\}_{\bullet \in I}$ denote framed double quiver representations of dimension $\bv$ and $\bv+\bs^i$, with framing given by $\{W_\bullet\}_{\bullet \in I}$. Let
$$
\tW_\bullet = \BC^{\delta_{i\bullet}} \oplus W_\bullet
$$
respectively. According to Corollary \ref{cor:full attr}, the full attracting set is given by 
\begin{align}
&\attr^f_F = \Big\{ 0 \rightarrow \tV_\bullet \xrightarrow{\cong} \tV_\bullet \Big\} \label{eqn:point 1} \\
&\attr^f_{F'} = \Big\{ \BC^{\delta_{i\bullet}} \rightarrow \tV_\bullet \twoheadrightarrow V_\bullet \Big\} \label{eqn:point 2}
\end{align}
with the implication that the arrows in \eqref{eqn:point 1}--\eqref{eqn:point 2} should be compatible via the $A,B$ maps with the split short exact sequence of framing vector spaces
$$
\BC^{\delta_{i\bullet}} \hookrightarrow \BC^{\delta_{i\bullet}} \oplus W_\bullet \twoheadrightarrow W_\bullet
$$
Clearly, the attracting set of the minimal component $F$ is 
$$
\attr^f_F = \attr_F = \text{Tot}_F(q^{-1}\tCV_i^\vee)
$$
since $\Nor^-_F \CM_{\bv+\bs^i, \bw+\bs^i} = q^{-1}\tCV^\vee_i$. As for the attracting set of the component $F'$, we will show that it is a local complete intersection in the ambient space $F' \times \CM_{\bv+\bs^i, \bw+\bs^i}$.

\medskip

\begin{claim}
\label{claim:lci}

The full attracting set $\eattr^f_{F'}$, namely
$$
\{\BC^{\delta_{i\bullet}} \rightarrow \tV_\bullet \twoheadrightarrow V_\bullet\} \hookrightarrow \CM_{\bs^i,\bs^i} \times \CM_{\bv,\bw} \times \CM_{\bv+\bs^i,\bw+\bs^i}
$$
is a local complete intersection, cut out by the section $s \in \Gamma(\CM_{\bv,\bw} \times \CM_{\bv+\bs^i, \bw+\bs^i},\CE)$ of \eqref{eqn:section}, together with $2g_i$ extra equations (where $g_i$ is the number of loops at $i$).

\end{claim}

\medskip

\begin{proof} \emph{of Claim \ref{claim:lci}:} The total number of equations imposed is
$$
\text{rank } \CE + 2g_i = \bw \cdot (\bv+\bs^i) + (\bw+\bs^i)\cdot \bv - (\bv,\bv+\bs^i) + 2g_i
$$
which is the same as the codimension of the full attracting set of $F'$
\begin{multline*}
\frac 12 \left(\CM_{\bs^i,\bs^i} \times \CM_{\bv,\bw} \times \CM_{\bv+\bs^i,\bw+\bs^i} \right) = \\ = g_i + \bw\cdot \bv - \frac {(\bv,\bv)}2 + (\bw +\bs^i)\cdot (\bv+\bs^i) - \frac {(\bv+\bs^i, \bv+\bs^i)}2
\end{multline*}
Therefore, it remains to show that the full attracting set matches the locus cut out by the equations in question. From Proposition \ref{prop:lci}, we know that the zero locus of the section $s$ parameterizes maps of framed double quiver representations
\begin{equation}
\label{eqn:tv v}
\left \{ \tV_j \stackrel{\psi_j}{\twoheadrightarrow} V_j \right\}_{j \in I}
\end{equation}
Define $\BC \xrightarrow{\phi} \tV_i$ to coincide with $\tilde{A}_i|_{\BC}$, with respect to the fixed decomposition
\begin{equation}
\label{eqn:decomposition fixed}
\tW_i = \BC \oplus W_i
\end{equation}
To make $\BC$ (the domain of the map $\phi$) into a double quiver representation, we require the $2g_i$ maps $X_e,Y_e$ to act on it by the same scalars as the same named maps act on the one-dimensional kernel of the surjection \eqref{eqn:tv v}; these are the $2g_i$ extra equations referenced in the statement of Claim \ref{claim:lci}. To show that $\phi$ yields a map of framed double quiver representations (thus allowing us to combine $\phi$ and $\{\psi_j\}_{j \in I}$ into a point \eqref{eqn:point 2} of the full attracting set), it suffices to prove the identity
\begin{equation}
\label{eqn:ba coefficient}
\langle \BC|\tB_i \tA_i |\BC \rangle = 0
\end{equation}
To this end, let $\rho_i : W_i \leftrightharpoons \tW_i : \pi_i$ denote the standard inclusion and projection maps with respect to \eqref{eqn:decomposition fixed}, and let $\rho_j = \pi_j = \text{Id}_{W_j}, \forall j \neq i$. Then we have
$$
\langle \BC|\tB_i \tA_i |\BC \rangle = \text{Tr} (\tB_i \tA_i - \tB_i \tA_i \rho_i\pi_i) = \text{Tr} \sum_{j \in I} (\tB_j \tA_j - \tB_j \tA_j \rho_j\pi_j) = 
$$
$$
= \text{Tr} \sum_{j \in I} \tB_j \tA_j - \text{Tr} \sum_{j \in I} \pi_j \tB_j \tA_j \rho_j = \text{Tr} \sum_{j \in I} \tB_j \tA_j - \text{Tr} \sum_{j \in I} B_j A_j
$$
The expression above is 0 because both terms are 0, as we have seen in \eqref{eqn:trace is zero}. \end{proof}

\medskip

\noindent As a consequence of Claim \ref{claim:lci} (and of the regularity of the section $s$), we have the following equality in the equivariant $K$-theory group of $F' \times \CM_{\bv+\bs^i, \bw+\bs^i}$
\begin{equation}
\label{eqn:structure sheaf attracting}
[\CO_{\text{Attr}^f_{F'}}] = \wedge^\bullet \left(- C_{-1}^\vee + C_0^\vee - C_{1}^\vee + \sum_{e = \oii} \left(t_e + \frac q{t_e} \right) \right) = 
\end{equation}
$$
\wedge^\bullet \left( \sum_{j \in I} \left( \frac {\tW_j}{\CV_j} + \frac {q\tCV_j}{W_j} \right)  + \sum_{e = \overrightarrow{jj'} \in E} \left( \frac {t_e \tCV_j}{\CV_{j'}}  + \frac {q \tCV_{j'}}{t_e \CV_j} \right) - \left(1+q\right) \sum_{j \in I} \frac {\tCV_j}{\CV_j} \right) \prod_{e = \oii} \left[ (1-t_e) \left(1-\frac q{t_e} \right) \right]
$$
Consider the following line bundles on $F' \times \CM_{\bv+\bs^i, \bw+\bs^i}$
\begin{align*}
&\varepsilon = \sdet \left( \sum_{e = \overrightarrow{jj'} \in E} \frac {t_e\CV_j}{q^{\frac 12}(\tCV_{j'} - \CV_{j'})} - \sum_{j \in I} \frac {\CV_j}{q^{\frac 12}(\tCV_j - \CV_j)} + \sum_{j \in I} \frac {W_j}{q^{\frac 12}(\tCV_j - \CV_j)} \right) \\
&\varepsilon' = \sdet \left( \sum_{e = \overrightarrow{jj'} \in E} \frac {q^{\frac 12}\CV_{j'}}{t_e(\tCV_j - \CV_j)} - \sum_{j \in I} \frac {q^{\frac 12}\CV_j}{\tCV_j - \CV_j} + \frac {q^{\frac 12}\CV_i}{u} \right) 
\end{align*}
The line bundles above are chosen so that the classes
$$
\Gamma = \varepsilon \cdot \left[ \CO_{\text{Attr}^f_{F'}} \right] \qquad \text{and} \qquad \Gamma' = \varepsilon' \cdot \left[ \CO_{\text{Attr}^f_{F'}} \right]
$$
satisfy the axioms of the stable envelopes $\stab|_{F' \times \CM_{\bv+\bs^i,\bw+\bs^i}}$ and $\stab'|_{F' \times \CM_{\bv+\bs^i,\bw+\bs^i}}$ (this is a straightforward check of the normalization and degree conditions, which we leave to the interested reader). Therefore, we may use formula \eqref{eqn:structure sheaf attracting} to compute the restrictions of $\Gamma$ and $\Gamma'$ to $F' \times F$. We obtain
\begin{align*}
&\Gamma \Big|_{F' \times F} \stackrel{\eqref{eqn:formula simple correspondence}}= \varepsilon\Big|_{F' \times F} \cdot \wedge^\bullet \left(\frac {u}{\CV_i}  \right) \cdot \left[\CO_{\CN_{\bv,\bv+\bs^i,\bw}} \right] (1-q) \prod_{e = \oii} \left[ (1-t_e) \left(1-\frac q{t_e} \right) \right]\\
&\Gamma' \Big|_{F' \times F} \stackrel{\eqref{eqn:formula simple correspondence}}= \varepsilon'\Big|_{F' \times F} \cdot \wedge^\bullet \left(\frac {u}{\CV_i}\right) \cdot \left[\CO_{\CN_{\bv,\bv+\bs^i,\bw}} \right] (1-q) \prod_{e = \oii} \left[ (1-t_e) \left(1-\frac q{t_e} \right) \right]
\end{align*}
Since $\bCN_{\bv,\bv+\bs^i,\bw}$ is obtained from $\CN_{\bv,\bv+\bs^i,\bw}$ by imposing the $2g_i$ extra equations given by the vanishing of the maps $\{X_e,Y_e\}_{e = \oii}$ on $\BC \cong \text{Ker}(\tCV_i \twoheadrightarrow \CV_i)$, we have
$$
\left[\CO_{\bCN_{\bv,\bv+\bs^i,\bw}} \right] = \left[\CO_{\CN_{\bv,\bv+\bs^i,\bw}} \right] \prod_{e = \oii} \left[ (1-t_e) \left(1-\frac q{t_e} \right) \right]
$$
and therefore 
\begin{align*}
&\Gamma \Big|_{F' \times F} = \sdet \left( \sum_{e = \overrightarrow{ji}} \frac {t_e\CV_j}{q^{\frac 12}\CL_i} - \frac {\CV_i}{q^{\frac 12}\CL_i} + \frac {W_i}{q^{\frac 12}\CL_i} \right) \cdot (1-q) \wedge^\bullet \left(\frac {u}{\CV_i}  \right) \cdot \left[\CO_{\bCN_{\bv,\bv+\bs^i,\bw}} \right] \\
&\Gamma' \Big|_{F' \times F} = \sdet \left( \sum_{e = \overrightarrow{ij}} \frac {q^{\frac 12} \CV_j}{t_e \CL_i} - \frac {q^{\frac 12}\CV_i}{\CL_i} + \frac {q^{\frac 12}\CV_i}u\right) \cdot (1-q) \wedge^\bullet \left(\frac {u}{\CV_i}\right) \cdot \left[\CO_{\bCN_{\bv,\bv+\bs^i,\bw}} \right] 
\end{align*}
With this in mind, we conclude that the compositions
\begin{align*}
&K(F) \xrightarrow{\stab} K(\CM_{\bv+\bs^i,\bw+\bs^i}) \xrightarrow{{\stab'}^T} K(F') \\
&K(F') \xrightarrow{\stab'} K(\CM_{\bv+\bs^i,\bw+\bs^i}) \xrightarrow{\stab^T} K(F)
\end{align*}
are given by (we write points of $\bCN_{\bv,\bv+\bs^i,\bw}$ as $\{L_i \hookrightarrow \tV_i \twoheadrightarrow V_i, \tv_j = V_j\}_{j \neq i}$)
\begin{multline}
\label{eqn:formula 1}
{_{\langle \bs^i|}R(u)_{|\b0 \rangle}} = \frac {\Gamma|_{F \times F} \Gamma'|_{F' \times F}}{\wedge^\bullet(\Tan \ \CM_{\bv+\bs^i,\bw+\bs^i}|_F)^\vee} =   \\
=q^{\frac {v_i}2} (1-q)  \frac {\wedge^\bullet \left(\frac {\CV_i}{u} \right)}{\wedge^\bullet \left(\frac {q\tCV_i}{u} \right)} \cdot\sdet \left( \sum_{e = \overrightarrow{ij}} \frac {q^{\frac 12} \CV_j}{t_e \CL_i} - \frac {q^{\frac 12}\CV_i}{\CL_i} \right) \left[\CO_{\bCN_{\bv,\bv+\bs^i,\bw}} \right] 
\end{multline}
\begin{multline}
\label{eqn:formula 2}
 {_{\langle \b0 |}R(u)_{|\bs^i \rangle}} = \frac {\Gamma|_{F' \times F} \Gamma'|_{F \times F}}{\wedge^\bullet(\Tan \ \CM_{\bv+\bs^i,\bw+\bs^i}|_F)^\vee} =  \\ =q^{\frac {v_i+1}2} (1-q) \frac {\left( - \frac {\CL_i}{u} \right) \wedge^\bullet \left(\frac {\CV_i}{u} \right)}{\wedge^\bullet \left(\frac {q\tCV_i}{u} \right)}  \cdot \sdet \left( \sum_{e = \overrightarrow{ji}} \frac {t_e\CV_j}{q^{\frac 12}\CL_i} - \frac {\CV_i}{q^{\frac 12}\CL_i} + \frac {W_i}{q^{\frac 12}\CL_i} \right) \left[\CO_{\bCN_{\bv,\bv+\bs^i,\bw}} \right]
\end{multline}
respectively. In the formula above, we use the following particular cases of the normalization conditions for the fixed component $F =\CM_{\b0,\bs^i} \times \CM_{\bv+\bs^i,\bw}$ (see \eqref{eqn:alpha alpha'})
$$
\Gamma\Big|_{F \times F} = \Delta_{F*} \left(\wedge^\bullet \left(\frac u{\CV_i} \right) \right) \qquad \text{and} \qquad \Gamma'\Big|_{F \times F} = q^{\frac {v_i+1}2} \Delta_{F*} \left(\wedge^\bullet \left(\frac {\CV_i}u \right) \right) 
$$
and, using \eqref{eqn:normal bundle},
$$
\wedge^\bullet ( \Nor_F \CM_{\bv+\bs^i,\bw+\bs^i})^\vee = \wedge^\bullet \left(\frac u{\tCV_i} + \frac {q\tCV_i}u \right)
$$
By Definition \ref{def:quantum group}, the power series coefficients (in the variable $u$) of the right-hand sides of the expressions above yield elements of $\UU$, when viewed as operators
\begin{align*}
&K(\bw)_{\bv+\bs^i} \rightarrow K(\bw)_{\bv} \\
&K(\bw)_{\bv} \rightarrow K(\bw)_{\bv+\bs^i}
\end{align*}
respectively. However, up to composition with the operator of multiplication \eqref{eqn:r(a) comp 2}, which we already know is in the algebra $\UU$, the operators above are simply
\begin{align}
&\frac {1-q}{1 - \frac {q\CL_i}{u}} \cdot\sdet \left( \sum_{e = \overrightarrow{ij}} \frac {q^{\frac 12}\CV_j}{t_e \CL_i} - \frac {q^{\frac 12}\CV_i}{\CL_i} \right) \left[\CO_{\bCN_{\bv,\bv+\bs^i,\bw}} \right] \label{eqn:function of u 1} \\
&\frac {q^{-\frac 12}(1-q)}{1 - \frac u{q\CL_i}}  \cdot \sdet \left( \sum_{e = \overrightarrow{j i}} \frac {t_e\CV_j}{q^{\frac 12}\CL_i} - \frac {\CV_i}{q^{\frac 12}\CL_i} + \frac {W_i}{q^{\frac 12}\CL_i} \right) \left[\CO_{\bCN_{\bv,\bv+\bs^i,\bw}} \right] \label{eqn:function of u 2}
\end{align}
The power series coefficients in $u$ of the operators above are (up to a scalar) simply $f_{i,d}$ and $e_{i,d}$ of \eqref{eqn:action f} and \eqref{eqn:action e}, respectively, as we needed to show. \end{proof}

\medskip

\subsection{} Because $K(\bw) = \bigoplus_{\bv \in \nn} K(\bw)_{\bv}$ is $\nn$-graded, then the assignment
$$
\deg \left(K(\bw) \xrightarrow{\phi} K(\bw) \right) = \bn \quad \text{if} \quad \phi (K(\bw)_{\bv}) \subset K(\bw)_{\bv+\bn}, \forall \bv \in \nn
$$
yields a $\zz$-grading on $\text{End}(K(\bw))$, and thus also on $\UU$. We can describe this grading explicitly by observing that the $R$-matrices \eqref{eqn:r-matrix} preserve overall degrees, i.e. they act non-trivially between the graded components
\begin{equation}
\label{eqn:graded components}
K(\bw^1)_{\bv^1} \otimes K(\bw^2)_{\bv^2} \longrightarrow K(\bw^1)_{\tilde{\bv}^1} \otimes K(\bw^2)_{\tilde{\bv}^2}
\end{equation}
only if $\bv^1+\bv^2 = \tilde{\bv}^1 + \tilde{\bv}^2$. Therefore, the element
$$
K(\bw) \xrightarrow{e \otimes \text{Id}} K(\bw^{\text{aux}})_{\bm} \otimes K(\bw) \xrightarrow{R_{\bw^{\text{aux}}, \bw}} K(\bw^{\text{aux}})_{\bn} \otimes K(\bw) \xrightarrow{f \otimes \text{Id}} K(\bw)
$$
of $\UU$ has degree $\bm - \bn$, for any  $e \in K(\bw^{\text{aux}})_{\bm}$ and $f \in K(\bw^{\text{aux}})_{\bn}^\vee$.

\medskip

\begin{definition}
\label{def:completion}

Consider the completion
\begin{equation}
\label{eqn:completion 1}
\UU \woo \UU 
\end{equation}
defined with respect to the nested family of subspaces
\begin{equation}
	\label{eqn:sum completion 1}
\left\{ \UU_{\bk^1} \otimes \UU_{\bk^2} \Big| k^1_i < -N \text{ and } k^2_j > N \text{ for some }i,j \in I \right\}_{N \in \BN}
\end{equation}
Also consider the smaller completion
\begin{equation}
\label{eqn:completion 2}
\UU \too \UU 
\end{equation}
consisting of all finite linear combinations of elements in 
\begin{equation}
\label{eqn:sum completion 2}
\sum_{\bn \in \nn} \UU_{\bm^1 - \bn} \otimes \UU_{\bm^2 + \bn} 
\end{equation}
for arbitrary $\bm^1,\bm^2 \in \nn$.

\end{definition}

\medskip

\noindent Note that the completion \eqref{eqn:completion 1} has a well-defined action on $K(\bw^1) \otimes K(\bw^2)$ for all $\bw^1,\bw^2 \in \nn$, because any vector in such a tensor product is annihilated by the subspace \eqref{eqn:sum completion 1} for $N$ large enough. The following result is straightforward, and we leave it as an exercise to the reader. 

\medskip

\begin{proposition}
\label{prop:compatible}

The usual multiplication makes $\UU \too \UU$ into an algebra, and $\UU \woo \UU$ into a bimodule for the aforementioned algebra.

\end{proposition}

\medskip

\subsection{} 
\label{sub:frt principle} 

The FRT principle asserts that the quantum Yang-Baxter equation \eqref{eqn:qybe} makes $\UU$ into a Hopf algebra. In particular, there exists a universal $R$-matrix
\begin{equation}
\label{eqn:universal r-matrix}
\CR \in \UU \woo \UU 
\end{equation}
which maps to $R_{\bw^1,\bw^2}$ under the action maps
$$
\UU \woo \UU \rightarrow \text{End}(K(\bw^1) \otimes K(\bw^2))
$$
for all $\bw^1,\bw^2 \in \nn$. Moreover, we have a coproduct
\begin{equation}
\label{eqn:coproduct}
\Delta : \UU \longrightarrow \UU \too \UU
\end{equation}
(the fact that the coproduct takes values in the completion \eqref{eqn:completion 2} is a straightforward consequence of the fact that elements of $\UU$ are matrix coefficients of $R$-matrices in terms of the $\nn$-graded vector spaces $K(\bw^{\text{aux}})$, cf. Definition \ref{def:quantum group}) induced by restricting the map 
\begin{equation}
\label{eqn:coproduct conjugation}
\prod_{\bw \in \nn} \text{End}(K(\bw)) \xrightarrow{\text{conjugation by }\stab} \prod_{\bw^1,\bw^2 \in \nn} \text{End}(K(\bw^1) \otimes K(\bw^2))
\end{equation}
The coproduct and universal $R$-matrix are related by the following properties
\begin{equation}
\label{eqn:coproduct commutes with r-matrix}
\CR \Delta(x) = \Delta^{\op}(x) \CR
\end{equation}
for any $x \in \UU$, as well as the following identity in $\UU \woo \UU \woo \UU$
\begin{equation}
\label{eqn:coproduct of r-matrix}
(\text{Id} \otimes \Delta)(\CR) = \CR_{02} \CR_{01}
\end{equation}
(the three indices in $\UU \woo \UU \woo \UU$ are labeled as $0,1,2$, left to right). The following result is a $K$-theoretic counterpart of \cite[Section 6]{Mc}.

\medskip

\begin{proposition}
\label{prop:coproduct e and f}

We have the following formulas
\begin{equation}
\label{eqn:coproduct e}
\Delta(e_{i,0}) = h_{i,0}^{-1} \otimes e_{i,0} + e_{i,0} \otimes 1  
\end{equation}
\begin{equation}
\label{eqn:coproduct f}
\Delta(f_{i,0}) = f_{i,0} \otimes h_{i,0} + 1 \otimes f_{i,0} 
\end{equation}
for all $i \in I$, where $h_{i,0}$ was defined in \eqref{eqn:g to h}.

\end{proposition}

\medskip

\begin{proof} Let us consider the identity \eqref{eqn:coproduct of r-matrix} and evaluate it in the representation
\begin{equation}
\label{eqn:triple representation}
K(\bs^i) \otimes K(\bw^1) \otimes K(\bw^2)
\end{equation}
for any $\bw^1,\bw^2 \in \nn$. Then for any $\bx,\by \in \nn$, we have the following identity in terms of the matrix coefficients \eqref{eqn:r(u) coefficient}
\begin{equation}
\label{eqn:foil out}
\Delta\left(_{\langle \bx |}R(u)_{|\by \rangle}\right) = \sum_{\bz \in \nn} {_{\langle \bx |}R_2(u)_{|\bz \rangle}} {_{\langle \bz |}R_1(u)_{|\by \rangle}}
\end{equation}
with the coproduct being taken in the $K(\bw)$ tensor factor of $R(u)$ (and $R_{\epsilon}(u)$ denotes the operator \eqref{eqn:r(u)} acting in the tensor factors $K(\bs^i) \otimes K(\bw^{\epsilon})$, for any $\epsilon \in \{1,2\}$). However, by the very definition of $R(u)$, for any fixed components
$$
F = \CM_{\bx, \bs^i} \times \CM_{\bv + \bs^i - \bx, \bw} \hookrightarrow \CM_{\bv + \bs^i, \bw + \bs^i} \hookleftarrow \CM_{\by, \bs^i} \times \CM_{\bv + \bs^i - \by, \bw} = F'
$$
we have 
\begin{equation}
\label{eqn:r(u) by localization}
R(u) \Big|_{F \times F'} = \sum_{G \leq F,F'} \frac {\Gamma |_{F \times G} \Gamma' |_{F' \times G}}{\wedge^\bullet (\Tan \ X|_G)^\vee}
\end{equation}
Let $k$ and $l$ be the integers that appear in expression \eqref{eqn:normal explicit} for the fixed component $G$. Conditions \eqref{eqn:degree explicit} and \eqref{eqn:degree prime explicit} imply that the numerator of the right-hand side of \eqref{eqn:r(u) by localization} has $A \cong \BC^*$ weights contained in the interval

\medskip

\begin{itemize}

\item $(-k-l,k+l)$ if $G < F,F'$

\medskip

\item $(-k-l,k+l]$ if $G = F < F'$

\medskip

\item $[-k-l,k+l)$ if $G = F' < F$

\medskip

\item $[-k-l,k+l]$ if $G = F' = F$

\medskip

\end{itemize}

\noindent Meanwhile, the denominator of the right-hand side of \eqref{eqn:r(u) by localization} has $A \cong \BC^*$ weights contained in the interval $[-k-l,k+l]$. Therefore, we conclude that
\begin{align*}
\lim_{u \rightarrow 0} R(u) \Big|_{F \times F'} =  \begin{cases} \displaystyle \frac {\Gamma |_{F \times F'} \Gamma' |_{F' \times F'}}{\wedge^\bullet (\Tan \ X|_{F'})^\vee} &\text{if } F' \leq F \\ 0 &\text{ otherwise} \end{cases} \\
\lim_{u \rightarrow \infty} R(u) \Big|_{F \times F'} = \begin{cases} \displaystyle  \frac {\Gamma |_{F \times F} \Gamma' |_{F' \times F}}{\wedge^\bullet (\Tan \ X|_F)^\vee} &\text{if } F \leq F' \\ 0 &\text{ otherwise} \end{cases}
\end{align*}
When $F = F'$, we can explicitly compute the limits above using \eqref{eqn:normalization explicit} and \eqref{eqn:normalization prime explicit}
\begin{align}
&\lim_{u \rightarrow 0} R(u) \Big|_{F \times F} = q^{-\frac {k+l}2} \left[ \CO_{\Delta_F} \right] \label{eqn:diagonal limit 1} \\
&\lim_{u \rightarrow \infty} R(u) \Big|_{F \times F} = q^{\frac {k+l}2}  \left[ \CO_{\Delta_F} \right]\label{eqn:diagonal limit 2}
\end{align}
Let us henceforth set $F = \CM_{\b0,\bs^i} \times \CM_{\bv+\bs^i,\bw}$ and $F' = \CM_{\bs^i,\bs^i} \times \CM_{\bv, \bw}$. Then we have the following consequences of formulas \eqref{eqn:formula 2} and \eqref{eqn:formula 1}, respectively
\begin{align*}
&\lim_{u \rightarrow 0}{ _{\langle \b0 |}R(u)_{|\bs^i \rangle}}  = q^{-\frac {v_i+1}2} \left(q^{-\frac 12}-q^{\frac 12} \right) \prod_{e = \oii} \left(-\frac {q^{\frac 12}}{t_e} \right) e_{i,0} \\
&\lim_{u \rightarrow \infty} {_{\langle \bs^i |}R(u)_{|\b0 \rangle}}  = q^{\frac {v_i}2}\left(1 - q \right) f_{i,0}
\end{align*}
Moreover, comparing \eqref{eqn:diagonal limit 1}--\eqref{eqn:diagonal limit 2} with \eqref{eqn:k nakajima}--\eqref{eqn:l nakajima} shows us that
\begin{align*}
&\lim_{u \rightarrow 0} R(u) \Big|_{F \times F} = q^{-\frac {v_i+1}2}  \left[ \CO_{\Delta_F} \right] \\
&\lim_{u \rightarrow \infty} R(u) \Big|_{F \times F} = q^{\frac {v_i+1}2}  \left[ \CO_{\Delta_F}  \right]
\end{align*}
and
\begin{align*}
&\lim_{u \rightarrow 0} R(u) \Big|_{F' \times F'} = h_{i,0}^{-1} \cdot q^{-\frac {v_i}2}  \left[ \CO_{\Delta_{F'}} \right] \\
&\lim_{u \rightarrow \infty} R(u) \Big|_{F' \times F'} = h_{i,0} \cdot q^{\frac {v_i}2}  \left[ \CO_{\Delta_{F'}} \right]
\end{align*}
where in the last pair of equalities we used \eqref{eqn:multiplication leading}. With the computations above in mind, taking the limits of \eqref{eqn:foil out} for $(\bx,\by)$ given by $(\b0,\bs^i)$ and $(\bs^i,\b0)$ (respectively) as $u \rightarrow 0$ and $u\rightarrow \infty$ (respectively) yields \eqref{eqn:coproduct e} and \eqref{eqn:coproduct f} (respectively). \end{proof}

\medskip

\subsection{}

In Proposition \ref{prop:coproduct e and f}, we showed that the coproduct of $e_{i,0}, f_{i,0} \in \CA \subset \UU$ lies in $\CA \otimes \CA$. As for the other generators of the algebra $\CA$ (as per Proposition \ref{prop:generate}), because the stable envelope is additive in $\bv$ and $\bw$, we have for all $i \in I$
\begin{align}
&\Delta \left( q^{\pm \frac {v_i}2} \right) = q^{\pm \frac {v_i}2} \otimes q^{\pm \frac {v_i}2} \label{eqn:coproduct v} \\ 
&\Delta \left( q^{\pm \frac {w_i}2} \right) = q^{\pm \frac {w_i}2} \otimes q^{\pm \frac {w_i}2} \label{eqn:coproduct w} 
\end{align}
Meanwhile, the additivity of stable envelopes in the vector spaces $W_i$ implies that
\begin{equation}
\label{eqn:coproduct b}
\Delta (b_{i, \pm d}) = b_{i, \pm d} \otimes 1 + 1 \otimes b_{i, \pm d}
\end{equation}
$\forall d$. Conjecture \ref{conj:coproduct h intro} posits that similar formulas exist for the generators $a_{i, \pm 1} \in \CA$.

\medskip

\begin{proposition}
\label{prop:coproduct preserve}

Assuming Conjecture \ref{conj:coproduct h intro}, we have
\begin{equation}
\label{eqn:coproduct preserve}
\Delta(\CA) \subset \CA \too \CA
\end{equation}
In other words, the coproduct preserves the subalgebra $\CA$ of  \eqref{eqn:iota prime}.

\end{proposition}

\medskip

\begin{proof} We will start by proving the following refinement of formula \eqref{eqn:coproduct preserve}
\begin{equation}
\label{eqn:coproduct preserve plus}
\Delta(e) \in e \otimes 1 + \CA \too \CA^>
\end{equation}
\begin{equation}
\label{eqn:coproduct preserve minus}
\Delta(f) \in 1 \otimes f + \CA^< \too \CA
\end{equation}
for any $e \in \CA^+$ and $f \in \CA^-$, where $\CA^>$ and $\CA^<$ are defined in \eqref{eqn:a big and small}. According to Theorem \ref{thm:generate}, the algebra $\CA^+$ (respectively $\CA^-$) is generated by $e_{i,d}$ (respectively $f_{i,d}$) as $(i,d)$ ranges over $I \times \BZ$. Since $\CA^>$ and $\CA^<$ are closed under multiplication, in order to prove \eqref{eqn:coproduct preserve plus}--\eqref{eqn:coproduct preserve minus}, it suffices to establish the following formulas 
\begin{equation}
\label{eqn:coproduct plus details}
\Delta(e_{i,d}) \in e_{i,d} \otimes 1 + \CA \too \CA^>  
\end{equation}
\begin{equation}
\label{eqn:coproduct minus details} 
\Delta(f_{i,d}) \in 1 \otimes f_{i,d} + \CA^< \too \CA 
\end{equation}
We will prove \eqref{eqn:coproduct plus details}, and leave the analogous formula \eqref{eqn:coproduct minus details} as an exercise to the reader. According to \eqref{eqn:rel double 3}, we have the following relation for all $i \in I$ and $d \in \BZ$
$$
[e_{i,d}, a_{i,\pm 1}] = (q^{\mp 1}-1) e_{i,d \pm 1}
$$
Taking the commutator of \eqref{eqn:coproduct h pm intro} and \eqref{eqn:coproduct plus details} yields
$$
\Delta\left( (q^{\mp 1}-1) e_{i,d \pm 1} \right) = \Delta(e_{i,d})\Delta(a_{i,\pm 1}) - \Delta(a_{i,\pm 1}) \Delta(e_{i,d}) \in  \underbrace{[a_{i,\pm 1}, e_{i,d}]}_{=(q^{\mp 1}-1) e_{i,d \pm 1} } \otimes 1 + 
$$
$$
+ \underbrace{[a_{i,\pm 1}, \CA]}_{\in \CA} \too \CA^>  + \CA \too \underbrace{[a_{i,\pm 1}, \CA^>]}_{\in \CA^>} + \underbrace{[\CA^<, e_{i,d}]}_{\in \CA} \too \CA^>  + \underbrace{\CA^< \CA \too \CA^> \CA^> + \CA \CA^< \too \CA^> \CA^>}_{\in \CA \too \CA^>}
$$
(above, we use the facts that $\CA^> \CA^> \subset \CA^>$ and $[a_{i,\pm 1}, \CA^>] \subset \CA^>$). Thus, formula \eqref{eqn:coproduct plus details} for a given $d$ implies the same formula for $d \pm 1$; together with \eqref{eqn:coproduct e}, we conclude \eqref{eqn:coproduct plus details} for all $d \in \BZ$ by induction on $|d|$. Finally, it remains to prove that
\begin{equation}
\label{eqn:coproduct zero details}
\Delta(x) \in \CA \too \CA
\end{equation}
for all $x \in \CA^0$. For $x \in \{q^{\pm \frac {v_i}2}, q^{\pm \frac {w_i}2},b_{i,d}\}_{i \in I, d > 0}$, this was established in \eqref{eqn:coproduct v}, \eqref{eqn:coproduct w} and \eqref{eqn:coproduct b}. For $x \in \{h_{i,\pm d}\}_{i \in I, d > 0}$, formula \eqref{eqn:coproduct zero details} follows from \eqref{eqn:rel double 5} and the already proved relations \eqref{eqn:coproduct plus details}--\eqref{eqn:coproduct minus details}. Since we may express the generators $a_{i,\pm d}$ in terms of the aforementioned $x$'s using \eqref{eqn:formal power series} (due to the invertibility of the square matrix with coefficients \eqref{eqn:square matrix}) we conclude formula \eqref{eqn:coproduct zero details} for all $x \in \CA^0$. \end{proof}

\medskip

\subsection{}

We are now poised to complete the proof of our main Theorem.

\medskip

\begin{proof} \emph{of Theorem \ref{thm:main}:} Having already constructed an injective map 
$$
\CA \hookrightarrow \UU
$$
in Proposition \ref{prop:maps to subalgebra}, it remains to show that this map is surjective. By Proposition \ref{prop:any coefficient}, any element of $\UU$ is a linear combination of power series coefficients in the equivariant parameters $u_{i1},\dots,u_{iw^{\text{aux}}_i}$ of expressions of the form
\begin{multline}
\label{eqn:coefficients of r-matrix}
R^{e,f}_{\bw^{\text{aux}}, \bw} : \underbrace{K(\bw^{\text{aux}})_{\b0}}_{\BF_{\bw^{\text{aux}}}} \otimes K(\bw) \xrightarrow{e \otimes \text{Id}} K(\bw^{\text{aux}}) \otimes K(\bw) \\ \xrightarrow{R_{\bw^{\text{aux}}, \bw}} K(\bw^{\text{aux}}) \otimes K(\bw) \xrightarrow{f \otimes \text{Id}} \underbrace{K(\bw^{\text{aux}})_{\b0}}_{\BF_{\bw^{\text{aux}}}} \otimes K(\bw)  
\end{multline}
for various $e \in \CA_{\bm}$ and $f \in \CA_{-\bn}$. Let us first deal with the case when $e = f = 1$, i.e. $\bm = \bn = \b0$. In this case, because of the minimality of the fixed component
$$
F = \CM_{\b0,\bw^{\text{aux}}} \times \CM_{\bv,\bw} \hookrightarrow \CM_{\bv,\bw^{\text{aux}}+\bw}
$$
we have 
\begin{equation}
\label{eqn:vacuum vacuum}
_{\langle \b0|}(R_{\bw^{\text{aux}}, \bw})_{|\b0\rangle} = \text{multiplication by } \prod_{i \in I} \prod_{a = 1}^{w^{\text{aux}}_i} q^{\frac {v_i}2} \wedge^\bullet\left(\frac {(1-q)\CV_i}{u_{ia}} \right)
\end{equation}
(see \eqref{eqn:r(a) comp 2} for the case when $\bw^{\text{aux}} = \bs^i$, which is completely representative for the general case). Therefore, arbitrary power series coefficients (with respect to the $u_{ia}$'s) of the expression above act on $K(\bw)$ as multiplication by tautological classes. Since the latter lie in $\CA$, we conclude that the $e = f = 1$ case of \eqref{eqn:coefficients of r-matrix} lies in $\CA$.

\medskip

\noindent Let us now consider \eqref{eqn:coefficients of r-matrix} when $e = 1$, i.e. $\bm = 0$. We will show that $R^{1,f}_{\bw^{\text{aux}}, \bw}$ lies in $\CA$, by induction on the degree $\bn$ of $f$ (specifically, the order on $\nn$ with respect to which we perform induction is $\bn \leq \bn'$ if $\bn' - \bn \in \nn$). The base case $\bn = \b0$ was treated in the previous paragraph. For the induction step, recall from \eqref{eqn:coproduct preserve minus} that
\begin{equation}
\label{eqn:coproduct f proof}
\Delta(f) = 1 \otimes f + \sum_i f'_i \otimes f_i'' + \dots
\end{equation}
for various $f'_i \in \CA^<$, $f''_i \in \CA^{\leq} = \CA^0 \otimes \CA^-$, and where the ellipsis denotes summands whose second tensor factor is in $\CA \backslash \CA^{\leq}$ (and thus does not map into $K(\bw)_{\b0}$ under the action \eqref{eqn:action 1 intro} for any $\bw$). Because coproduct preserves degrees, $\deg f_i'' < \bn$ for all terms which appear in \eqref{eqn:coproduct f proof} (not considering the ellipsis terms). By the induction hypothesis, showing that $R^{1,f}_{\bw^{\text{aux}}, \bw}$ lies in $\CA$ is therefore equivalent to showing that the composition 
$$
K(\bw^{\text{aux}})_{\b0} \otimes K(\bw) \xrightarrow{R_{\bw^{\text{aux}}, \bw}} K(\bw^{\text{aux}}) \otimes K(\bw) \xrightarrow{\Delta^{\op}(f)} K(\bw^{\text{aux}})_{\b0} \otimes K(\bw) 
$$
lies in $\CA$ (note that the ellipsis terms in \eqref{eqn:coproduct f proof} do not contribute either to the display above, nor to the display below, for degree reasons). However, the fact that $R_{\bw^{\text{aux}}, \bw}$ intertwines $\Delta(f)$ with $\Delta^{\op}(f)$ (according to \eqref{eqn:coproduct commutes with r-matrix}) implies that the expression above is
$$
K(\bw^{\text{aux}})_{\b0} \otimes K(\bw) \xrightarrow{\Delta(f)} K(\bw^{\text{aux}}) \otimes K(\bw) \xrightarrow{R_{\bw^{\text{aux}}, \bw}} K(\bw^{\text{aux}})_{\b0} \otimes K(\bw) 
$$
Since any element in $\CA^<$ annihilates $K(\bw^{\text{aux}})_{\b0}$, the expression above is $R^{1,1}_{\bw^{\text{aux}}, \bw} \circ f$, and thus lies in $\CA$ by the $\bm=\bn=\b0$ case of \eqref{eqn:coefficients of r-matrix}, which we already proved.

\medskip

\noindent Finally, let us prove that $R^{e,f}_{\bw^{\text{aux}}, \bw}$ lies in $\CA$ by induction on the degree $\bm$ of $e$. The base case $\bm = \b0$ was treated in the previous paragraph. Recall from \eqref{eqn:coproduct preserve plus} that
\begin{equation}
\label{eqn:coproduct e proof}
\Delta(e) = e \otimes 1 + \sum_i e_i' \otimes e_i'' + \dots
\end{equation}
for various $e_i' \in \CA^{\geq} = \CA^+ \otimes \CA^0$, $e_i'' \in \CA^>$, and where the ellipsis denotes terms whose first tensor summand is in $\CA \backslash \CA^{\geq}$ (and thus annihilates $K(\bw)_{\b0}$). Because coproduct preserves degrees, $\deg e_i' < \bm$ for all terms which appear in \eqref{eqn:coproduct e proof} (not considering the ellipsis terms). By the induction hypothesis, the fact that $R^{e,f}_{\bw^{\text{aux}}, \bw} \in \CA$ is equivalent to the composition
\begin{multline*}
K(\bw^{\text{aux}})_{\b0} \otimes K(\bw) \xrightarrow{\Delta(e)} K(\bw^{\text{aux}}) \otimes K(\bw) \xrightarrow{R_{\bw^{\text{aux}}, \bw}} \\ K(\bw^{\text{aux}}) \otimes K(\bw) \xrightarrow{f \otimes \text{Id}} K(\bw^{\text{aux}})_{\b0} \otimes K(\bw)
\end{multline*}
lying in $\CA$ (note that the ellipsis terms in \eqref{eqn:coproduct e proof} do not contribute either to the display above, nor to the display below, for degree reasons). However, the fact that $R_{\bw^{\text{aux}}, \bw}$ intertwines $\Delta(e)$ with $\Delta^{\op}(e)$ (according to \eqref{eqn:coproduct commutes with r-matrix}) implies that the expression above is
$$
K(\bw^{\text{aux}})_{\b0} \otimes K(\bw) \xrightarrow{R_{\bw^{\text{aux}}, \bw}} K(\bw^{\text{aux}}) \otimes K(\bw) \xrightarrow{(f \otimes \text{Id}) \cdot \Delta^{\op}(e)} K(\bw^{\text{aux}})_{\b0} \otimes K(\bw)
$$
The expression above is equal to various elements of $\CA$ composed with expressions of the form $R^{1,f'}_{\bw^{\text{aux}}, \bw}$ for $f' \in \CA^{\leq}$, which we have already showed lie in $\CA$. \end{proof}
  
\medskip

\section{Appendix}
\label{sec:appendix}

\medskip

\subsection{} In the present Section, we will use equivariant localization and shuffle algebra computations to prove Proposition \ref{prop:inj}, i.e. the injectivity of the map
$$
\CA \xrightarrow{\iota} \prod_{\bw \in \nn} \text{End}(K(\bw))
$$
Recall the torus $T$ of \eqref{eqn:small torus} and its action on $\CM_{\bv,\bw}$ from \eqref{eqn:small torus acts}. In the proof of Proposition \ref{prop:any coefficient}, we saw that a quadruple $(X_e,Y_e,A_i,B_i) \in \CM_{\bv,\bw}$ is fixed under the aforementioned $T$-action precisely when there exists a decomposition of the underlying vector spaces
\begin{equation}
\label{eqn:fixed 0}
V_i = \mathop{\bigoplus_{\text{characters}}}_{\chi : T \rightarrow \BC^*} V_{i,\chi}
\end{equation}
such that the maps $(X_e,Y_e,A_i,B_i)$ are only non-zero between the direct summands
\begin{equation}
\label{eqn:fixed 1}
X_e : V_{i,\chi} \rightarrow V_{j,\frac {\chi}{t_e}}, \quad Y_e : V_{j,\chi} \rightarrow V_{i,\frac {\chi t_e}q}
\end{equation}
\begin{equation}
\label{eqn:fixed 2}
A_i : W_i \rightarrow V_{i,1}, \qquad \quad \ B_i : V_{i,q} \rightarrow W_i
\end{equation}
The stability condition (i.e. the fact that the vector spaces $\{V_i\}_{i\in I}$ are generated by the maps $X_e, Y_e$ acting on $\{\text{Im }A_i\}_{i\in I}$) implies that the decomposition \eqref{eqn:fixed 0} must necessarily satisfy the property that
$$
V_{i,\chi} \neq 0
$$
only if $\chi$ is a product of $t_e$ and $\frac q{t_e}$ raised to non-positive powers.

\medskip

\begin{proposition}
\label{prop:projective}

For any $\bv,\bw \in \nn$, the fixed point set $\CM_{\bv,\bw}^T$ is projective.

\end{proposition}

\medskip

\begin{proof} Recall the proper map studied in \cite{Nak 0}
$$
\CM_{\bv,\bw} \rightarrow \CM_{\bv,\bw}^0 = \text{Spec} \left(\BC[\mu_{\bv,\bw}^{-1}(0)]^{G_{\bv}} \right)
$$
It suffices to prove that $\CM_{\bv,\bw}^T$ lies in the inverse image of 0 under this map, or in other words, that any non-constant $G_{\bv}$-invariant algebraic function on the set of quadruples $(X_e,Y_e,A_i,B_i)_{e \in E, i\in I}$ vanishes on $T$-fixed quadruples. However, all such functions are polynomial expressions in the following particular types of functions

\medskip

\begin{itemize}

\item trace of any polynomial in $X_e,Y_e$, and

\medskip

\item matrix coefficients of compositions of the form $A_i (\text{polynomial in }X_e,Y_e )B_j$

\end{itemize}

\medskip

\noindent Both functions above vanish on a $T$-fixed quadruple, because with respect to the decomposition \eqref{eqn:fixed 0}, any polynomial $P$ in $X_e,Y_e$ only maps non-trivially 
$$
\text{from } V_{i,\chi} \text{ to } V_{i', \chi'}
$$ 
if $\chi' / \chi$ is a product of $t_e$ and $\frac q{t_e}$ raised to non-negative powers. If $P$ has no constant term, the assumption in Remark \ref{rem:smaller torus} precludes $\chi = \chi'$, which implies that $\text{Tr}(P) = 0$. Similarly, the fact that $A_iPB_j$ must equal $0$ is due to the fact that the maps $A_i$ and $B_j$ are non-zero only on the direct summands \eqref{eqn:fixed 2}, and the fact that the character $q$ cannot be expressed  as a product of $t_e$ and $\frac q{t_e}$ raised to non-positive powers. \end{proof}

\medskip

\subsection{} A particular connected component of $\CM_{\bv,\bw}^T$ is the one where the only non-zero vector spaces in \eqref{eqn:fixed 0} are those corresponding to $\chi = 1$. For such a quadruple, equation \eqref{eqn:fixed 1} implies that $X_e = Y_e = 0$ for all edges $e$, while  \eqref{eqn:fixed 2} implies that $B_i = 0$ for all $i \in I$. Thus, the connected component in question is none other than
\begin{equation}
\label{eqn:grassmannian}
\Gr(\bv,\bw) := \prod_{i \in I} \Gr(v_i,w_i) \stackrel{\tau}{\hookrightarrow} \CM_{\bv,\bw}
\end{equation}
where $\Gr(v,w)$ denotes the Grassmannian of $v$-dimensional quotients of a given $w$-dimensional vector space. The fact that $\Gr(\bv,\bw)$ is a component of the $T$-fixed locus of $\CM_{\bv,\bw}$ implies that the push-forward of $\tau$ induces an injection
\begin{equation}
\label{eqn:inj grassmannian 1}
K_{T \times G_{\bw}}(\Gr(\bv,\bw))_{\loc} \stackrel{\tau_*}{\hookrightarrow} K_{T \times G_{\bw}}(\CM_{\bv,\bw})_{\loc} = K(\bw)_{\bv}
\end{equation}
where we recall that ``$\loc$" refers to localization with respect to the $T$-action only, i.e. tensoring with $\BQ(q^{\frac 12},t_e)_{e \in E}$. Let us now consider a maximal torus
$$
T_{\bw} \subset G_{\bw}
$$
and we will assume that $\{u_{i1},\dots,u_{iw_i}\}_{i \in I}$ denote the elementary characters of $T_{\bw}$. With this in mind, we have for any variety $X$ endowed with a $T \times G_{\bw}$ action
\begin{equation}
\label{eqn:weyl}
K_{T \times G_{\bw}}(X)_{\loc} \cong K_{T \times T_{\bw}}(X)_{\loc}^{\sym} \hookrightarrow K_{T \times T_{\bw}}(X)_{\loc}
\end{equation}
Explicitly, $K_{T \times T_{\bw}}(X)_{\loc}$ is a $S(\bw) = \prod_{i \in I} S(w_i)$ equivariant module over $\BF[u^{\pm 1}_{ia}]_{i \in I}^{a \leq w_i}$. Then \eqref{eqn:weyl} simply states that $K_{T \times G_{\bw}}(X)_{\loc}$ is the $S(\bw)$ invariant part of this action, which is a module over 
$$
\BF_{\bw} = \BF[u^{\pm 1}_{i1},\dots,u^{\pm 1}_{iw_i}]_{i \in I}^{\sym}
$$
Our main reason for working with $T_{\bw}$-equivariant instead of $G_{\bw}$-equivariant $K$-theory is that we may define the classes of fixed points
\begin{equation}
\label{eqn:fixed point}
I_{\bS} = \left[ \CO_{\CI_{\bS}} \right] \in K_{T \times T_{\bw}}(\CM_{\bv,\bw})_{\loc}
\end{equation}
for any collection of subsets
$$
\bS = (S_i)_{i \in I} \qquad \text{with} \quad S_i \subset \{1,\dots,w_i\}, \ |S_i| = v_i, \forall i \in I
$$
Explicitly, $\CI_S$ denotes the $T \times T_{\bw}$ fixed point of $\CM_{\bv,\bw}$ in which $X_e = Y_e = 0$ for all edges $e$, $B_i = 0$ for all vertices $i \in I$, while $A_i : W_i \rightarrow V_i$ annihilates the $a$-th coordinate vector for all $a \in \{1,\dots,w_i\}\backslash S_i$ (on the other hand, the stability condition forces the $a$-th coordinate vectors of $W_i$ for $a \in S_i$ to be mapped to a basis of $V_i$). Since $\CI_S$ are also the standard fixed points of the Grassmannian $\Gr(\bv,\bw)$, the injectivity of the map $\tau_*$ implies that the classes \eqref{eqn:fixed point} are linearly independent. 

\medskip

\subsection{} 

We will now compute the action of the operators $f_{i,d}$ of \eqref{eqn:action f} on the classes \eqref{eqn:fixed point}. Recall the diagram
\begin{equation}
\label{eqn:ext diagram 3}
\xymatrix{& \bar{\CN}_{\bv-\bs^i,\bv,\bw} \ar[ld]_{\bar{\pi}_+} \ar[rd]^{\bar{\pi}_-} & \\
\CM_{\bv,\bw} & & \CM_{\bv-\bs^i,\bw}}
\end{equation}
with respect to which 
\begin{equation}
\label{eqn:action f appendix}
f_{i,d} = \bar{\pi}_{-*} \left( \CL_i^d \cdot \sdet \left[ \sum_{e = \oij} \frac {q^{\frac 12}\CV^-_j}{t_e \CL_i} - \frac {q^{\frac 12} \CV_i^-}{\CL_i} \right] \cdot \bar{\pi}_+^* \right) 
\end{equation}
Above, we write points of $\bar{\CN}_{\bv-\bs^i,\bv,\bw}$ as short exact sequences
\begin{equation}
\label{eqn:ses 3}
\Big( 0\rightarrow \BC^{\delta_{i\bullet}} \rightarrow V^+_\bullet \rightarrow V^-_\bullet \rightarrow 0 \Big)
\end{equation}
of framed double quiver representations. If $V^+_\bullet$ corresponds to a fixed point $\CI_{\bS} \in \CM_{\bv,\bw}$ (i.e. $X_e = Y_e = 0$ for all $e$, $B_i = 0$ for all $i$, while $A_i$ is non-zero only on those coordinate vectors which correspond to indices in $S_i$) then the choice of long exact sequence in \eqref{eqn:ses 3} boils down to the choice of a line in $V^+_i$. Thus, as a set
$$
\bar{\pi}_+^{-1}(\CI_{\bS}) = \BP\left( (V^+_i)^\vee \right)
$$
The torus fixed points in the projective space above are $\CI_{\bS - a_i}$, where for any $a \in S_i$ 
$$
\bS - a_i := (T_j)_{j \in I}, \quad \text{where} \quad T_j = \begin{cases} S_j &\text{if } j \neq i \\ S_j \backslash a &\text{if } j = i \end{cases}
$$
However, as a scheme, the map $\bar{\pi}_+$ is a local complete intersection morphism, explicitly described by the projectivization of the two-step complex \eqref{eqn:projectivization 2}. Let $I_{\bS,\bS - a_i}$ denote the skyscraper sheaf at the torus fixed point $(\CI_{\bS}, \CI_{\bS-a_i}) \in \bar{\CN}_{\bv-\bs^i, \bv,\bw}$. Then a straightforward fixed point computation reveals that
$$
\bar{\pi}_+^{*}(I_{\bS}) = \sum_{a \in S_i} \frac {I_{\bS,\bS-a_i}}{\wedge^\bullet \left( \text{Tan}_{\bar{\pi}_+} |_{(\CI_{\bS},\CI_{\bS-a_i})} \right)^\vee} = \sum_{a \in S_i} I_{\bS,\bS-a_i} \cdot \wedge^\bullet \left(1 + \frac {q \CL_i}{\CU_i} \Big|_{(\CI_{\bS},\CI_{\bS-a_i})} \right)
$$
where the last equality is due to the usual (virtual) tangent bundle computation of a projectivization. To compute the right-hand side of the preceding equation, we use \eqref{eqn:universal complex} to obtain the following formula in $K$-theory
$$
[\CU_i] = W_i + \sum_{e = \oij} \frac {q}{t_e} [\CV_j] + \sum_{e = \oji} t_e [\CV_j] - (1+q) \cdot[\CV_i]
$$
and the obvious restriction formula
$$
\CV_j |_{\CI_{\bS}} = \sum_{b \in S_j} u_{jb}
$$
Therefore, we obtain
$$
\bar{\pi}_+^{*}(I_{\bS}) = \sum_{a \in S_i} I_{\bS,\bS-a_i} \cdot \wedge^\bullet \left(-q+\sum_{b=1}^{w_i} \frac {qu_{ia}}{u_{ib}} + \sum^{e = \oij}_{b \in S_j} \frac {t_e u_{ia}}{u_{jb}} + \sum^{e = \oji}_{b \in S_j} \frac {qu_{ia}}{t_eu_{jb}} - \sum_{b \in S_i \backslash a} \frac {(1+q)u_{ia}}{u_{ib}} \right) 
$$
Since $\bar{\pi}_{-*} (I_{\bS,\bS-a_i}) = I_{\bS-a_i}$, formula \eqref{eqn:action f appendix} implies
\begin{equation}
\label{eqn:single f in fixed points}
f_{i,d}(I_{\bS}) = \sum_{\bT = \bS - a_i} I_{\bT} \cdot  u_{ia}^d
\end{equation}
$$
\frac {\sigma_i \prod^{e = \oij}_{b \in T_j} q^{\frac 12} \left(1 - \frac {u_{jb}}{u_{ia}t_e} \right) \prod^{e = \oji}_{b \in T_j} \left(1 - \frac {qu_{ia}}{t_eu_{jb}}\right) \prod_{b \in \{1,\dots,w_i\} \backslash S_i} \left(1 - \frac {qu_{ia}}{u_{ib}} \right)}{\prod_{b \in T_i} q^{\frac 12} \left(1 - \frac {u_{ib}}{u_{ia}} \right)}
$$
where $\sigma_i = \prod_{e = \oii} (1-t_e)\left(1 - \frac q{t_e}\right)$ and we write $\bT = (T_j)_{j \in I}$ above. 

\medskip

\subsection{} We may iterate formula \eqref{eqn:single f in fixed points} and infer the following expression
\begin{equation}
\label{eqn:many f in fixed points}
f_{i_1,d_1} \dots f_{i_n,d_n} (I_{\bS}) = \sum_{\bT = \bS - (a_1)_{i_1} - \dots - (a_n)_{i_n}} I_{\bT} \cdot \prod_{k=1}^n u_{i_ka_k}^{d_k} \prod_{1 \leq k < l \leq n} \zeta_{i_l i_k} \left(\frac {u_{i_la_l}}{u_{i_ka_k}} \right)
\end{equation}
$$
\prod_{k=1}^n \frac {\sigma_{i_k} \prod^{e = \overrightarrow{i_kj}}_{b \in T_j} q^{\frac 12} \left(1 - \frac {u_{jb}}{u_{i_ka_k}t_e} \right) \prod^{e = \overrightarrow{ji_k}}_{b \in T_j} \left(1 - \frac {qu_{i_ka_k}}{t_eu_{jb}}\right) \prod_{b \in \{1,\dots,w_{i_k}\} \backslash S_{i_k}} \left(1 - \frac {qu_{i_ka_k}}{u_{i_kb}} \right)}{\prod_{b \in T_{i_k}} q^{\frac 12}\left(1 - \frac {u_{i_kb}}{u_{i_ka_k}} \right)}
$$
for all $i_1,\dots,i_n \in I$ and $d_1,\dots,d_n \in \BZ$, where
\begin{equation}
\label{eqn:def zeta}
\zeta_{ij}(x) = \left( \frac {x-q}{q^{\frac 12}(x-1)} \right)^{\delta_{ij}}\prod_{e = \oij} q^{\frac 12} \left(1 - \frac 1{xt_e} \right) \displaystyle \prod_{e = \oji} \left(1 - \frac {qx}{t_e} \right)
\end{equation}
To develop a systematic treatment of operators of the form \eqref{eqn:many f in fixed points}, we recall the concept of shuffle algebras.

\medskip

\begin{definition}
\label{def:shuffle}

The big shuffle algebra is defined as
\begin{equation}
\label{eqn:def shuf}
\CV = \bigoplus_{\bn = (n_i)_{i \in I} \in \nn} \BF[z^{\pm 1}_{i1},\dots,z^{\pm 1}_{in_i}]_{i \in I}^{\emph{sym}}
\end{equation}
endowed with the multiplication 
\begin{equation}
\label{eqn:shuf mult}
R(z_{i1},\dots,z_{in_i})_{i \in I} * R'(z_{i1},\dots,z_{in_i'})_{i \in I} = 
\end{equation}
$$
 = \emph{Sym} \left[\frac {R(z_{i1},\dots,z_{in_i})_{i \in I} R'(z_{i,n_i+1},\dots,z_{i,n_i+n_i'})_{i \in I}}{\bn! \bn'!} \prod_{i,j \in I} \prod_{a=1}^{n_i} \prod_{b=n_j+1}^{n_j+n_j'} \zeta_{ji} \left(\frac {z_{jb}}{z_{ia}} \right) \right]
$$
where $\bn! = \prod_{i \in I} n_i!$. The shuffle algebra 
$$
\CS^- \subset \CV
$$
is defined as the subalgebra generated by $\{z_{i1}^d\}_{i \in I}^{d \in \BZ}$. 

\end{definition}

\medskip

\noindent It was shown in \cite{Wheel} that
\begin{equation}
\label{eqn:wheel}
R (\dots,z_{ia},\dots) \in \CS^- \qquad \Leftrightarrow \qquad R\Big|_{z_{ia} = \frac {qz_{jb}}{t_e} = qz_{ic}} = R\Big|_{z_{ja} = t_e z_{ib} = qz_{jc}} = 0
\end{equation}
for all $i,j \in I$, all edges $e = \oij$ and all indices $a,b,c$ such that $a \neq c$ (if $i = j$ then we also require $a \neq b \neq c$ in the right-hand side of \eqref{eqn:wheel}). One of the main results of \cite{Wheel} is that we have an isomorphism
\begin{equation}
\label{eqn:iso shuffle}
\CA^- \xrightarrow{\sim} \CS^-
\end{equation}
$$
f_{i_1,d_1} \dots f_{i_n,d_n} \mapsto \Sym \left[ z_{i_1a_1}^{d_1} \dots z_{i_na_n}^{d_n} \prod_{1 \leq k < l \leq n} \zeta_{i_li_k} \left( \frac {z_{i_la_l}}{z_{i_ka_k}} \right) \right]
$$
where for all $i \in I$, the multiset $\{a_k \text{ s.t. } i_k = i\}$ is required to be some permutation of $\{1,2,\dots\}$. Thus, the action \eqref{eqn:action minus} yields via the isomorphism \eqref{eqn:iso shuffle} an action
\begin{equation}
\label{eqn:shuffle acts}
\CS^- \curvearrowright K(\bw)
\end{equation}
Due to the particular form of the isomorphism \eqref{eqn:iso shuffle}, it is clear that formula \eqref{eqn:many f in fixed points} translates into the following formula for any $R (z_{i1},\dots,z_{in_i})_{i \in I} \in \CS^-$ 
\begin{equation}
\label{eqn:shuffle in fixed points}
R(I_{\bS}) = \sum^{\bT \subset \bS}_{|T_i| = S_i - n_i, \forall i \in I} I_{\bT} \cdot \frac {R(\dots,u_{ia},\dots)_{i \in I, a \in S_i \backslash T_i}}{\bn!}
\end{equation}
$$
\prod^{i\in I}_{a \in S_i \backslash T_i} \frac {\sigma_{i} \prod^{e = \overrightarrow{ij}}_{b \in T_j} q^{\frac 12} \left(1 - \frac {u_{jb}}{u_{ia}t_e} \right) \prod^{e = \overrightarrow{ji}}_{b \in T_j} \left(1 - \frac {qu_{ia}}{t_eu_{jb}}\right) \prod_{b \in \{1,\dots,w_i\} \backslash S_i} \left(1 - \frac {qu_{ia}}{u_{ib}} \right)}{\prod_{b \in T_i} q^{\frac 12}\left(1 - \frac {u_{ib}}{u_{ia}} \right)}
$$
Thus, up to a product of predictable linear factors, elements $R \in \CS^-$ keep track of matrix coefficients of the corresponding endomorphisms of $K(\bw)$ in the linearly independent set $\{I_{\bS}\}_{\bS = (S_i)_{i \in I}}$.

\medskip

\subsection{} 

Using formula \eqref{eqn:shuffle in fixed points}, we are almost ready to give a proof of Proposition \ref{prop:inj}. The last ingredient we need is the following slight modification of the Euler characteristic pairing of Subsection \ref{sub:notation 4}
\begin{equation}
\label{eqn:modified pairing}
K_{T \times G_{\bw}}(\CM_{\bv,\bw})_{\loc} \otimes K_{T \times G_{\bw}}(\CM_{\bv,\bw})_{\loc} \xrightarrow{(\cdot,\cdot)'} \BF_{\bw}
\end{equation}
$$
(\alpha,\beta)' = \chi_{T \times G_{\bw}} \left(\CM_{\bv,\bw},  \alpha \cdot \beta \cdot \sdet \left[ \sum_{e = \oij \in E} \frac {q^{\frac 12} \CV_j}{t_e \CV_i}  - \sum_{i \in I} \frac {q^{\frac 12} \CV_i}{\CV_i} + \sum_{i \in I} \frac {q^{\frac 12} \CV_i}{W_i} \right] \right)
$$
as well as the analogous pairing for $T_{\bw}$ instead of $G_{\bw}$ equivariant $K$-theory (with the caveat that in the latter case the pairing takes values in $\BF[u_{i1}^{\pm 1},\dots,u^{\pm 1}_{iw_i}]_{i \in I}$). 

\medskip

\begin{proposition}
\label{prop:adjoint}

To any $\Gamma \in \CA_{\bn} = \CA_{-\bn}$ (the latter equality is one of $\BF$-vector spaces), formulas \eqref{eqn:action plus} and \eqref{eqn:action minus} assign the following endomorphisms of $K(\bw)$
\begin{align}
&e_\Gamma = \pi_{+*} \left(\emph{sdet} \left[ \sum_{e = \oji \in E} \frac {t_e\CV_j^+}{q^{\frac 12}\CK_i} - \sum_{i\in I} \frac {\CV_i^+}{q^{\frac 12} \CK_i} + \sum_{i\in I}\frac {W_i}{q^{\frac 12}\CK_i} \right] \cdot p^!(\Gamma) \cdot \pi^*_- \right) \label{eqn:e gamma} \\
&f_\Gamma = \pi_{-*} \left(\emph{sdet} \left[ \sum_{e = \oij \in E} \frac {q^{\frac 12}\CV^-_j}{t_e \CK_i} - \sum_{i\in I} \frac {q^{\frac 12}\CV_i^-}{\CK_i} \right] \cdot p^!(\Gamma) \cdot \pi^*_+ \right) \label{eqn:f gamma}
\end{align}
Then $e_{\Gamma}$ and $f_{\Gamma}$ are adjoint with respect to the modified pairing \eqref{eqn:modified pairing}
\begin{equation}
\label{eqn:modified adjointness}
\Big(e_\Gamma(\alpha), \beta \Big)' = \Big(\alpha, f_\Gamma(\beta) \Big)'
\end{equation}
for any $\alpha \in K_{T\times G_{\bw}}(\CM_{\bv,\bw})_{\emph{loc}}$ and $\beta \in K_{T\times G_{\bw}}(\CM_{\bv+\bn,\bw})_{\emph{loc}}$.

\end{proposition}

\medskip

\noindent The natural analogue of Proposition \ref{prop:adjoint} holds with $G_{\bw}$-equivariant replaced by $T_{\bw}$-equivariant $K$-theory. We leave the proof of either ($G_{\bw}$ or $T_{\bw}$-equivariant) case of Proposition \ref{prop:adjoint} as an exercise to the interested reader. It is predicated on the fact that the operators $\pi_{\pm *}$ and $\pi_\pm^*$ are adjoint with respect to the usual Euler characteristic pairing of Subsection \ref{sub:notation 4}. The need for modifying this pairing by the line bundle that features in \eqref{eqn:modified pairing} is due to the different line bundles that appear in formulas \eqref{eqn:e gamma} and \eqref{eqn:f gamma}.

\medskip

\begin{proof} \emph{of Proposition \ref{prop:inj}:} Let us consider an arbitrary non-zero element
\begin{equation}
\label{eqn:non-zero phi}
\phi = \sum_{\kappa} e^{\kappa} P_{\kappa} f_{\kappa} \in \CA^+ \otimes \CA^0 \otimes \CA^- = \CA
\end{equation}
for various $e^\kappa \in \CA^+$, $P_{\kappa} \in \CA^0$, $f_{\kappa} \in \CA^-$. Assume for the purpose of contradiction that $\phi$ acts by the 0 endomorphism on $K(\bw)$ for all $\bw \in \nn$. Because the endomorphism associated to $\phi$ is given by a correspondence (i.e. a $K$-theory class) in
$$
K_{T \times G_{\bw}}\left(\bigsqcup_{\bv^1,\bv^2 \in \nn} \CM_{\bv^1,\bw} \times \CM_{\bv^2,\bw} \right)_{\loc}
$$
then by mapping this correspondence to $T \times T_{\bw}$ equivariant $K$-theory, we would conclude that $\phi$ also acts by 0 on 
\begin{equation}
\label{eqn:tw vector space}
\bigoplus_{\bv \in \nn} K_{T \times T_{\bw}} \left( \CM_{\bv,\bw} \right)_{\loc}
\end{equation}
Our reason for replacing $G_{\bw}$ by $T_{\bw}$ equivariant $K$-theory above is that we may now use the linearly independent elements $I_{\bS}$ for various $\bS = (S_i)_{i \in I}$, $S_i \subset \{1,\dots,w_i\}$. More specifically, for any such $I$-tuples of sets $\bS,\bS'$, we conclude that
$$
0 = \Big(\phi(I_{\bS}), I_{\bS'} \Big)' = \sum_{\kappa} \Big(e^{\kappa} P_{\kappa} f_{\kappa}(I_{\bS}), I_{\bS'} \Big)'
$$
for the modified pairing \eqref{eqn:modified pairing}. If we let $f^\kappa$ denote the adjoint operator of $e^\kappa$ (i.e. the same element viewed under the equality of vector spaces $\CA^+ = \CA^-$, as shown in Proposition \ref{prop:adjoint}), then we have
$$
0 = \sum_{\kappa} \Big(P_{\kappa} f_{\kappa}(I_{\bS}), f^\kappa(I_{\bS'}) \Big)'
$$
However, we may now use formula \eqref{eqn:shuffle in fixed points} applied to the elements $R_\kappa, R^\kappa \in \CS^-$ that correspond to $f_\kappa, f^\kappa$ under the isomorphism \eqref{eqn:iso shuffle}
$$
0 = \sum_\kappa \mathop{\sum_{\bT \subset \bS}}_{\bT' \subset \bS'} \Big(P_\kappa (I_{\bT}), I_{\bT'} \Big)' R_\kappa(\dots, u_{ia}, \dots)_{i \in I, a \in S_i \backslash T_i}  R^\kappa(\dots, u_{ia}, \dots)_{i \in I, a \in S'_i \backslash T'_i}  \Omega 
$$
where $\Omega$ is a non-zero product of linear factors that depends on $\bS, \bS', \bT, \bT'$ but not on $\kappa$. Finally, let us note that any element $P_\kappa \in \CA^0$ is a symmetric Laurent polynomial in the Chern classes of $W_i$ and $\CV_i$, and therefore
$$
P_{\kappa}(I_{\bT}) = \rho_\kappa(\underbrace{\dots, u_{ia}, \dots}_{a \in \{1,\dots,w_i\}} | \underbrace{\dots, u_{ia}, \dots}_{a \in T_i}) \cdot I_{\bT} 
$$
where the Laurent polynomial $\rho_\kappa(\dots, x_i, \dots | \dots y_j \dots)$ is symmetric in the variables $x_i$ and $y_j$ separately. Finally, we note that the classes $I_{\bT}$ are orthogonal under the pairing \eqref{eqn:modified pairing} (as they are skyscraper sheaves of different torus fixed points) and that the pairing of $I_{\bT}$ with itself is given by a non-zero product of linear factors (which is none other than the exterior algebra of the cotangent space to $\CI_{\bT}$ in $\CM_{\bv,\bw}$. We thus conclude that
\begin{multline}
0 = \sum_{\bS \supset \bT \subset \bS'} \sum_\kappa R_\kappa(\dots, u_{ia}, \dots)_{i \in I, a \in S_i \backslash T_i}  R^\kappa(\dots, u_{ia}, \dots)_{i \in I, a \in S'_i \backslash T_i} \cdot  \\ \cdot \rho_\kappa(\underbrace{\dots, u_{ia}, \dots}_{a \in \{1,\dots,w_i\}} | \underbrace{\dots, u_{ia}, \dots}_{a \in T_i})  \Omega'  \label{eqn:multline above}
\end{multline}
where $\Omega'$ is another product of linear factors that depends on $\bS,\bS',\bT$ but not on $\kappa$. We may obtain a contradiction from \eqref{eqn:multline above} as follows. Let us choose $\bm,\bn \in \nn$ such that there exist non-zero terms in \eqref{eqn:non-zero phi} with $\deg e^\kappa = \bm$ and $\deg f_\kappa = -\bn$, but $\bm$ and $\bn$ are maximal in $\nn$ with respect to this property. Then let us choose $\bp,\bw \in \nn$ such that $w_i \gg p_i+n_i+m_i$ for all $i \in I$. If 
\begin{align*}
&\bS = (S_i)_{i \in I} \quad \text{with } \quad S_i = \{1,\dots,p_i,p_i+1,\dots,p_i+n_i\} \\
&\bS' = (S'_i)_{i \in I} \quad \text{with} \quad S'_i = \{1,\dots,p_i,p_i+n_i+1,\dots,p_i+n_i+m_i\}
\end{align*}
then the only $\bT$ which can appear in the right-hand side of \eqref{eqn:multline above} corresponds to $T_i = \{1,\dots,p_i\}$ for all $i \in I$, and we have
\begin{multline*}
0 = \sum^{\kappa \text{ corresponding to}}_{\deg e^\kappa = \bm \text{ and } \deg f_\kappa = -\bn} R_{\kappa}(u_{i,p_i+1},\dots,u_{i,p_i+n_i})_{i \in I} \cdot \\ R^\kappa(u_{i,p_i+n_i+1},\dots,u_{i,p_i+n_i+m_i})_{i \in I} \cdot \rho_\kappa(u_{i1},\dots,u_{iw_i}|u_{i1},\dots,u_{ip_i})_{i \in I}
\end{multline*}
However, the elements $R_\kappa, R^\kappa \in \CS^-$ and the polynomials $\rho_\kappa$ all have coefficients in the field $\BF = \BQ(q^{\frac 12},t_e)_{e \in E}$. The fact that the identity above holds for all formal parameters $u_{i1},\dots,u_{iw_i}$ (assuming $\bp$ is chosen large enough compared to the finitely many polynomials $\rho_\kappa$ which appear) implies that 
$$
0 = \sum^{\kappa \text{ corresponding to}}_{\deg e^\kappa = \bm \text{ and } \deg f_\kappa = -\bn} R^\kappa \otimes p_\kappa \otimes R_\kappa
$$
identically, which contradicts the fact that the element $\phi$ is non-zero. \end{proof}

\end{document}